\documentclass[a4paper,english,10pt,twoside,dvips]{amsart}
\advance\oddsidemargin by -0.1cm
\advance\evensidemargin by -0.7cm
\advance\topmargin by -1.0cm %.6 initially
\usepackage{amssymb} 
\usepackage{babel}
\usepackage{amstext}
\usepackage{diagrams}
\usepackage{amscd}   
\usepackage{epsfig}
\usepackage{url}  
\usepackage{rotating}
\usepackage{psfrag}  
\usepackage{multicol}
\usepackage{array}
\def\haken{\mathbin{\hbox to 6pt{%
                 \vrule height0.4pt width5pt depth0pt
                 \kern-.4pt
                 \vrule height6pt width0.4pt depth0pt\hss}}}
\theoremstyle{plain}
\newtheorem{cor}{Corollary}[section]
\newtheorem{lem}{Lemma}[section]
\newtheorem{thm}{Theorem}[section]
\newtheorem{prop}{Proposition}[section]
\theoremstyle{definition}
\newtheorem{exa}{Example}[section]
\newtheorem{NB}{Remark}[section]
\newtheorem{dfn}{Definition}[section]
\newtheorem*{NB*}{Remark}
\newcommand{\be}{\begin{equation}}
\newcommand{\ee}{\end{equation}}
\newcommand{\bdm}{\begin{displaymath}}
\newcommand{\edm}{\end{displaymath}}
\newcommand{\ba}[1]{\begin{array}{#1}}
\newcommand{\ea}{\end{array}}
\newcommand{\bea}[1][]{\begin{eqnarray#1}}
\newcommand{\eea}[1][]{\end{eqnarray#1}}
\newcommand{\bcen}{\begin{center}}
\newcommand{\ecen}{\end{center}}
\newcommand{\btab}{\begin{tabular}}
\newcommand{\etab}{\end{tabular}}

\newcommand{\x}{\times}
\newcommand{\op}{\oplus}
\newcommand{\ox}{\otimes}

\newcommand{\ra}{\rightarrow}
\newcommand{\lra}{\longrightarrow}

\newcommand{\qqs}{\forall}

\newcommand{\Id}{\ensuremath{\mathrm{Id}}}

\newcommand{\ad}{\ensuremath{\mathrm{ad}}}

\newcommand{\Ad}{\ensuremath{\mathrm{Ad}}}

\newcommand{\del}{\partial}

\renewcommand{\Im}{\ensuremath{\mathrm{Im\,}}}
\renewcommand{\Re}{\ensuremath{\mathrm{Re\,}}}

\newcommand{\diag}{\ensuremath{\mathrm{diag}}}
\newcommand{\grad}{\ensuremath{\mathrm{grad}}}

\newcommand{\cyclic}[1]{\stackrel{{\scriptsize #1}}{\mathfrak{S}}}

\newcommand{\lan}{\left\langle}
\newcommand{\ran}{\right\rangle}

\newcommand{\ovn}{\overline{\nabla}}

\newcommand{\kr}{\ensuremath{\mathcal{R}}}
\newcommand{\Ric}{{\mathrm{Ric}}}

\newcommand{\rank}{{\mathrm{rank}}}
\newcommand{\scal}{{\mathrm{Scal}}}
\newcommand{\Scal}{{\mathrm{scal}}}

\newcommand{\Jacm}{\ensuremath{\mathrm{Jac}_{\mathfrak{m}}}}
\newcommand{\Jach}{\ensuremath{\mathrm{Jac}_{\mathfrak{h}}}}

\newcommand{\arctanh}{\ensuremath{\mathrm{arctanh}\,}}
\newcommand{\Adtilde}{\ensuremath{\widetilde{\mathrm{A}\,}\!\mathrm{d}\,}}
\newcommand{\adtilde}{\ensuremath{\widetilde{\mathrm{a}\,}\!\mathrm{d}\,}}
\newcommand{\C}{\ensuremath{\mathbb{C}}}
\newcommand{\R}{\ensuremath{\mathbb{R}}}

\newcommand{\Z}{\ensuremath{\mathbb{Z}}}

\renewcommand{\P}{\ensuremath{\mathbb{P}}}

\newcommand{\G}{\ensuremath{\mathcal{G}}}
\newcommand{\g}{\ensuremath{\mathfrak{g}}}

\newcommand{\V}{\ensuremath{\mathcal{V}}}               % vector fields
\newcommand{\W}{\ensuremath{\mathcal{W}}}

\newcommand{\X}{\ensuremath{\mathcal{X}}}

\newcommand{\h}{\ensuremath{\mathfrak{h}}}
\newcommand{\m}{\ensuremath{\mathfrak{m}}}
\newcommand{\n}{\ensuremath{\mathfrak{n}}}

\newcommand{\D}{\ensuremath{\mathcal{D}}}
\newcommand{\T}{\ensuremath{\mathrm{T}}}

\newcommand{\eps}{\ensuremath{\varepsilon}}             % nice epsilon
\newcommand{\vphi}{\ensuremath{\varphi}}                % nice phi
\newcommand{\vrho}{\ensuremath{\varrho}}                % nice rho

\newcommand{\End}{\ensuremath{\mathrm{End}}}
\newcommand{\Hol}{\ensuremath{\mathrm{Hol}}}
\newcommand{\hol}{\ensuremath{\mathfrak{hol}}}
\newcommand{\GL}{\ensuremath{\mathrm{GL}}}

\newcommand{\SL}{\ensuremath{\mathrm{SL}}}

\newcommand{\un}{\ensuremath{\mathfrak{u}}}
\newcommand{\su}{\ensuremath{\mathfrak{su}}}
\newcommand{\SU}{\ensuremath{\mathrm{SU}}}
\newcommand{\U}{\ensuremath{\mathrm{U}}}
\newcommand{\Sympl}{\ensuremath{\mathrm{Sp}}}

\newcommand{\so}{\ensuremath{\mathfrak{so}}}
\newcommand{\SO}{\ensuremath{\mathrm{SO}}}
\newcommand{\Orth}{\ensuremath{\mathrm{O}}}
\newcommand{\Spin}{\ensuremath{\mathrm{Spin}}}
\newcommand{\spin}{\ensuremath{\mathfrak{spin}}}

\begin{document}
\thispagestyle{empty}

\hbox to \hsize{%
  \vtop{} \hfill
  \vtop{\hbox{\small to appear in Rend.~del Circ.~Mat.~di Palermo, 2006}}}
\title[The Srni lectures on non-integrable geometries with torsion]
{The Srni lectures on non-integrable geometries with torsion}
%\title{The Srn\'{\i} lectures on non-integrable  geometries with torsion}
%----------------------------------------------------------------
\author[Ilka Agricola]{Ilka Agricola}
%\author{\vspace{6mm}\small Version: \today}
%
\address{\hspace{-5mm}
{\normalfont\ttfamily agricola@mathematik.hu-berlin.de}\newline
Institut f\"ur Mathematik \newline
Humboldt-Universit\"at zu Berlin\newline
Unter den Linden 6\newline
Sitz: John-von-Neumann-Haus, Adlershof\newline
D-10099 Berlin, Germany}
\thanks{Supported by the SFB 647 "Space---Time---Matter" of the DFG and
the Junior Research Group "Special Geometries in Mathematical Physics" of
the Volkswagen Foundation.}
\begin{abstract}
%---------------
This review article intends to introduce the reader to non-integrable 
geometric structures on Riemannian manifolds and
invariant metric connections with torsion, and to discuss
recent aspects of mathematical
physics---in particular superstring theory---where these naturally appear. 

Connections with skew-symmetric torsion are exhibited as one of the main 
tools to understand  non-integrable geometries.  To this aim a 
a series of key examples is presented and successively dealt with
using the notions of  intrinsic
torsion and characteristic connection of a $G$-structure as unifying
principles.
The General Holonomy Principle bridges over to parallel objects, thus 
motivating the discussion of geometric stabilizers, with emphasis on 
spinors and differential forms. Several Weitzenb\"ock formulas for
Dirac operators associated with  torsion connections enable us 
to discuss spinorial field equations, such as those governing the common
sector of type  II superstring theory. They also provide the link to
Kostant's cubic Dirac operator.

\end{abstract}
\maketitle
\tableofcontents
\pagestyle{headings}
%
%------------- body of the document -------------------------------------------
%
%------------------------------------------------------------------------------
%
% Summary
%
% I. Agricola  
%
%----------------------------------------------------------------------------- 
%
%---------------------------------------------------------------------------
\section{Background and motivation}
%---------------------------------------------------------------------------
%  
\subsection{Introduction} 
%-----------------------------------------------------------------------------
%
Since Paul Dirac's  formulation in 1928 of the field equation for a quantized 
electron in flat Minkowski space, Dirac operators on Riemannian manifolds 
have become  a powerful tool for the treatment of various problems in
geometry, analysis and theoretical physics. Meanwhile, starting from the 
fifties
the French school founded by M.\,Berger had developed  the idea that manifolds
should be subdivided into different classes according to their holonomy group.
The name \emph{special (integrable) geometries} has become customary for those
which are not of general type. Already at  that early stage
there were hints that parallel spinor fields would induce special geometries,
but this idea was not further investigated. At the beginning of the
seventies, A.\,Gray  generalized the classical holonomy concept by
introducing a classification principle for \emph{non-integrable
special Riemannian geometries} 
%(\cite{Gray71}) 
and  studied
the defining differential equations of each class. The connection between
these two lines of research in mathematical physics became clear
in the eighties  in the context of twistor theory and the study of small
eigenvalues of the Dirac operator,  mainly developed by
the Berlin school around Th.\,Friedrich. 
%(\cite{Friedrich80},\cite{BFGK}) 
In the case of homogeneous manifolds, integrable geometries correspond to 
symmetric spaces, whose classification by E.\,Cartan is a milestone in 20th 
century differential geometry. The much richer class of
homogeneous reductive 
spaces---which is inaccessible to any kind of classification---has been
studied intensively since the mid-sixties, and is a main
source of examples for non-integrable geometries.

\smallskip
The interest in non-integrable geometries was revived in the past years
through developments of superstring theory. Firstly, integrable
geometries (Calabi-Yau manifolds, Joyce manifolds etc.) are exact
solutions of the Strominger model (1986), though
with vanishing $B$-field. If one deforms these vacuum equations and looks for
models with non-trivial $B$-field, a new mathematical approach going back 
a decade implies that solutions can be constructed geometrically  from non
integrable geometries with torsion. In this way, manifolds \emph{not} 
belonging to the field
of algebraic geometry (integrable geometries) become 
candidates for interesting models in theoretical physics.

\smallskip
Before discussing the deep mathematical and physical backgrounds,
let us give a---very intuitive---explanation of why the traditional
Yang-Mills approach needs modification in string theory and how torsion
enters the scene.
Point particles move along world-lines, and physical quantities are typically
computed as line integrals of some potential that is, mathematically 
speaking, just a $1$-form. The associated field strength is then its
differential---a $2$-form---and  interpreted as the curvature
of some connection. In contrast, excitations of extended $1$-dimensional
objects (the `strings') are `world-surfaces', and  physical quantities 
have to be
surface  integrals of certain potential $2$-forms. Their field strengths
are thus  $3$-forms and cannot be interpreted as curvatures anymore.
The key idea is to supply the  (pseudo)-Riemannian manifold underlying 
the physical model with a non-integrable $G$-structure admitting a `good' 
metric $G$-connection $\nabla$ with torsion, which in turn will
play the role of a $B$-field strength; and the art is to choose
the $G$-structure so that the connection  $\nabla$ admits the desired
parallel objects, in particular  spinors,  interpreted as
supersymmetry transformations.

\smallskip
My warmest thanks go to all colleagues who, by their countless remarks and
corrections,  helped improving the quality of this text, in particular
to Simon Chiossi, Richard Cleyton, Thomas Friedrich,  Mario Kassuba,
Nils Schoemann (Humboldt University Berlin) as well as Pawe{\l} 
Nurowski and Andrzej Trautman (Warsaw University).
\subsection{Mathematical motivation}\label{mathmot}
%---------------------------------------------------------------------------
From classical mechanics, it is a well-known fact that symmetry
considerations can simplify the study of geometric problems---for
example, Noether's theorem tells us how to construct first integrals,
like momentum, from invariance properties of the underlying mechanical
system.  In fact, beginning from 1870, it became clear that the principle organizing geometry ought 
to be the study of its symmetry groups. In his inaugural lecture at the 
University of Erlangen in 1872, which later became known as the ``Erlanger 
Programm'', Felix Klein said \cite[p. 34]{Klein1872}: 
\begin{quote}
\emph{Let a manifold 
and on it a transformation group be given; 
the objects belonging to the manifold ought to be studied with respect 
to those properties which are not changed by the transformations of the 
group}\footnote{ \glqq Es ist eine Mannigfaltigkeit und in derselben 
eine Transformationsgruppe
gegeben; man soll die der Mannigfaltigkeit angeh\"origen Gebilde
hinsichtlich solcher Eigenschaften untersuchen, die durch die
Transformationen der Gruppe nicht ge\"andert werden.\grqq}.
\end{quote}
Hence the classical symmetry approach in differential geometry was based
on the \emph{isometry group} of a manifold, that is, the group of all
transformations acting on the given manifold. 

By the mid-fifties, a second intrinsic group
associated to a Riemannian manifold turned out to be deeply related to
fundamental features like  curvature and parallel objects. 
The so-called 
\emph{holonomy group} determines how a vector can change under 
parallel transport along a closed loop inside the manifold (only
in the flat case  the transported vector will coincide with the
original one). Berger's theorem (1955) classifies all possible 
restricted holonomy groups of a simply connected, irreducible and non-symmetric
Riemannian manifold $(M,g)$ (see \cite{Berger55}, \cite{Simon62} for
corrections and simplifications in the proof and \cite{Bryant96} for a
status report). The holonomy group  can be either $\SO(n)$ in the
generic case or one of the groups listed in Table~\ref{Berger-list}
(here and in the sequel, $\nabla^g$ denotes the Levi-Civita connection).
\begin{table}
\bdm
\setlength{\extrarowheight}{4pt}
\begin{array}{|c|c|c|c|c|c|c|c|}
\hline
\dim M & 4n & 2n & 2n & 4n & 7 & 8 & [16]\\[2mm] \hline\hline
\Hol(M) & \Sympl(n)\Sympl(1) & \U(n) & \SU(n) & \Sympl(n) &
G_2 & \Spin(7) & [\Spin(9)]\\[1mm] \hline
\text{\small name} & 
\ba{c}\text{\small quatern.-} \\ \text{\small K\"ahler m.}\ea & 
\ba{c}\text{\small K\"ahler}\\ \text{\small manifold}\ea & 
\ba{c}\text{\small Calabi-}\\ \text{\small Yau m.}\ea & 
\ba{c} \text{\small hyper-}\\  \text{\small K\"ahler m.}\ea  & 
\ba{c} \text{\small parallel }\\  \text{\small $G_2$-m.}\ea  &
\ba{c} \text{\small parallel }\\  \text{\small $\Spin(7)$-m.}\ea  &
\ba{c}\text{\small [parallel }\\ \text{\small $\Spin(9)$-m.}]\ea\\[1mm]\hline
\text{\small par.\,object} & - & \nabla^g J = 0 & \nabla^g J = 0 & \nabla^g J = 0 & 
\nabla^g \omega^3 =0 & \nabla^g\beta^4 =0& - \\[1mm]\hline
\text{\small curvature} & \Ric = \lambda g & - & \Ric = 0 & \Ric = 0 & \Ric = 0 & 
\Ric = 0 & - \\ \hline
\end{array}
\edm
\bigskip
\caption{Possible Riemannian holonomy groups (`Berger's list').}
\label{Berger-list}
\end{table}
Manifolds having one of these holonomy groups are 
called manifolds with special (integrable) holonomy, or 
\emph{special (integrable) geometries} for short. 
We put the case $n=16$ and $\Hol(M)=\Spin(9)$ into brackets,
because Alekseevski and Brown/Gray showed independently that 
such a  manifold is necessarily symmetric (\cite{Alekseevski68},
\cite{Brown&G72}). The point is that Berger proved that the
groups on this list were the only possibilities, but he was not able to show 
whether they actually occurred as holonomy groups of compact manifolds. It 
took another thirty years to find out that---with the exception of
$\Spin(9)$---this is indeed the case: The existence
of metrics with holonomy $\SU(m)$ or $\Sympl(m)$ on compact manifolds 
followed from
Yau's solution of the Calabi Conjecture posed in 1954 \cite{Yau78}.
Explicit non-compact metrics with holonomy $G_2$ or $\Spin(7)$
are due to R.~Bryant \cite{Bryant87} and R.~Bryant and S.~Salamon
\cite{Bryant&S89}, while  compact manifolds with
holonomy $G_2$ or $\Spin(7)$ were constructed by D.\ Joyce only in 1996
(see \cite{Joyce96a}, \cite{Joyce96b} \cite{Joyce96c}
and the book \cite{Joyce}, which also contains a proof of the
Calabi Conjecture). Later, compact exceptional holonomy manifolds have also 
been constructed by other  methods by Kovalev (\cite{Kovalev03}).

As we will explain later, the General Holonomy Principle relates manifolds 
with $\Hol(M)=\SU(n), \Sympl(n), G_2$ or $\Spin(7)$
with $\nabla^g$-parallel spinors (see Section \ref{geom-stab}).
Already in the sixties it  had been observed 
that the existence of such a spinor implies in turn the vanishing of the 
Ricci curvature (\cite{Bonan66} and Proposition \ref{curv-parallel}) and 
restricts the holonomy group of the manifold 
(\cite{Hitchin74}, \cite{McKW89}),
but the difficulties in constructing explicit compact manifolds with
special integrable Ricci-flat metrics inhibited further research on the
deeper meaning of this result. 

There was progress in this direction
only in the homogeneous case. Symmetric spaces
are the "integrable" geometries inside the much larger class of
homogeneous reductive spaces. Given a non-compact semisimple Lie 
group $G$ and a maximal compact subgroup $K$ such that $\rank\, G=\rank\, K$,
consider the associated symmetric space $G/K$. The Dirac operator can  be
twisted by a finite-dimensional irreducible unitary
representation $\tau$ of $K$, and it was shown by Parthasarathy,
Wolf, Atiyah and Schmid that for suitable $\tau$
most of the discrete series representations of $G$ can be
realized on the $L^2$-kernel of this twisted Dirac operator 
(\cite{Parthasarathy72}, \cite{Wolf74}, \cite{Atiyah&S77}). The crucial
step  is to relate the square of the Dirac operator with
the Casimir operator $\Omega_G$ of $G$; for trivial $\tau$, the corresponding
formula reads
\be\label{Kos-Parth-D2-symm}
D^2\ =\ \Omega_G + \frac{1}{8}\,\scal\,.
\ee
Meanwhile many people began looking for suitable generalizations of the 
classical holonomy concept. One motivation for this was that the
notion of Riemannian holonomy  is too restrictive for vast classes of
interesting Riemannian manifolds; for example, contact  or
almost Hermitian manifolds cannot be distinguished merely by their holonomy
properties (they have generic holonomy $\SO(n)$), and the Levi-Civita
connection is not adapted to the underlying geometric structure (meaning
that the defining objects are not parallel).

In 1971 A.~Gray introduced the notion of \emph{weak holonomy}
(\cite{Gray71}), "one of his most original concepts"  and "an idea much 
ahead of its time" in the words of N.\ Hitchin \cite{Hitchin01}. This concept
turned out to yield interesting non-integrable geometries in 
dimensions  $n \leq 8$ and $n=16$. 
In particular, manifolds with weak holonomy $\U(n)$ and $G_2$ became known
as \emph{nearly K\"ahler} and \emph{nearly parallel $G_2$}-manifolds,
respectively. But whereas  metrics of  compact Ricci-flat integrable
geometries have not been realized explicitly (so far), there are many 
well-known
homogeneous reductive examples of non-integrable geometries (\cite{Gray70}, 
\cite{Fernandez87}, \cite{BFGK}, \cite{FKMS}, \cite{BG}, \cite{Fino03} and
many others). The relation to Dirac operators emerged shortly after
Th.~Friedrich proved in 1980  a seminal inequality
for the first eigenvalue $\lambda_1$ of the Dirac operator on a compact 
Riemannian manifold $M^n$ of non-negative curvature \cite{Friedrich80}
 \be\label{friedrich-80}
  (\lambda_1)^2\ \geq\ \frac{n}{4 (n-1)}\,\min_{M^n} (\Scal)\,, 
 \ee
In this estimate,
equality occurs  precisely if the corresponding eigenspinor $\psi$
satisfies the Killing equation 
 \bdm
 \nabla^g_X\psi \ =\ \pm \frac{1}{2}\sqrt{\frac{\min(\Scal)}{n(n-1)}} 
X\cdot \psi \ =: \ \mu X\cdot\psi.
 \edm
The first non-trivial compact examples of Riemannian manifolds with Killing 
spinors were found in dimensions $5$ and $6$ in 1980 and 1985, respectively 
(\cite{Friedrich80}, \cite{Friedrich&G85}). The link to 
non-integrable geometry was established shortly after; 
for instance, a compact, connected and simply connected $6$-dimensional 
Hermitian manifold is nearly K\"ahler if and only if it admits 
a Killing spinor with real Killing number $\mu$ \cite{Grunewald90}.
Similar results hold for Einstein-Sasaki structures in dimension $5$ and
nearly parallel $G_2$-manifolds in dimension $7$ (\cite{FK89}, \cite{FK}).
Remarkably, the proof of  inequality (\ref{friedrich-80}) relies
on introducing a suitable spin connection---an idea much in line with 
recent developments.
A.\ Lichnerowicz  established the link to  twistor theory
by showing  that on a compact manifold the space
of twistor spinors coincides---up to a conformal change of the
metric---with the space of  Killing spinors \cite{Lichnerowicz88}.

%-----------------------------------------------------------------------------
\subsection{Physical motivation -- torsion in gravity}\label{physmotgrav}
%-----------------------------------------------------------------------------
%
The first attempts to introduce torsion as an additional 'datum' for 
describing physics in general relativity goes back to Cartan himself 
\cite{Cartan24a}. Viewing torsion as some intrinsic angular momentum,
 he derived a set of gravitational field equations from a variational
principle, but postulated that the energy-momentum tensor should still
be divergence-free, a condition too  restrictive
for making this approach useful. The idea was taken up again in 
broader context in the late fifties. The variation of the scalar curvature
(and an additional Lagrangian generating the energy-momentum tensor) on a 
space-time endowed with a metric connection with torsion
yielded the two fundamental equations of \emph{Einstein-Cartan theory},
first formulated by Kibble \cite{Kibble61} and Sciama (see his article 
in \cite{Infeld62}). The first equation is (formally) Einstein's classical 
field equation of general relativity with an effective energy 
momentum tensor $T_{\mathrm{eff}}$ depending on torsion, the second one  
can be written in index-free notation as
\bdm
Q(X,Y)+ \sum_{i=1}^n \big(Q(Y,e_i)\haken e_i\big)\cdot X - 
 \big(Q(X,e_i)\haken e_i\big)\cdot Y \ =\ 8\pi S(X,Y).
\edm
Here, $Q$ denotes the torsion of the new connection $\nabla$, $S$ the spin 
density and $e_1,\ldots,e_n$ any orthonormal frame. A.~Trautman provided an 
elegant formulation of  Einstein-Cartan theory in the
language of principal fibre bundles \cite{Trautman73a}. 
The most striking predictions of Einstein-Cartan theory are in cosmology.
In the presence of very dense spinning matter, nonsingular cosmological 
models may be constructed because the effective energy 
momentum tensor $T_{\mathrm{eff}}$ does not fulfill the conditions of the 
Penrose-Hawking singularity theorems anymore \cite{Trautman73b}. The first 
example of such a model was provided by W.~Kopczy\'nski \cite{Kopczynski73}, 
while J.~Tafel found  a large class of such
models with homogeneous spacial sections \cite{Tafel75}. For a general 
review of gravity with spin and torsion including extensive references, we 
refer to the article \cite{Hehl&H&K&N76}.  

In the absence of spin, the torsion vanishes and the whole theory reduces
to Einstein's original formulation of general relativity. In practice,
torsion turned out to be hard to detect experimentally, since
all  tests of general relativity are based on experiments in empty space.  
Einstein-Cartan theory is pursued no longer, although some concepts that 
it inspired are still of relevance  (see 
\cite{Hehl&M&M&N95} for a generalization with additional currents and shear,
\cite{Trautman99} for optical aspects, \cite{Ruggiero&T03} for the link to 
the classical theory of defects in elastic media). Yet, it may be possible 
that Einstein-Cartan theory will prove to be a better classical  limit 
of a future quantum theory of gravitation than the theory without torsion.

%----------------------------------------------------------------------------
\subsection{Physical motivation -- torsion in superstring 
theory}\label{physmot}
%----------------------------------------------------------------------------
%
Superstring theory (see for example \cite{Green&S&W87}, \cite{Luest&T89}) 
is a physical theory 
aiming at describing nature at small distances ($\simeq~10^{-25}$~m). 
The concept of point-like elementary particles is 
replaced by one-dimensional objects as building blocks of matter---the
so-called strings. Particles are then understood as resonance states of 
strings and can be described together with their interactions up to
very high energies (small distances) without internal contradictions.
Besides gravitation, string theory incorporates many other gauge interactions
and hence is an excellent candidate for a more profound description
of matter than the standard model of elementary particles.
Quantization of superstrings is only possible in the critical
dimension $10$, while $M$-theory is a non-perturbative description of
superstrings with "geometrized" coupling, and lives in dimension $11$.
Mathematically speaking,   a   $10$- or $11$-dimensional 
configuration space $Y$ (a priori not necessarily smooth) 
is assumed to be the product
\bdm
Y^{10,11}\ =\ V^{3\text{-}5}\times M^{5\text{-}8}
\edm
of a low-dimensional spacetime $V$ describing the `external' part
of the theory (typically, Minkowski space or a space motivated from
general relativity like anti-de-Sitter space), and a higher-dimen\-sion\-al
`internal space' $M$ with some special geometric structure.
 The metric is typically a direct or a warped product.
On $M$, internal symmetries of particles are described by parallel spinor 
fields, the most important of which being the existing
supersymmetries: a spinor field has spin $1/2$, so  tensoring with it
swaps bosons and  fermions. By the General Holonomy Principle
(see Theorem \ref{hol-principle}), the holonomy group has to be a
subgroup of the  stabilizer  of the 
set of parallel spinors inside $\Spin(9,1)$. These  are well known and 
summarized in Table \ref{par-spin-dim10}.
We shall explain how to derive this result and how to understand the occurring
semidirect products in Section \ref{stab-spin7}.

\begin{table}
\begin{center}
\setlength{\extrarowheight}{4pt}
\begin{tabular}{|c|c|}\hline
\# of inv. spinors  & stabilizer groups \\ \hline\hline
1 & $\Spin(7)\ltimes \R^8$\\ \hline
2& $G_2$, $\SU(4)\ltimes\R^8$ \\ \hline
3 & $\Sympl(2)\ltimes\R^8$ \\ \hline
4 & $\SU(3)$, $(\SU(2)\times\SU(2))\ltimes\R^8$  \\ \hline
8 & $\SU(2)$, $\R^8$ \\ \hline
16 & $\{e\}$ \\ \hline
\end{tabular}
\bigskip
\end{center}
\caption{Possible stabilizers of invariant spinors inside $\Spin(9,1)$.}
\label{par-spin-dim10}
\end{table}

Since its early days, string theory has been intricately related with 
some branches of algebraic geometry. This is due to the fact that
the integrable, Ricci-flat geometries with a parallel spinor field 
with respect to the Levi-Civita connection are  exact solutions
of the Strominger model for a string vacuum with vanishing $B$-field
and constant dilaton.
This rich and active area of mathematical research lead to
interesting developments such as the discovery of mirror symmetry.

%
%----------------------------------------------------------------------------
\subsection{First developments since 1980}\label{since-1980}
%----------------------------------------------------------------------------
%
In the early eighties, several physicists independently tried to
incorporate torsion into superstring and supergravity theories in order
to get a  more physically flexible model, possibly inspired by the 
developments in classical gravity (\cite{VanNieuwenhuizen81}, 
\cite{Gates&H&R84}, \cite{Howe&P87}, \cite{deWit&S&HD87}, \cite{Rocek92}). 
In fact, simple supergravity is  equivalent to Einstein-Cartan theory with
a massless, anticommuting  Rarita-Schwinger field as source.  But
contrary to general relativity, one difficulty  
stems from the fact that there are several models in superstring theory (type
I, II, heterotic\ldots) that vary in the excitation spectrum and the 
possible interactions. 

In his article  "Superstrings with torsion" \cite{Strominger86}, 
A.~Strominger  describes the basic model in the common sector of type II 
superstring theory as a $6$-tuple $(M^n,g,\nabla,T,\Phi, \Psi)$ 
consisting of a Riemannian spin manifold $(M^n,g)$, a $3$-form $T$, a dilaton 
function $\Phi$ and a spinor field $\Psi$. The field equations can be 
written in the following form (recall that $\nabla^g$ denotes
the Levi-Civita connection):
\bdm
\Ric_{ij} - \frac{1}{4}T_{imn}T_{jmn} + 2 \nabla^g_i\partial_j \Phi
\ =\ 0,\quad \delta(e^{-2\Phi}T)\ =\ 0, \ 
\edm
\bdm
(\nabla^{g}_X +\frac{1}{4}X\haken T)\psi\ =\ 0,\quad 
(2 d\Phi-T)\cdot\psi\ =\ 0\,.
\edm
If one introduces a new metric
connection $\nabla$ whose torsion is given by the $3$-form $T$,
\bdm
\nabla_X Y\ :=\  \nabla^{g}_X Y + \frac{1}{2} T(X,Y,-),
\edm
one sees that the third equation is equivalent to $\nabla\Psi=0$.
The remaining  equations can similarly be rewritten in terms of $\nabla$. For 
constant dilaton $\Phi$, they take the particularly simple form
\cite{Ivanov&P01}
\be\label{string-eq}
\Ric^{\nabla} \ = \ 0, \quad  \delta^g(T) \ = \ 0, 
\quad \nabla \Psi \ = \ 0, \quad T \cdot \Psi \ = \ 0 \, ,
\ee
and the second equation ($\delta^g(T) = 0$) now follows from the first
 equation ($\Ric^{\nabla}=0$).
For $M$ compact, it was shown in \cite[Theorem 4.1]{Agri} that a solution
of all equations  necessarily forces $T=0$, i.\,e.\ an integrable
Ricci-flat geometry with classical holonomy given by Berger's list. 
By a careful analysis of the integrability conditions, this result
could later be extended to the non-compact case (\cite{Agri&F&N&P05},
see also Section \ref{super}).
Together with the well-understood Calabi-Yau manifolds, Joyce manifolds with 
Riemannian holonomy $G_2$ or $\Spin(7)$ thus became  of interest in recent 
times (see \cite{Atiyah&W01}, \cite{Curio&K&L01}). From a mathematical
point of view, this result stresses the importance of tackling easier
problems first, for example partial solutions. As first
step in the investigation of metric connections 
with totally skew-symmetric torsion, Dirac operators, parallel 
spinors etc., Th.~Friedrich and S.~Ivanov proved that many non-integrable 
geometric structures (almost contact metric structures, nearly K\"ahler and 
weak $G_2$-structures) admit a unique invariant connection $\nabla$ 
with totally skew-symmetric torsion \cite{Friedrich&I1}, thus being a
natural replacement for the Levi-Civita connection.
Non-integrable geometries could then be studied by their 
 holonomy properties. 

\smallskip
In fact, in mathematics the times had ben  ripe
for a new look at the intricate relationship between holonomy, special 
geometries, spinors and differential forms: in 1987, R.~Bryant found the first 
explicit local examples of metrics with exceptional Riemannian holonomy
(see \cite{Bryant87} and \cite{Bryant&S89}), Chr.~B\"ar
described their relation to Killing spinors via the cone
construction \cite{Baer93}. Building on the insightful vision of 
Gray, S.~Salamon realized the centrality of the concept of 
\emph{intrinsic torsion} (\cite{Salamon89} and, for recent results,
\cite{Fino98}, \cite{Chiossi&S02}, \cite{Cleyton&S04}).
Swann successfully tried \emph{weakening holonomy} \cite{Swann00}, and
N.\,Hitchin  characterized non-integrable geometries  as critical points 
of some linear functionals
on differential forms \cite{Hitchin01}. In particular, he motivated
a  generalization of Calabi-Yau-manifolds \cite{Hitchin00}
and of $G_2$-manifolds \cite{Witt04}, and discovered 
a new, previously unknown special geometry in dimension $8$
("weak $\mathrm{PSU}(3)$-structures", see also \cite{Witt06}).
Friedrich reformulated the concepts of non-integrable geometries in
terms of principle fiber bundles \cite{Fri2} and discussed the
exceptional dimension $16$ suggested by A.~Gray  years before
(\cite{Friedrich01}, \cite{Friedrich02a}). Analytic problems---in particular,
the investigation of the Dirac operator---on
non-integrable Riemannian manifolds contributed to a further
understanding of the underlying geometry (\cite{Bis}, \cite{Alexandrov&I00},
\cite{Gauduchon97}). Finally,  the Italian school and  collaborators
devoted over the past years a lot of effort to the explicit 
construction of
homogeneous examples of non-integrable geometries with special
properties in small dimensions (see for example \cite{Abbena&G&S00},
\cite{Fino&G03}, \cite{Fino&P&S04}, \cite{Salamon01} 
and the literature cited therein), making it possible to test the 
different concepts on explicit examples.

\smallskip
The first non-integrable geometry that raised the interest of 
string theorists was the squashed $7$-sphere with its weak $G_2$-structure,
although the first steps in this direction were still marked by confusion
about the different holonomy concepts.
A good overview about $G_2$ in string theory is the survey article 
by M.\,Duff (\cite{Duf02}). It  includes speculations about possible
applications of weak $\Spin(9)$-structures in dimension $16$, which 
a priori are of too high dimension to be 
considered in physics. In dimension three, it is well known
(see for example \cite{Spindel&S&T&vP88}) that the Strominger
equation $\nabla\Psi=0$ can  be solved only on a compact Lie
group with biinvariant metric, and that the torsion of the 
invariant connection $\nabla$ coincides with the Lie bracket. In
dimension four, the Strominger model leads to a HKT structure (see Section 
\ref{exa-HKT} for more references), i.\,e.~a
hyper-Hermitian structure that is parallel with respect to $\nabla$, 
and---in the compact case---the 
manifold is either a Calabi-Yau manifold or a Hopf surface \cite{Ivanov&P01}.
Hence, the first interesting dimension for further mathematical investigations
is five. 

 Obviously, besides the basic correspondence
outlined here, there is still much more going on between 
special geometries and detailed properties of physical models constructed 
from them. Some weak geometries have been rederived 
by physicists looking for partial solutions by  numerical analysis of 
ODE's and heavy special function machinery \cite{Gauntlett01}.

\smallskip
As an example of the many interesting mathematical problems appearing in
the context of string theory, the physicists Ramond and Pengpan
observed  that there is an infinite set of irreducible 
representations of $\Spin(9)$ partitioned into triplets
$\mathcal{S}=\cup_i \{\mu^i,\sigma^i,\tau^i \}$, whose
representations are related in a remarkable way.
For example, the infinitesimal character value of the Casimir operator
is constant on triplets, and $\dim \mu^i + \dim\sigma^i=\dim\tau^i$
if numbered appropriately. These triplets are used to describe massless 
supermultiplets, for example  $N=2$ hypermultiplets in $(3+1)$ dimensions 
with helicity $\U(1)$ or $N=1$ supermultiplets in eleven dimensions, where 
$\SO(9)$ is the light-cone little group \cite{Brink&R&X02}. 
To explain this fact, B.~Kostant introduced 
an element in the tensor product of the Clifford algebra and the 
universal enveloping algebra of a Lie group called "Kostant's cubic 
Dirac operator", and derived a striking formula
for its square (\cite{Gross&K&R&S98}, \cite{Kostant99}). The
triplet structure of the representations observed for  $\Spin(9)$
is due to the fact that the Euler characteristic of $F_4/\Spin(9)$ is 
three, hence the name "Euler multiplets" has become common for describing
this effect. In Section \ref{kostant}, we will show that Kostant's
operator may be interpreted as the symbol of a usual Dirac operator which is 
induced by a
non-standard connection on a homogeneous naturally reductive space 
(\cite{Agricola02},\cite{Agri}).
In particular, this Dirac operator satisfies a  remarkably simple 
formula  which is a direct generalization of 
Parthasarathy's formula on symmetric spaces \cite{Parthasarathy72}.
This established the link between Kostant's algebraic considerations
and recent models in string theory; in particular, it made 
homogeneous naturally reductive spaces to key examples for string theory
 and allowed  the derivation of strong vanishing theorems on them.
In representation theory this opened
the possibility to realize infinite-dimensional representations in kernels
of twisted Dirac operators on homogeneous spaces (\cite{Huang&P02}, 
\cite{Mehdi&Z04}), as it had 
been carried out on symmetric spaces in the seventies (\cite{Parthasarathy72},
\cite{Wolf74}, \cite{Atiyah&S77}). 

%-----------------------------------------------------------------------------
\section{Metric connections with torsion}
%-----------------------------------------------------------------------------
%  
%-----------------------------------------------------------------------------
\subsection{Types of connections and their lift into the spinor bundle}
%-----------------------------------------------------------------------------
%  
Let us begin with some general remarks on torsion. The notion of torsion of 
a connection was invented by Elie Cartan,
and appeared for the first time in a short note at the Acad\'emie des Sciences
de Paris in 1922 \cite{Cartan1}. Although the article contains no 
formulas, Cartan observed that such a connection may or may not preserve 
geodesics, and initially turns his attention to those who  do so. 
In this sense, Cartan was the first to investigate this class of 
connections. At that time, it was not  customary---as it became 
in the second half of the 20th century---to assign to a Riemannian manifold
only its Levi-Civita connection. Rather, Cartan demands (see \cite{Cartan4}):
\begin{quote}  
\emph{Given a manifold embedded in affine (or projective or conformal etc.)
space, attribute to this manifold the affine (or projective or conformal etc.)
connection that reflects in the simplest possible way the relations of this
manifold with the ambient space}\footnote{ \flqq \'Etant donn\'e une 
vari\'et\'e plong\'ee dans l'espace affine (ou
projectif, ou conforme etc.), attribuer \`a cette vari\'et\'e la
connexion affine (ou projective, ou conforme etc.) qui rende le plus
simplement compte des relations de cette vari\'et\'e avec l'espace 
ambiant.\frqq}.
\end{quote}
He then goes on to explain in very general terms how the connection should
be \emph{adapted} to the geometry under consideration. We believe that  this 
point of
view should be taken into account in Riemannian geometry, too. 

We now give a short review of the $8$ classes of geometric torsion tensors. 
Consider a Riemannian manifold $(M^n, g)$. The difference between its 
Levi-Civita connection $\nabla^g$ and any linear 
connection $\nabla$ is  a $(2,1)$-tensor field $A$,
\bdm 
\nabla_X Y\ =\ \nabla^g_X Y + A(X,Y),\quad X,Y \in TM^n.
\edm
The vanishing of the symmetric or the antisymmetric part of $A$ 
has immediate geometric interpretations.
The connection $\nabla$ is torsion-free if and only if $A$ is symmetric.
The connection $\nabla$ has the same geodesics as the Levi-Civita connection 
$\nabla^g$ if and only if $A$ is skew-symmetric.
Following Cartan, we study the algebraic 
types of the torsion tensor for a metric connection.
Denote by the same symbol the $(3,0)$-tensor
derived from a $(2,1)$-tensor  by contraction with the metric.
We identify $TM^n$ with $(TM^n)^*$ using $g$ from now on. 
Let $\mathcal{T}$ be the $n^2(n-1)/2$-dimensional space of all 
possible torsion tensors,
\bdm
\mathcal{T}\ =\ \{T\in\ox^3 TM^n \ | \ T(X,Y,Z)= - \, T(Y,X,Z) \}
\ \cong \ \Lambda^2 TM^n\ox TM^n \, .
\edm
A connection $\nabla$ is metric if and only if $A$ belongs to the space
\bdm
\mathcal{A}^g\ :=\ TM^n\ox(\Lambda^2 TM^n) \ = \ \{A \in\ox^3 TM^n \ | 
\ A(X,V,W) +  A(X,W,V) \ =\ 0\} \, .
\edm 
In particular, $\dim \mathcal{A}^g = \dim \mathcal{T}$, reflecting
the fact that metric connections can be uniquely characterized by their
torsion.
\begin{prop}[{\cite[p.51]{Cartan25a}, \cite{Tricerri&V1}}, \cite{Salamon89}]
\label{classessum}
%---------------------------------------------------------------------------
The spaces $\mathcal{T}$ and $\mathcal{A}^g$ are isomorphic
as $\Orth(n)$ representations, an  equivariant bijection being 
\begin{eqnarray*} 
T(X,Y,Z) & = & A(X,Y,Z)-A(Y,X,Z)\, ,\\ 
2 \, A(X,Y,Z) & = & T(X,Y,Z)-T(Y,Z,X)+T(Z,X,Y) \, .
\end{eqnarray*} 
For $n\geq 3$,  they split under the action of $\Orth(n)$
into the sum of three irreducible representations,
\bdm
\mathcal{T}\cong TM^n \op \Lambda^3(M^n) \op \mathcal{T}'.
\edm
The last module  will also be denoted $\mathcal{A}'$ if viewed as
a subspace of  $\mathcal{A}^g$ and is equivalent to the
Cartan product  of representations  $TM^n\ox \Lambda^2 TM^n$,
\bdm
\mathcal{T}'\ =\ \{ T\in\mathcal{T}\ |\  \cyclic{X,Y,Z}T(X,Y,Z)=0,\ 
\sum_{i=1}^nT(e_i,e_i,X)=0\ \qqs X,Y,Z \}
\edm
for any orthonormal frame $e_1,\ldots,e_n$. For $n=2$,  
$\mathcal{T}\cong\mathcal{A}^g\cong\R^2$ is $\Orth(2)$-irreducible.
\end{prop}
The eight classes of linear connections are now defined by the possible
components of their torsions $T$ in these spaces.  The nice lecture notes 
by Tricerri and Vanhecke \cite{Tricerri&V1} use a similar approach in order
to classify homogeneous spaces by the algebraic properties of the torsion
of the canonical connection. They construct homogeneous examples of all
classes, study their ``richness'' and give explicit formulas for the
projections on every irreducible component of $\mathcal{T}$ in terms
of $\Orth(n)$-invariants.
\begin{dfn}[Connection with vectorial torsion]
%---------------------------------------------
The connection $\nabla$ is said to have  \emph{vectorial torsion} if
its torsion tensor lies in  the first space of the
decomposition in Proposition \ref{classessum}, i.\,e.~if it is essentially
defined by some vector field $V$ on $M$. The tensors $A$ and $T$
can then  be directly expressed through  $V$  as
\bdm
A(X,Y)\ =\ g(X,Y)V-g(V,Y)X,\quad
T(X,Y,Z)\ =\ g\big(g( V, X)Y- g(V, Y)X,Z\big).
\edm
These connections are particularly interesting on surfaces, in as much
that \emph{every} metric connection on a surface is of this type.

In \cite{Tricerri&V1}, F.~Tricerri and L.~Vanhecke  showed
that if $M$ is connected, complete,  simply-connected and
$V$ is $\nabla$-parallel, then $(M,g)$ has to be
isometric to the hyperbolic space.
V.~Miquel studied in \cite{Miquel82} and \cite{Miquel01} the
growth of geodesic balls of such connections, but did not
investigate the detailed shape of  geodesics. The study of the latter 
was outlined in \cite{Agricola&T04} (see Example \ref{exa-surfaces}), whereas 
\cite{Agri&F05} and \cite{Ivanov&P&P05} 
are devoted to holonomy aspects and a possible  role in superstring theory.

Notice that there is some similarity to Weyl geometry.
In both cases, we consider a Riemannian manifold with a fixed
vector field $V$ on it (\cite{Calderbank&P99}, \cite{Gauduchon95}).
A Weyl structure is a pair consisting of a conformal
class of metrics and a \emph{torsion-free non-metric connection} preserving 
the conformal structure. This connection is constructed by choosing 
a metric $g$ in the conformal class and is then defined by the formula
\bdm
\nabla^{\mathrm{w}}_X Y\ :=\ \nabla^g_X Y+g(X ,V)\cdot Y+ g(Y,V)X - g(X,Y)V .
\edm
Weyl geometry deals with the geometric properties of these connections, but in
spite of the resemblance, it
turns out to be a rather different topic. Yet in special geometric
situations it may happen that ideas from Weyl geometry can be useful.
\end{dfn}
\begin{dfn}[Connection with skew-symmetric torsion]
%--------------------------------------------------
The connection $\nabla$ is said to have  \emph{(totally) skew-symmetric 
torsion} if its torsion tensor lies in  the second component of the
decomposition in Proposition \ref{classessum}, i.\,e.~it is
given by a $3$-form. They are by now---for reasons to be
detailed later---a well-established tool in superstring theory and
weak holonomy theories (see  for example \cite{Strominger86},
\cite{Luest&T89}, \cite{Gauntlett01}, \cite{Curio&K&L01}, 
\cite{FigueroaPapadopoulos}, \cite{Duf02}, \cite{Agricola&F03a} etc.).
In Examples \ref{exa-nat-red} to \ref{exa-contact},  we describe large 
classes of
interesting manifolds that carry natural connections with skew-symmetric 
torsion. Observe that we can characterize these connections geometrically as 
follows:
\begin{cor}
%----------
A connection $\nabla$ on $(M^n,g)$ is metric and geodesic-preserving if 
and only if its torsion $T$ lies in $\Lambda^3(TM^n)$.
In this case, $2\, A=T$ holds,
\bdm
\nabla_X Y \ = \ \nabla^g_X Y \, + \, \frac{1}{2} \, T(X,Y,-),
\edm 
and the $\nabla$-Killing vector fields coincide with the Riemannian Killing 
vector fields. 
\end{cor}
In contrast to the case of vectorial torsion,
manifolds admitting invariant metric connections $\nabla$ with 
$\nabla$-parallel skew-symmetric torsion form a vast class that is worth
a separate investigation (\cite{Alexandrov03}, \cite{Cleyton&S04},
\cite{Schoemann06}). 
\end{dfn}

\smallskip
Suppose now that we are given a metric connection $\nabla$ with torsion on a
Riemannian spin manifold $(M^n,g)$ with spin bundle $\Sigma M^n$. We slightly 
modify our notation and write $\nabla$ as 
\bdm
\nabla_X Y\ :=\ \nabla^g_X Y+A_X Y,
\edm
where $A_X$ defines an endomorphism $TM^n\ra TM^n$ for every $X$. The
condition for $\nabla$ to be metric
\bdm
g(A_X Y, Z)+g(Y,A_X Z)\ =\ 0
\edm
means that $A_X$ preserves the scalar product $g$, which 
can be expressed as $A_X\in \so(n)$. After identifying $\so(n)$ with
$\Lambda^2(\R^n)$, $A_X$ can be written  relative to some orthonormal frame
\bdm
A_X\ =\ \sum_{i<j}\alpha_{ij}e_i\wedge e_j.
\edm
Since the lift into $\spin(n)$ of $e_i\wedge e_j$ is 
$E_i\cdot E_j /2$, $A_X$ defines an element in $\spin(n)$ and hence
an endomorphism of the spinor bundle. In fact, we need not introduce
a different notation for the lift of $A_X$. Rather, observe that
if $A_X$ is written as a $2$-form, 
\begin{enumerate}
\item its action on a vector $Y$ as
an element of $\so(n)$ is just $A_X Y = Y\haken A_X$, so our connection
takes on vectors the form
\bdm
\nabla_X Y\ =\ \nabla^g_X Y + Y\haken A_X,
\edm
\item the action of $A_X$ on a spinor $\psi$ as an element of $\spin(n)$ is 
just $A_X\psi = (1/2)\,A_X\cdot\psi$, where $\cdot$ denotes the Clifford 
product of a $k$-form by a spinor. The lift of the connection 
$\nabla$ to the spinor bundle $\Sigma M^n$ (again denoted by $\nabla$) is 
thus given by
\bdm
\nabla_X\psi \ =\ \nabla^g_X\psi+\frac{1}{2}\,A_X\cdot\psi.
\edm
\end{enumerate}
We denote by $(-,-)$ the Hermitian product on the spinor bundle $\Sigma M^n$ 
induced by $g$.
When lifted to the spinor bundle, $\nabla$ satisfies
the following properties that are well known for the lift of the
Levi-Civita connection. In fact, the proof easily follows from the 
corresponding properties for the Levi-Civita connection 
\cite[p.~59]{Dirac-Buch-00}
and the Hermitian product  \cite[p.~24]{Dirac-Buch-00}.
\begin{lem}
%----------
The lift of any metric connection $\nabla$ on $TM^n$ into the spinor bundle
$\Sigma M^n$ satisfies 
\bdm
\nabla_X(Y\cdot\psi)\ =\ (\nabla_X Y)\cdot\psi + Y\cdot (\nabla_X\psi),\quad
X ( \psi_1,\psi_2)\ =\ ( \nabla_X\psi_1,\psi_2) + (\psi_1,\nabla_X \psi_2).
\edm
\end{lem}
Any spinorial connection with the second property is again called 
\emph{metric}. The first property (chain rule for Clifford products) makes
only sense for spinorial connections that are lifts from the tangent bundle,
not for arbitrary spin connections.
%
%\begin{proof}
%%------------
%The proof follows from the corresponding properties of the Levi-Civita
%connection $\nabla^g$ and its lift to $SM$ (see \cite[p.~59]{Dirac-Buch-00}).
%For the chain rule, observe that for any $2$-form $\omega$, 
%$\omega\cdot Y = Y\cdot\omega + 2 Y\haken \omega$, hence we obtain
%%
%\bea[*]
%\nabla_X(Y\cdot\psi) & = & \nabla^g _X(Y\cdot\psi)+ \frac{1}{2}A_X\cdot Y
%\cdot\psi\\
%& =& (\nabla^g_X Y)\cdot\psi + Y\cdot (\nabla^g_X\psi)
%+\frac{1}{2} [ Y\cdot A_X + 2 Y\haken A_X]\cdot\psi\\
%&=& [\nabla^g_X Y+Y\haken A_X]\cdot\psi + Y\cdot[\nabla^g_X\psi+A_X\cdot\psi].
%\eea[*]
%%
%As for the property for $\nabla$ of being `metric' on $SM$,
%it is obviously equivalent to the condition 
%%
%\bdm
%(A_X\cdot\psi_1,\psi_2)+(\psi_1,A_X\cdot\psi_2)\ =\ 0
%\edm
%%
%for the $2$-form $A_X$. But the Hermitian product satisfies
%$(X\cdot\psi_1,\psi_2)=-(\psi_1,X\cdot\psi_2)$ \cite[p.~24]{Dirac-Buch-00},
%hence, for $i\neq j$:
%%
%\bdm
%(E_i E_j\psi_1,\psi_2)\ =\ - (E_j\psi_1,E_i\psi_2)\ =\ +
%(\psi_1,E_jE_i\psi_2)\ =\ - (\psi_1,E_iE_j\psi_2).\qedhere 
%\edm
%%
%\end{proof}
%
\begin{exa}[Connection with vectorial torsion]
%---------------------------------------------
For a metric connection with vectorial torsion given by $V\in TM$,
$A_X=2\,X\wedge V$, since
\bdm
Y\haken(2\, X\wedge V)\ =\ 2(X\wedge V)(Y,-)\ =\
(X\ox V)(Y,-) -(V\ox X)(Y,-) \ =\ g(X,Y)V-g(Y,V)X.
\edm
\end{exa}
\begin{exa}[Connection with skew-symmetric torsion]
%--------------------------------------------------
For a  metric connection with skew-symmetric torsion defined by some
$T\in \Lambda^3(M)$, $A_X=X\haken T$. Examples of manifolds
with a geometrically defined torsion $3$-form are given in the next section.
\end{exa}
\begin{exa}[Connection defined by higher order differential forms]
\label{exa-higher-forms}
%-------------------------------------------------------------------------
As example of a metric spinorial connection  not
induced from the tangent bundle, consider
\bdm
\nabla_X\psi\ :=\ \nabla^g_X\psi + (X\haken \omega^k)\cdot\psi
+ (X\wedge \eta^l)\cdot\psi
\edm
for some forms $\omega^k\in\Lambda^k(M)$, $\eta^l\in\Lambda^l(M)$
($k,l\geq 4$). These are of particular interest in string theory 
 as they are used for the description of higher 
dimensional membranes (\cite{AgriFri1}, \cite{Puhle06}).
\end{exa}
\begin{exa}[General case]
%------------------------
The class $\mathcal{A}'$ of Proposition \ref{classessum} cannot  be directly
interpreted as vectors or forms of a given degree, but it is not complicated 
to construct elements in $\mathcal{A}'$ either. For simplicity,
assume $n=3$, and let $\nabla$ be  the metric connection with vectorial
torsion $V=e_1$. Then
\bdm
A_X \ =\ X\wedge e_1\ =\ \big(\sum_{i=1}^3 g(X,e_i)e_i\big)\wedge e_1\ =\
g(X,e_2)\,e_2\wedge e_1 + g(X,e_3)\,e_3\wedge e_1.
\edm
Thus, the new form $\tilde{A}_X:=g(X,e_2)e_2\wedge e_1 - g(X,e_3)e_3
\wedge e_1 $ defines  a metric linear connection as well. One easily
checks that, as a tensor in $\mathcal{A}^g$, it is orthogonal to
$\Lambda^3(M) \oplus \mathfrak{X}(M)$, hence it lies in
$\mathcal{A}'$. Connections of this type have not yet been
investigated as a class of their own, but they are used as an interesting
tool in several contexts---for example, in closed $G_2$-geometry
(\cite{Bryant03}, \cite{Cleyton&I03}). The canonical connection of an
almost K\"ahler manifold is also of this type.
\end{exa}
What makes metric connections with torsion so interesting is the huge
variety of geometric situations that they unify in a  mathematically 
useful way. Let us illustrate this fact by some examples.
%
%-----------------------------------------------------------------------------
\subsection{Naturally reductive spaces}\label{exa-nat-red}
%---------------------------------------------------------------------------
%
Naturally reductive spaces are a key example of  manifolds with a 
metric connection with skew-symmetric torsion.

Consider a  Riemannian homogeneous space $M=G/H$. We suppose that
$M$ is \emph{reductive}, i.\,e.~the Lie algebra $\g$ of $G$ splits as
vector space direct sum of the Lie algebra $\h$ of $H$ and an
$\Ad(H)$-invariant subspace $\m$: $\g=\h\oplus\m$ and $\Ad(H)\m\subset\m$,
where $\mathrm{Ad}:\, H\ra \SO(\m)$ is the isotropy representation of $M$.
We identify $\m$ with $T_0M$ and we pull back
the Riemannian metric $\lan\ ,\ \ran_0$ on $T_0M$ to an inner product
$\lan\ ,\ \ran$ on $\m$. By a
theorem of Wang (\cite[Ch. X, Thm 2.1]{Kobayashi&N2}), there is a one-to-one
correspondence between the set of $G$-invariant metric affine connections
and the set of linear mappings $\Lambda_{\m}:\ \m\ra\so(\m)$ such that
 \bdm
 \Lambda_{\m}(hXh^{-1})\ =\ \Ad(h)\Lambda_{\m}(X)\Ad(h)^{-1} \
 \text{ for }X\in\m \text{ and } h\in H\,.
 \edm
A homogeneous Riemannian metric $g$ on $M$ is said to be \emph{naturally
reductive} (with respect to $G$) if the map 
$[X,-]_{\m}:\m\ra\m$ is skew-symmetric,
 \bdm
 g( [X,Y]_{\m},Z) + g( Y, [X,Z]_{\m}) \ =\ 0 \text{ for all }
 X,Y,Z\in\m\,.
 \edm
The family of metric connections $\nabla^t$ defined by
$\Lambda^t_\m (X)Y:=(1-t)/2\,[X,Y]_{\m}$ has then skew-symmetric torsion 
$T^t(X,Y)=-t[X,Y]_{\m}$. The connection $\nabla^1$ is of particular interest
and is called the \emph{canonical connection}. Naturally
reductive homogeneous spaces equipped with their canonical connection
 are a well studied (see for example \cite{DAtri&Z79}) generalization of 
symmetric spaces since they satisfy $\nabla^1 T^1=\nabla^1\kr^1=0$,
where $\kr^1$ denotes the curvature tensor of $\nabla^1$ (Ambrose-Singer,
\cite{Ambrose&S58}). In fact, a converse holds:
\begin{thm}[{\cite[Thm 2.3]{Tricerri&V84b}}]
%-------------------------------------------
A connected, simply connected and complete Riemannian manifold $(M,g)$
is a naturally reductive homogeneous space if and only if there exists
a skew-symmetric tensor field $T$ of type $(1,2)$ such that 
$\nabla:=\nabla^g-T$ 
 is a metric connection with $\nabla T = \nabla \kr =0$.
\end{thm}
The characterization of naturally reductive homogeneous spaces
given in \cite{Ambrose&S58} through the property that their geodesics are 
orbits of one-parameter subgroups of isometries is actually wrong; Kaplan's
$6$-dimensional Heisenberg group is the most prominent counterexample 
(see \cite{Kaplan83} and \cite{Tricerri&V84b}). Naturally reductive spaces
have been classified in small dimensions by Kowalski, Tricerri and Vanhecke, 
partially in the larger context of commutative spaces
(in the sense of Gel'fand): the  $3$-dimensional naturally reductive 
homogeneous spaces  are $\SU(2)$, the universal covering group of $\SL(2,\R)$ 
and the
Heisenberg group $H_3$, all with special families of left-invariant metrics 
(\cite{Tricerri&V1}). A simply connected four-dimensional naturally 
reductive space is either symmetric or 
decomposable as direct product (\cite{Kowalski&Vh83}). In dimension $5$, 
it is either symmetric, decomposable or locally isometric to
$\SO(3)\x\SO(3)/\SO(2)$,
 $\SO(3)\x H_3/\SO(2)$ (or any of these with $\SO(3)$ replaced by 
$\SL(2,\R)$), to the five-dimensional Heisenberg group $H_5$ or
to the Berger sphere $\SU(3)/\SU(2)$ (or $\SU(2,1)/\SU(2)$), all endowed 
with special families of metrics (\cite{Kowalski&Vh85}).

Other standard examples of naturally reductive spaces are
\begin{itemize}
\item Geodesic spheres in two-point homogeneous spaces,with the exception 
of the complex and quaternionic Cayley planes \cite{Ziller82}, 
\cite{Tricerri&V84a}
\item  Geodesic hyperspheres,
horospheres and tubes around totally geodesic non-flat complex space forms,
described and classified in detail by S.~Nagai \cite{Nagai95}, \cite{Nagai96},
\cite{Nagai97}
\item Simply connected $\vphi$-symmetric spaces \cite{Blair&V87}.
They are Sasaki manifolds with complete characteristic field
for which reflections
with respect to the integral curves of that field are global isometries.
\item All known left-invariant Einstein metrics on compact Lie groups
\cite{DAtri&Z79}. In fact, every simple Lie group apart from $\SO(3)$ and
$\SU(2)$
carries at least one naturally reductive  Einstein metric other than
the biinvariant metric. Similarly, large families of naturally reductive
Einstein metrics on  compact homogeneous spaces were constructed in
\cite{Wang&Z85}.
In contrast, non-compact naturally reductive Einstein manifolds
are necessarily symmetric \cite{Gordon&Z84}.
\end{itemize}
%
%----------------------------------------------------------------------------
\subsection{Almost Hermitian manifolds}\label{exa-almost-hermitian}
%----------------------------------------------------------------------------
%
An \emph{almost Hermitian manifold} $(M^{2n},g,J)$ is a manifold with a
Riemannian metric $g$ and a $g$-compatible almost complex structure
$J:TM^{2n}\ra TM^{2n}$. We denote by $\Omega(X,Y):=g(JX,Y)$ its
K\"ahler form and by   $N$ the Nijenhuis tensor of $J$, defined by
\bea[*]
N(X,Y)& :=& [JX,JY]-J[X,JY] -J[JX,Y]-[X,Y]\\
& =&  (\nabla^g_X J)JY -(\nabla^g_Y J)JX+(\nabla^g_{JX}J)Y- (\nabla^g_{JY}J)X,
\eea[*]
where the second expression follows directly from the vanishing of the torsion
of $\nabla^g$ and the identity 
\be\label{id-J}
(\nabla^g_X J)(Y)\ =\ \nabla^g_X(JY) - J(\nabla^g_X Y).
\ee
The reader is probably acquainted with the \emph{first canonical  Hermitian
connection}\footnote{By definition, a connection $\nabla$ is called 
\emph{Hermitian} if it is metric and has $\nabla$-parallel almost complex 
structure $J$.} (see the nice article \cite{Gauduchon97} by Paul Gauduchon,
which we strongly recommend for further reading on Hermitian connections)
\bdm
\ovn_X Y\ :=\ \nabla^g_X Y + \frac{1}{2}(\nabla^g_X J) JY.
\edm
Indeed, the condition $\ovn J = 0$ is equivalent to the identity 
($\ref{id-J}$),
and the antisymmetry of the difference tensor $g((\nabla^g_X J) JY,Z) $
in $Y$ and $Z$ can for example be seen from the standard identity\footnote{
For a proof, see \cite[Prop. 4.2.]{Kobayashi&N2}. Be aware of the 
different conventions in this book: $\Omega$ is defined with $J$ in the
second argument, the Nijenhuis tensor is twice our $N$ and 
derivatives of $k$-forms differ by a  multiple of  $1/k$, see 
\cite[Prop. 3.11.]{Kobayashi&N1}.}
\be\label{diff-J}
2\, g((\nabla^g_X J) Y,Z)\ =\ d\Omega(X,Y,Z) - d\Omega(X,JY,JZ)
+g(N(Y,Z),JX).
\ee
Sometimes, one finds the alternative formula
$- 1/2 J(\nabla^g_X J)Y$ for the difference tensor of $\ovn$; but
this is the same, since $J^2=-1$ implies $\nabla^g_X J^2=0=(\nabla^g_XJ)J
+J(\nabla_X J)$, i.\,e.~$\nabla^g J\in \un(n)^\perp\subset\so(2n)$.
Let us now express the difference tensor of the connection $\ovn$ using
the Nijenhuis tensor and the K\"ahler form. Since 
$\nabla^g_X\Omega(Y,Z)=g\big((\nabla^g_XJ)Y,Z\big)$, the differential
$d\Omega$ is just
\be\label{diff-omega}
d\Omega(X,Y,Z)\ =\ g\big((\nabla^g_XJ)Y,Z\big) - g\big((\nabla^g_YJ)X,Z\big)
+g\big((\nabla^g_ZJ)X,Y\big).
\ee
Together with the expression for
$N$ in terms of covariant derivatives of $J$, this yields 
\bdm
g((\nabla^g_X J) JY,Z)\ = \ g\big(N(X,Y),Z\big) + d\Omega(JX,JY,JZ)
- g((\nabla^g_Y J)JZ,X) - g((\nabla^g_{JZ} J) Y,X).
\edm
A priori, $(\nabla^g_XJ)Y$ has no particular symmetry properties in $X$ 
and $Y$, hence the last two terms cannot be simplified any further
(in general, they are a mixture of the two other Cartan types).
An exceptional situation occurs if $M$ is nearly K\"ahler 
($(\nabla^g_XJ)X=0$), 
for then $(\nabla^g_XJ)Y=-(\nabla^g_YJ)X$ and the last two terms  
cancel each other. Furthermore, this antisymmetry property
implies that the difference tensor is totally skew-symmetric, hence
we can conclude:
\begin{lem}\label{char-connection-nK}
%------------------------------------
On a nearly K\"ahler manifold $(M^{2n},g,J)$, the formula
\bdm
\ovn_X Y\ :=\ \nabla^g_X Y + \frac{1}{2}(\nabla^g_X J) JY\ =\ 
\nabla^g_X Y + \frac{1}{2}\left[N(X,Y) +d\Omega(JX,JY,J-)\right]
\edm
defines a  Hermitian connection with totally skew-symmetric torsion.
\end{lem}
This connection was first defined and studied by Alfred Gray
(see \cite[p.~304]{Gray70} and \cite[p.~237]{GR}). It is a non-trivial
result of Kirichenko that it has $\ovn$-parallel torsion
(\cite{Kirichenko77}, see also \cite{Alexandrov&F&S04} for a modern index-free
proof). Furthermore, it is shown in \cite{GR} that any $6$-dimensional
nearly K\"ahler manifold is Einstein and of constant type, i.\,e.~it
satisfies 
\bdm
||(\nabla^g_XJ)(Y)||^2\ =\ \frac{\Scal^g}{30}\left[||X||^2\cdot ||Y||^2 
-g(X,Y)^2 - g(X,JY)^2\right].
\edm
Together with Lemma \ref{char-connection-nK}, this identity yields by
a direct calculation that \emph{any $6$-dimensional nearly
K\"ahler manifold is also $\ovn$-Einstein with} 
$\Ric^{\ovn}=2(\Scal^g/15)\,g$ (see Theorem \ref{curvature-identities}
for the relation between Ricci tensors).

Now let us look for  a  Hermitian connection with
totally skew-symmetric torsion on a larger class of Hermitian
manifolds generalizing nearly K\"ahler manifolds. 
\begin{lem}\label{lem-almost-Herm-conn}
%--------------------------------------
Let $(M^{2n},g,J)$ be an almost Hermitian manifold with
 skew-symmetric Nijenhuis tensor $N(X,Y,Z):=g(N(X,Y),Z)$. Then
the formula
\bdm
g(\nabla_X Y, Z) \ :=\  
g(\nabla^g_X Y,Z) + \frac{1}{2}\left[N(X,Y,Z) +d\Omega(JX,JY,JZ)\right]
\edm
defines a  Hermitian connection with skew-symmetric torsion.
\end{lem}
\begin{proof}
%------------
Obviously, only $\nabla J=0$ requires a calculation. 
By (\ref{id-J}) and the definition of $\nabla$, we have
\bea[*]
2\nabla_XJ(Y) &=& 2g\big(\nabla_X(JY)-J(\nabla_XY),Z\big)\ =\ 
2g\big(\nabla_X(JY),Z\big)+ 2 g\big(\nabla_XY,JZ\big)\\
&=& 2\nabla^g_XJ(Y) +N(X,JY,Z)-d\Omega(JX,Y,JZ)+N(X,Y,JZ)-d\Omega(JX,JY,Z).
\eea[*]
But from the symmetry properties of the Nijenhuis tensor and the metric, 
one sees
\bdm
N(X,Y,JZ)\ =\ g(N(X,Y),JZ)\ =\ -g(JN(X,Y),Z)\ =\ g(N(JX,Y),Z)\ =\ N(JX,Y,Z)
\edm
and repeated application of the identity (\ref{diff-omega}) for $d\Omega$
yields 
\bdm
3 N(JX,Y,Z)\ =\ d\Omega(X,JY,JZ)-d\Omega(X,Y,Z)+d\Omega(JX,JY,Z)+
d\Omega(JX,Y,JZ).
\edm
Together, these two equations show that the previous expression for
$2\nabla_XJ(Y)$ vanishes by equation (\ref{diff-J}).
\end{proof}
Besides nearly K\"ahler manifolds, Hermitian manifolds ($N=0$) trivially
fulfill the condition of the preceding lemma and $\nabla$ coincides then
with the Bismut  connection; however, in the non-Hermitian situation, $\nabla$
is \emph{not} in the standard family of canonical Hermitian connections
that is usually considered (see \cite[2.5.4]{Gauduchon97}).
Proposition 2 in this same reference gives the decomposition of the torsion
of any Hermitian connection in its $(p,q)$-components and gives
another justification for this precise form for the torsion.   Later, we shall 
see that $\nabla$
is the \emph{only possible} Hermitian connection with skew-symmetric torsion 
and that the class of almost Hermitian manifolds with  skew-symmetric 
Nijenhuis tensor is the largest possible where it is defined. 

We will put major emphasis  on almost Hermitian manifolds
of dimension $6$, although one will find some general results 
formulated independently of the dimension. Two reasons for this choice
are that nearly K\"ahler manifolds are of interest only in dimension $6$,  
and that $6$ is also the relevant dimension in superstring theory. 
%
%------------------------------------------------------------------------
\subsection{Hyper-K\"ahler manifolds with torsion (HKT-manifolds)}
\label{exa-HKT}
%------------------------------------------------------------------------
%
We recall that a manifold $M$ is called \emph{hypercomplex} if it is endowed
with three (integrable) complex structures $I,J,K$ satisfying the
quaternionic identities $IJ=-JI=K$. A metric $g$ compatible with these
three complex structures (a so-called \emph{hyper-Hermitian} metric) is said 
to be \emph{hyper-K\"ahler with torsion}
or just a \emph{HKT-metric} if the K\"ahler forms satisfy the identity
\be\label{HKT}
I\, d\Omega_I\ =\ J\, d\Omega_J\ =\ K\,d\Omega_K.
\ee
Despite  the misleading name, these manifolds are not
K\"ahler (and hence even less hyper-K\"ahler). HKT-metrics
were introduced by Howe and Papadopoulos as target spaces of some
two-dimensional sigma models with $(4,0)$ supersymmetry with  Wess-Zumino
term \cite{Howe&P96}.
Their mathematical description was given by Grantcharov and Poon in 
\cite{Grantcharov&P00} and further investigated by several authors
since then (see for example \cite{Verbitsky02}, \cite{Dotti&F},
\cite{Poon&S03}, \cite{Fino&G03}, \cite{Fino&P&S04}, 
\cite{Ivanov&M04}). From the 
previous example, we can conclude immediately that 
\bdm
g(\nabla_X Y, Z) \ :=\  
g(\nabla^g_X Y,Z) + \frac{1}{2}\, d\Omega_I(IX,IY,IZ)
\edm
defines a metric connection with skew-symmetric torsion such that 
$\nabla I = \nabla J = \nabla K=0$; one easily checks that $\nabla$ is again
the \emph{only} connection fulfilling these conditions.
 Equation (\ref{HKT}) implies that we could equally well have chosen  $J$ 
or $K$ in the last term. In  general, a hyper-Hermitian manifold will not
carry an HKT-structure, except in dimension $4$ where this is proved in
\cite{Gauduchon&T98}. Examples of  homogeneous HKT-metrics can be
constructed using a family of homogeneous hypercomplex structures
associated with compact semisimple Lie groups constructed by Joyce
\cite{Joyce92}. Inhomogeneous HKT-structures exist for
example on $S^1\times S^{4n-3}$ \cite{Grantcharov&P00}.
The question of existence of suitable potential functions for HKT-manifolds
was first raised and discussed in the context of super-conformal quantum 
mechanics by the physicists Michelson and Strominger 
\cite{Michelson&S00} (a maximum principle argument shows that compact 
HKT-manifolds do not admit global potentials); Poon and Swann 
discussed potentials for some symmetric HKT-manifolds \cite{Poon&S01},
while Banos and Swann were able to show local existence
\cite{Banos&S04}.
%
%-------------------------------------------------------------------------
\subsection{Almost contact metric structures}\label{exa-contact}
%-------------------------------------------------------------------------
%
An odd-dimensional manifold $M^{2n+1}$ is said to carry an
\emph{almost contact structure} if it admits a $(1,1)$-tensor field
$\vphi$ and a vector field $\xi$ (sometimes called the \emph{characteristic}
or \emph{Reeb vector field}) with dual $1$-form $\eta$ ($\eta(\xi)=1$)
such that $\vphi^2=-\Id+\eta\ox\xi$. Geometrically, this means
that $M$ has a preferred direction (defined by $\xi$) on which
$\vphi^2$ vanishes, while $\vphi$ behaves like an almost complex
structure on any linear complement of $\xi$. An easy argument
shows that $\vphi(\xi)=0$ \cite[Thm 4.1]{Blair02}.
If there exists in
addition a $\vphi$-compatible Riemannian metric $g$ on $M^{2n+1}$, 
i.\,e.~satisfying
\bdm
g(\vphi X, \vphi Y)\ =\ g(X,Y)-\eta(X)\eta(Y),
\edm
then we say that $(M^{2n+1}, g,\xi,\eta,\vphi)$ carries an  \emph{almost 
contact metric structure} or that it is an \emph{almost contact 
metric manifold}. The condition says that $\xi$ is
a vector field of unit length with respect to $g$ and that $g$ is 
$\vphi$-compatible in the sense of Hermitian geometry on the
orthogonal complement $\xi^\perp$. Unfortunately, this relatively
intuitive structural concept  splits into a myriad of
subtypes and leads to complicated equations in the defining 
data $(\xi,\eta,\vphi)$, making the investigation of almost contact 
metric structures
look rather unattractive at first sight (see \cite{Aleksiev&G86},
\cite{Chinea&G90}, \cite{Chinea&M92}, and \cite{Fino95} for a
classification). Yet, they constitute
a rich and particularly interesting class of non-integrable geometries,
as they have no integrable analogue on Berger's list.
An excellent general
source on contact manifolds with extensive references
are the books by David Blair, \cite{Blair76} and \cite{Blair02}
(however, the classification is not treated in these).
\begin{exa}\label{exa-Hermitian-contact}
%---------------------------------------
For every almost Hermitian manifold $(M^{2n},g,J)$, there exists an almost 
metric contact structure $(\tilde{g},\xi,\eta,\vphi)$ on the cone 
$M^{2n}\x\R$ with the product metric. For any function
$f\in C^\infty(M^{2n}\x\R)$ and vector field
$X\in\mathfrak{X}(M^{2n})$, it is defined by
\bdm
\vphi\left(X,f \del_t\right)\ =\ (JX,0),\quad
\xi=\left(0,\del_t\right),\quad \eta\left(X,f\del_t\right)\ =\ f.
\edm
Conversely, an almost metric contact structure $(g,\xi,\eta,\vphi)$
on $M^{2n+1}$ induces an almost Hermitian structure $(\tilde{g},J)$ on its cone
$M^{2n+1}\x\R$ with product metric by setting for $f\in
C^\infty(M^{2n+1}\x\R)$ and $X\in\mathfrak{X}(M^{2n+1})$
\bdm
J(X,f\del_t)\ =\ (\vphi X-f\xi,\eta(X)\del_t).
\edm
In fact, one shows that
every smooth orientable hypersurface $M^{2n-1}$ in an
almost Hermitian manifold $(M^{2n},g,J)$ carries a canonical almost contact
metric structure \cite[4.5.2]{Blair02}. In this way, one
easily constructs almost contact metric structures on compact manifolds, 
for example on all odd dimensional spheres. 
\end{exa}
The fundamental form $F$ of an almost contact metric structure
is defined by $F(X,Y)=g(X,\vphi(Y))$, its Nijenhuis tensor
is given by a similar, but slightly more complicated formula as in
the almost  Hermitian case and can also be written in terms of
covariant derivatives of $\vphi$,
\bea[*]
N(X,Y)& :=& [\vphi (X),\vphi (Y)]-\vphi[X,\vphi (Y)] -\vphi [\vphi (X),Y]
+\vphi^2[X,Y]+ d\eta(X,Y)\xi\\
& =&  (\nabla^g_X \vphi)\vphi(Y) -(\nabla^g_Y \vphi)\vphi(X)
+(\nabla^g_{\vphi (X)}\vphi)Y- (\nabla^g_{\vphi (Y)}\vphi)X
+\eta(X)\nabla^g_Y \xi - \eta(Y)\nabla^g_X \xi .
\eea[*]
Let us emphasize some particularly interesting cases. A  manifold
with an almost 
contact metric structure $(M^{2n+1}, g,\xi,\eta,\vphi)$ is called
\begin{enumerate}
\item a \emph{normal almost contact metric manifold} if $N=0$,
\item a \emph{contact metric manifold} if $2F=d\eta$. 
\end{enumerate}
Furthermore, a contact metric structure is said to be
a \emph{K-contact metric structure} if $\xi$ is in addition a Killing 
vector field, and  a \emph{Sasaki structure} if it is normal (it is then
automatically K-contact, see \cite[Cor. 6.3.]{Blair02}). 
\emph{Einstein-Sasaki manifolds} are just Sasaki manifolds whose
$\vphi$-compatible Riemannian metric $g$ is Einstein.
Without doubt, the forthcoming monograph by Ch.~Boyer and Kr.~Galicki
on Sasakian geometry \cite{Boyer&G07} is set to become \emph{the}
standard reference for this area of contact geometry in the future; in the 
meantime, the reader will have to be contented with the shorter 
reviews \cite{BG} and \cite{BG01}. 

Much less is known about metric connections on almost contact metric 
manifolds than on almost Hermitian manifolds. In fact, only the so-called 
\emph{generalized Tanaka connection} (introduced by S.~Tanno) has been 
investigated. It is a metric connection defined on the class of contact 
metric manifolds by the formula
\bdm
\nabla^*_X Y\ :=\ \nabla^g_X Y +\eta(X)\vphi(Y) - \eta(Y)\nabla^g_X\xi
+(\nabla_X\eta)(Y)\xi 
\edm
satisfying the additional conditions $\nabla^*\eta=0$ (which is of course
equivalent to $\nabla^*\xi=0$), see \cite{Tanno89} and \cite[10.4.]{Blair02}.
One easily checks that its torsion is not skew-symmetric, not even in
the Sasaki case. In fact, from the point of view of non-integrable
structures, it seems appropriate to require in  addition 
 $\nabla^*\vphi=0$ (compare with
the almost Hermitian situation). 

Following a similar but more complicated line of arguments as in the
almost Hermitian case,  Th.~Friedrich and S.~Ivanov showed:
\begin{thm}[{\cite[Thm 8.2.]{Friedrich&I1}}]\label{thm-contact}
%--------------------------------------------------------------
Let $(M^{2n+1}, g,\xi,\eta,\vphi)$ be an almost contact metric
manifold. It admits a metric connection $\nabla$ with totally skew-symmetric 
torsion $T$ and $\nabla\eta=\nabla\vphi=0$  if and only if
the Nijenhuis tensor $N$ is skew-symmetric and if $\xi$ is a Killing vector
field. Furthermore, $\nabla=\nabla^g+ (1/2)T$  is uniquely determined by 
\bdm
T\ =\ \eta\wedge d\eta+ d^\vphi F+N - \eta\wedge(\xi\haken N),
\edm
where $d^\vphi F$ stands for the $\vphi$-twisted derivative,
$d^\vphi F(X,Y,Z):= - dF(\vphi(X),\vphi(Y), \vphi(Z))$.
\end{thm}
For a Sasaki structure, $N=0$ and $2F=d\eta$ implies $d^\vphi F=0$, hence
$T$ is given by the much simpler formula $T=\eta\wedge d\eta$. 
This connection had been noticed before, for example in \cite{Kowalski&W87}.
In fact, one  sees that  $\nabla T=0$ holds, hence Sasaki
manifolds endowed with this connection are examples of 
non-integrable geometries with parallel torsion. A.~Fino studied
naturally reductive almost contact metric structures such that $\vphi$
is parallel with respect to the canonical connection in \cite{Fino94}.
In general, potentials are hardly studied in contact geometry (compare
with the situation for HKT-manifolds), but 
a suitable analogue of the K\"ahler potential was constructed on Sasaki 
manifolds by M.~Godlinski, W.~Kopczynski and P.~Nurowski 
\cite{Godlinski&K&N00}.
%
%--------------------------------------------------------------------------
\subsection{$3$-Sasaki manifolds}\label{exa-3-Sasaki}
%---------------------------------------------------------------------------
Similarly to HKT-manifolds and quaternionic-K\"ahler manifolds, it makes
sense to investigate configurations with three `compatible' almost
metric contact structures $(\vphi_i,\xi_i,\eta_i)$, $i=1,2,3$ on 
$(M^{2n+1},g)$ for some fixed metric $g$.
The compatibility condition may be formulated as
\bdm
\vphi_k\ =\ \vphi_i\vphi_j-\eta_j\ox\xi_i\ =\ 
-\vphi_j\vphi_i+\eta_i\ox\eta_j,\quad
\xi_k\ =\ \vphi_i\xi_j\ =\ -\vphi_j\xi_i
\edm
for any even permutation $(i,j,k)$ of $(1,2,3)$, and  
such a structure is called an \emph{almost contact metric $3$-structure}.
By defining on the cone $M^{2n+1}\x\R$ three almost complex structures 
$J_1,J_2,J_3$ as outlined in Example 
\ref{exa-Hermitian-contact}, one  sees that the cone carries an almost
quaternionic structure and hence has dimension divisible by $4$.
Consequently, almost contact metric $3$-structures exist only in
dimensions $4n+3, n\in\mathbb{N}$, and it is  no surprise that
the structure group of its tangent bundle turns out to be
contained in $\Sympl(n)\x\{1\}$. What is surprising is the recent result
by T.~Kashiwada that if all three structures  $(\vphi_i,\xi_i,\eta_i)$
are contact metric structures, they automatically have to be Sasakian
\cite{Kashiwada01}. A manifold with such a structure will be called a 
\emph{$3$-Sasaki(an) manifold}. An earlier result by T.~Kashiwada claims that
any $3$-Sasaki manifold is Einstein \cite{Kashiwada71}.
The canonical example of a $3$-Sasaki manifold is the sphere
$S^{4n+3}$ realized as a hypersurface in $\mathbb{H}^{n+1}$: each of the
three almost complex structures forming the quaternionic structure of
$\mathbb{H}^{n+1}$ applied to the exterior normal vector field of the sphere
yields a vector field $\xi_i$ ($i=1,2,3$) on $S^{4n+3}$, leading thus to
three orthonormal vector fields on $S^{4n+3}$.
Th.~Friedrich and I.~Kath showed that every compact simply connected 
$7$-dimensional spin manifold with regular $3$-Sasaki structure
is isometric to $S^7$ or the Aloff-Wallach space $N(1,1)=\SU(3)/S_1$
(see \cite{FK} or \cite{BFGK}). By now,  it is possible to
list all homogeneous $3$-Sasaki manifolds:
\begin{thm}[{\cite{Boyer&G&M94}}]
%--------------------------------
A homogeneous $3$-Sasaki manifold is isometric to one of the following:
\begin{enumerate}
\item Four families: $\frac{\Sympl(n+1)}{\Sympl(n)}\cong S^{4n+3}$,
$S^{4n+3}/\mathbb{Z}_2\cong \R\P^{4n+3}$,
$\frac{\SU(m)}{\mathrm{S}(\U(m-2)\x\U(1))}$ for $m\geq 3$,
$\frac{\SO(k)}{\SO(k-4)\x\Sympl(1)}$ for $k\geq 7$.
\item Five exceptional spaces: $G_2/\Sympl(1)$, $F_4/\Sympl(3)$,
$E_6/\SU(10)$, $E_7/\Spin(12)$, $E_8/E_7$.
\end{enumerate}
All these spaces fibre over a quaternionic K\"ahler manifold; the
fibre is $\Sympl(1)$ for   $S^{4n+3}$ and $\SO(3)$ in all other cases.
\end{thm}
Many non-homogeneous examples have also been constructed. The analogy
between $3$-Sasaki manifolds and HKT-manifolds breaks down when one starts
looking at connections, however. For an arbitrary $3$-Sasaki structure
with $1$-forms $\eta_i$ ($i=1,2,3$), each of the three underlying
Sasaki structures yields one possible choice of a  metric connection 
$\nabla^i$ with $\nabla^i\eta_i=0$ and torsion $T^i=\eta_i\wedge d\eta_i$ 
as detailed in Theorem \ref{thm-contact}. \emph{However, these three
connections do not coincide}; hence, the $3$-Sasaki structure itself is not
preserved by any metric connection with skew-symmetric torsion. 
A detailed discussion of these three connections and their spinorial 
properties in dimension $7$ can be found in \cite{Agricola&F03a}.
Nevertheless, $3$-Sasaki manifolds have recently appeared and been 
investigated in the context of the AdS/CFT-correspondence 
by Martelli, Sparks and Yau \cite{Martelli&S&Y06}.
%
%--------------------------------------------------------------------------
\subsection{Metric connections on surfaces}\label{exa-surfaces}
%---------------------------------------------------------------------------
%
Classical topics of surface theory like the Mercator projection can be 
understood in a different light with the help of metric connections with 
torsion. 

In \cite[\S\,67, p.\,408--409]{Cartan23a}, Cartan describes the two-dimensional
sphere with its flat metric connection, and observes (without proof) that 
``on this manifold, the straight lines are the \emph{loxodromes}, which 
intersect the
meridians at a constant angle. The only straight lines realizing shortest
paths are those  normal to the torsion at every point: \emph{these
are the meridians}\footnote{ \flqq Sur cette vari\'et\'e, les lignes droites sont
les \emph{loxodromies}, qui font un angle constant avec les m\'eridiennes.
Les seules lignes droites qui r\'ealisent les plus courts chemins sont celles
qui sont normales  en chaque point \`a la torsion~: \emph{ce sont les 
m\'eridiennes.\frqq}}''. 

%------------------------------------------------------------
\begin{figure}%\label{drehflaeche}
\bdm
\psfrag{x}{$x$}\psfrag{y}{$y$}\psfrag{z}{$z$}
\psfrag{a}{$h(s)$}\psfrag{r}{$r(s)$}
\psfrag{p}{$\varphi$}
\psfrag{v1}{$\nu_1$}\psfrag{v2}{$\nu_2$}
\includegraphics[width=6.6cm]{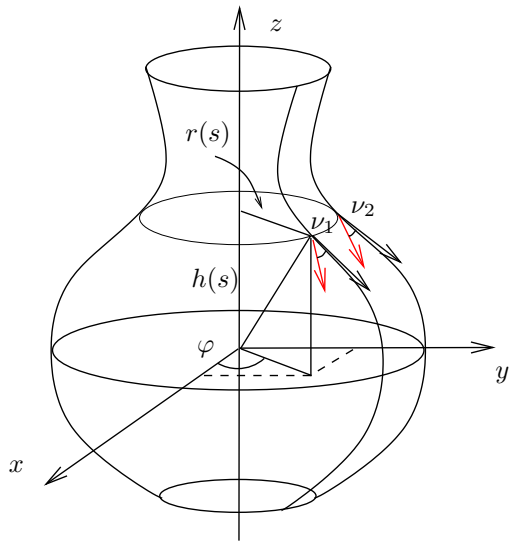}
\edm
\caption{Surface of revolution generated by a curve $\alpha$.}
\label{drehflaeche}    
\end{figure}
%-------------------------------------------------------------
%

This suggests that  there exists a class of  metric connections on surfaces 
of revolution whose geodesics admit a generalization of Clairaut's theorem,
yielding loxodromes in the case of the flat connection. Furthermore,
it is well known that the Mercator projection maps loxodromes to 
straight lines in the plane (i.\,e., Levi-Civita geodesics of the Euclidian
metric), and that
this mapping is conformal. Theorem~\ref{thm-grad-VF} provides the right
setting to understand both effects:

\begin{thm}[\cite{Agricola&T04}]\label{thm-grad-VF}
%---------------------------------------------------
Let $\sigma$ be a function on the Riemannian manifold $(M,g)$, 
$\nabla$ the metric connection with vectorial torsion defined by
$V=-\grad(\sigma)$, and consider the conformally equivalent metric
$\tilde{g}=e^{2\sigma}g$. Then:
\begin{enumerate}
\item Any $\nabla$-geodesic $\gamma(t)$ is, up to
a reparametrization $\tau$, a $\nabla^{\tilde{g}}$-geodesic, and 
the function $\tau$ is the unique solution of the differential
equation $\ddot{\tau}+\dot{\tau}\dot{\sigma}=0$, where we set
$\sigma(t):=\sigma\circ\gamma\circ\tau (t)$;
\item If $X$ is a Killing field for the metric $\tilde{g}$, the function
$e^{\sigma}g(\dot{\gamma}, X)$ is a constant of motion for the
$\nabla$-geodesic $\gamma(t)$.
\item The connection forms of $\nabla$ and $\nabla^{\tilde{g}}$
coincide; in particular, they have the same curvature.
\end{enumerate}
\end{thm}
%
%--------------------------------------------------------------
\begin{figure}
\bdm
\includegraphics[width=8cm, bbllx=50, bblly=230, bburx=560,bbury=590]{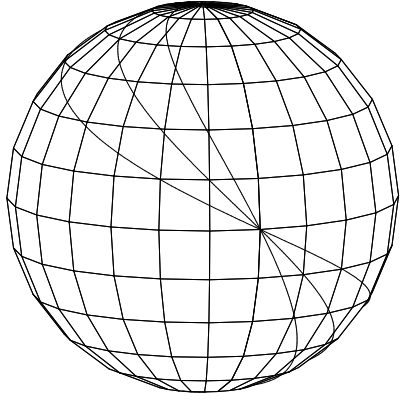}
\edm
\caption{Loxodromes on the sphere.}\label{sphaere}
\end{figure}
%--------------------------------------------------------------
%
We discuss Cartan's example in the light of Theorem \ref{thm-grad-VF}.
Let $\alpha=(r(s),h(s))$ be a curve parametrized by arclength,
and $M(s,\vphi)=(r(s)\cos\vphi,r(s)\sin\vphi,h(s))$ the surface of
revolution generated by it. The first fundamental form is 
$g=\diag(1,r^2(s))$, so we can choose the orthonormal frame
$e_1=\del_s$, $e_2=(1/r)\del_\vphi$ with dual $1$-forms
$\sigma^1=ds$, $\sigma^2=r\,d\vphi$. We define a connection $\nabla$ by
calling two tangential vectors
$v_1$ and $v_2$ parallel if the angles $\nu_1$ and $\nu_2$ with
the meridian through their foot point coincide (see Figure \ref{drehflaeche}).
Hence $\nabla e_1=\nabla e_2=0$, and the connection $\nabla$ is flat.
But for a flat connection, the torsion $T$ is can be derived from
$d\sigma^i(e_j,e_k)=\sigma^i(T(e_j,e_k))$. Since $d\sigma^1=0$ and
$d\sigma^2=(r'/r)\sigma^1\wedge\sigma^2$, one obtains 
\bdm
T(e_1,e_2)\ =\ \frac{r'(s)}{r(s)}\,e_2 \ \text{ and }\
V\ =\ \frac{r'(s)}{r(s)}\,e_1\ =\ -\grad \big(- \ln r(s)\big).
\edm
Thus, the metric connection $\nabla$ with vectorial torsion $T$ is
determined by the gradient of the function $\sigma:=- \ln r(s)$. 
By Theorem \ref{thm-grad-VF}, we conclude that its geodesics are
the Levi-Civita geodesics of the conformally equivalent metric
$\tilde{g}=e^{2\sigma}g=\diag(1/r^2, 1)$. This coincides with the standard
Euclidian metric by changing  variables
$x=\vphi$, $y=\int ds/r(s)$. For example, the sphere is obtained
for $r(s)=\sin s, h(s)=\cos s$, hence $y=\int ds/\sin s = \ln\tan(s/2)$
 ($|s|<\pi/2$), and
this is precisely the coordinate change of the Mercator projection.
Furthermore, $X=\del_\vphi$ is a Killing vector field for $\tilde{g}$, 
hence the second part of Theorem \ref{thm-grad-VF}
yields for a $\nabla$-geodesic $\gamma$ the invariant of motion
\bdm
\mathrm{const}\ =\ e^\sigma g(\dot{\gamma},X)\ =\
\frac{1}{r(s)} g(\dot{\gamma}, \del_\vphi)\ = \ g(\dot{\gamma},e_2),
\edm
i.\,e.~the cosine of the angle between $\gamma$ and a parallel
circle. This shows that $\gamma$ is a loxodrome on $M$, as claimed
(see Figure \ref{sphaere} for loxodromes on the sphere). In the same way,
one obtains a ``generalized Clairaut theorem'' for any gradient
vector field on a surface of revolution. For the
pseudosphere, one chooses
\bdm
r(s)\ =\ e^{-s},\quad h(s)\ =\ \arctanh \sqrt{1-e^{-2s}} - \sqrt{1-e^{-2s}},
\edm
hence $V=-e_1$ and $\nabla V=0$, in accordance with the results by
\cite{Tricerri&V1} cited before. Notice that $X$ is also
a Killing vector field for the metric $g$ and \emph{does} commute with $V$;
nevertheless, $g(\dot{\gamma},X)$ is \emph{not} an invariant of motion.
 
The catenoid is another interesting example:
since it is a minimal surface, the Gauss map to the sphere is 
conformal, hence it maps loxodromes to loxodromes. Thus,
Beltrami's theorem (``If a portion of a surface $S$ can be mapped 
LC-geodesically onto a portion of a surface $S^*$ of constant Gaussian
curvature, the Gaussian curvature of $S$ must also be constant'',
see for example \cite[\S 95]{Kreyszig2}) does \emph{not} hold for metric
connections with vectorial torsion---the sphere is a surface
of constant Gaussian curvature, but the catenoid is not.
\begin{NB}
%---------
The unique flat metric connection $\nabla$ does not have to be of vectorial 
type. For example, on the compact Lie group $\SO(3)$ the torsion is 
a $3$-form: Fix an orthonormal basis $e_1,e_2,e_3$ with
commutator relations $[e_1,e_2]=e_3, [e_2,e_3]=e_1$ and $[e_1,e_3]=-e_2$.
Cartan's structural equations then read $d\sigma^1=\sigma^2\wedge\sigma^3,
 d\sigma^2=-\sigma^1\wedge\sigma^3,d\sigma^3=\sigma^1\wedge\sigma^2$,
from which we deduce $T=2A=\sigma^1\wedge\sigma^2\wedge\sigma^3$.
In particular, $\nabla$ has the same geodesics as $\nabla^g$.
\end{NB}
%

%-----------------------------------------------------------------------------
\subsection{Holonomy theory}
%-----------------------------------------------------------------------------
%
Let  $(M^n,g)$ be a (connected) Riemannian manifold equipped with \emph{any}
connection $\nabla$. For a  curve $\gamma(s)$ from $p$ to $q$, parallel 
transport along $\gamma$ is the linear mapping $P_\gamma: T_pM\ra T_qM$ 
such that the vector field $\V$ defined by 
\bdm
\V(q)\ :=\ P_\gamma\V(p) \  \text{ along }\gamma
\edm
is parallel along $\gamma$, $\nabla\V(s) /ds= \nabla_{\dot{\gamma}}\V=0$. 
$P_\gamma$ is always an invertible endomorphism, hence,  for a closed loop 
$\gamma$ through
$p\in M$, it can be viewed as an element of $\GL(n,\R)$
(after choice of some basis).
\begin{figure}
\bdm
\psfrag{p}{$p$}\psfrag{TpM}{$T_pM$}\psfrag{A}{$P_\gamma$}
\psfrag{g}{$\gamma$}\psfrag{M}{$M$}
\includegraphics[width=5.5cm]{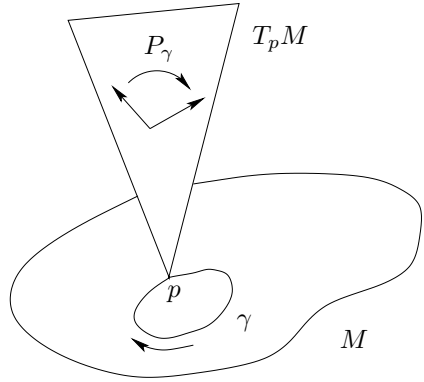}
\edm
\caption{Schematic concept of holonomy.}
\end{figure}
Consider the \emph{loop space} $C(p)$ of all closed, piecewise smooth
curves through $p$, and therein the subset $C_0(p)$ of   curves
that are homotopic to the identity. The set of parallel 
translations along loops in $C(p)$ or $C_0(p)$ forms
a group acting on $\R^n\cong T_pM$, called the \emph{holonomy group 
$\Hol(p;\nabla)$ of $\nabla$}   or the \emph{restricted holonomy group 
$\Hol_0(p;\nabla)$ of $\nabla$}  at the point $p$.
 Let us now change the point of view from $p$ to $q$, $\gamma$ a path joining
them;
then   $\Hol(q;\nabla)= P_\gamma \Hol(p;\nabla)P_\gamma^{-1}$ and similarly for
$\Hol_0(p;\nabla)$. Hence, all holonomy groups are isomorphic, so we  
drop the base point from now on. Customary notation for them is 
$\Hol(M;\nabla)$ and $\Hol_0(M;\nabla)$. Their action on $R^n\cong T_pM$ 
shall be called the \emph{(restricted) holonomy representation}.

In general, it is only known that (\cite[Thm IV.4.2]{Kobayashi&N1})
\begin{enumerate}
\item $\Hol(M;\nabla)$ is a Lie subgroup of $\GL(n,\R)$,
\item $\Hol_0(p)$ is the connected component of the identity of
$\Hol(M;\nabla)$.
\end{enumerate}
If one assumes in addition---as we will do through this text---that
$\nabla$ be metric, parallel transport becomes an isometry: for
any two parallel vector fields $\V(s)$ and $\W(s)$, being metric implies
\bdm
\frac{d}{ds}g\big( \V(s),\W(s)\big)\ =\ 
g\big( \frac{\nabla\V(s)}{ds},\W(s)\big) + 
\big( \V(s),\frac{\nabla\W(s)}{ds}\big)\ =\ 0.
\edm 
Hence, $\Hol(M;\nabla)\subset \Orth(n)$ and $\Hol_0(M;\nabla)\subset
\SO(n)$. 
For convenience, we shall henceforth speak of the
\emph{Riemannian (restricted) holonomy group} if $\nabla$ is
the Levi-Civita connection, to distinguish it from holonomy groups in our
more general setting. 
\begin{exa}
%----------
This is a good moment to discuss Cartan's first example of a space
with torsion (see \cite[p.~595]{Cartan1}). Consider $\R^3$ with its usual
Euclidean metric and the connection
\bdm
\nabla_X Y\ =\ \nabla^g_X Y - X\x Y,
\edm
corresponding, of course, to the choice $T=-2 \cdot e_1\wedge e_2\wedge e_3$.
Cartan observed correctly that this connection has same geodesics than
$\nabla^g$, but induces a different parallel transport\footnote{
``Deux tri\`edres [\ldots] de $\mathcal{E}$ seront parall\`eles lorsque les 
tri\`edres
correspondants de E [l'espace euclidien classique] pourront se d\'eduire
l'un de l'autre par un d\'eplacement h\'elico\"{\i}dal de pas donn\'e, de sens
donn\'e[\ldots]. L'espace $\mathcal{E}$ ainsi d\'efini admet un groupe de 
transformations \`a 6 param\`etres : ce serait notre espace ordinaire vu par
des observateurs dont toutes les perceptions seraient tordues.'' loc.cit.}.
Indeed, consider the $z$-axis $\gamma(t)=(0,0,t)$, a geodesic, and the
vector field $V$ which, in every point $\gamma(t)$, consists of the
vector $(\cos t,\sin t,0)$. Then one checks immediately that
$\nabla^g_{\dot{\gamma}}V= \dot{\gamma}\x V$, that is, the vector
$V$ is parallel transported according to a helicoidal movement.
If we now transport the vector along the edges of a closed triangle,
it will be rotated around three linearly independent axes, hence
the holonomy algebra is $\hol(\nabla)=\so(3)$.
\end{exa}
\begin{exa}[Holonomy of naturally reductive spaces]
%--------------------------------------------------
Consider a naturally reductive space $M^n=G/H$ as in Example \ref{exa-nat-red}
with its canonical connection $\nabla^1$, whose torsion 
is $T^1(X,Y):=-[X,Y]_{\m}$. Recall that $\ad: \h \ra\so(\m)$ denotes its
isotropy representation. The holonomy algebra
$\hol(\nabla^1)$ is the Lie subalgebra of $\ad(\h)\subset\so(\m)$ generated 
by the images under $\ad$ of all projections  of commutators
$[X,Y]_\h$  on $\h$ for $X,Y\in\m$,
\bdm
\hol(\nabla^1)\ =\ \mathrm{Lie}\big(\ad([X,Y]_\h)\big) \subset
\ad(\h)\subset\so(\m).
\edm
For all other connections $\nabla^t$ in this family, the general expression 
for the holonomy is considerably more complicated 
\cite[Thm. X.4.1]{Kobayashi&N2}.
\end{exa}
\begin{NB}[Holonomy \& contact properties]
%------------------------------------------
As we observed earlier, all contact structures are non-integrable
and therefore not covered by Berger's holonomy theorem. Via the cone
construction, it is nevertheless possible to characterize them
by a Riemannian holonomy property (see \cite{Boyer&G&M94}, \cite{BG}).
Consider a Riemannian manifold $(M^n,g)$ and its cone over the positive
real numbers $N:=\R_+\x M^n$ with the warped product metric 
$g_N:=dr^2+r^2g$. Then, $(M^n,g)$ is 
\begin{enumerate}
\item \emph{Sasakian} if and only if  
$\Hol(N;\nabla^g)\subset\U(\frac{n+1}{2})$, that is, its positive
cone is K\"ahler,
\item \emph{Einstein-Sasakian} if and only if
 $\Hol(N;\nabla^g)\subset\SU(\frac{n+1}{2})$, that is, its positive
cone is a (non-compact) Calabi-Yau manifold,
\item \emph{$3$-Sasakian} if and only if 
$\Hol(N;\nabla^g)\subset\Sympl(\frac{n+1}{4})$, that is, its positive
cone is hyper-K\"ahler.
\end{enumerate}
\end{NB} 
A holonomy group can be determined by computing  curvature.
\begin{thm}[{Ambrose-Singer, 1953 \cite{Ambrose&S53}}]
%-----------------------------------------------------
For any connection $\nabla$ on the tangent bundle of a Riemannian
manifold $(M,g)$, the Lie algebra $\hol(p)$ of $\Hol(p)$ in $p\in M$
is exactly the subalgebra of $\so(T_pM)$ generated by the elements
\bdm
P_\gamma^{-1}\circ \kr(P_\gamma V,P_\gamma W)\circ P_\gamma
\edm
where $V,W\in T_pM$ and $\gamma$ runs through all
piecewise smooth curves starting from $p$. 
\end{thm}
Yet, the practical use of this result is severely restricted by the fact
that the properties of the curvature transformation of a metric connection with
torsion are more complicated than  the Riemannian ones.
For example, $\kr(U,V)$ is still skew-adjoint with respect to the metric $g$,
\bdm
g(\kr(U,V)W_1,W_2)\ =\ -g(\kr(U,V)W_2,W_1),
\edm
but there  is  in general no relation between $g(\kr(U_1,U_2)W_1,W_2)$ and
$g(\kr(W_1,W_2)U_1,U_2)$ (but see Remark \ref{curv-par-torsion} below); in 
consequence, the Bianchi identities are less
tractable, the Ricci tensor is not necessarily symmetric etc. As an example
of these complications, 
we cite (see \cite{Ivanov&P01} for the case of skew-symmetric torsion):
\begin{thm}[First Bianchi identity]\label{Bianchi-I}
%---------------------------------------------------
\begin{enumerate}
\item[]
\item A metric connection $\nabla$ with vectorial torsion $V\in TM^n$ 
satisfies 
\bdm
\cyclic{X,Y,Z}\kr(X,Y)Z\ =\  \cyclic{X,Y,Z} dV(X,Y)Z.
\edm
\item A metric connection $\nabla$ with skew-symmetric torsion $T\in \Lambda^3(M^n)$ 
satisfies 
\bdm
\cyclic{X,Y,Z}\kr(X,Y,Z,V)\ =\  dT(X,Y,Z,V)-\sigma_T(X,Y,Z,V)
+(\nabla_V T)(X,Y,Z),
\edm
where $\sigma_T$ is a  $4$-form that is quadratic in $T$ defined by
$2\,\sigma_T=\sum\limits_{i=1}^n (e_i\haken T)\wedge (e_i\haken T)$ for any
orthonormal frame $e_1,\ldots, e_n$. 
\end{enumerate}
\end{thm}
\begin{NB}\label{curv-par-torsion}
%---------------------------------
Consequently, if the torsion $T\in \Lambda^3(M^n)$ of a metric connection with
skew-symmetric torsion happens to be $\nabla$-parallel, 
$\cyclic{X,Y,Z}\kr(X,Y,Z,V)$ is a $4$-form and thus antisymmetric. 
Since the cyclic sum over all four arguments of any $4$-form vanishes, 
we obtain 
\bdm
\cyclic{X,Y,Z,V} \big[\cyclic{X,Y,Z}\kr(X,Y,Z,V) \big]
\ =\ 2\, \kr(Z,X,Y,V) - 2\, \kr(Y,V,Z,X)\ =\ 0,
\edm 
as for the Levi-Civita connection. Thus the $\nabla$-curvature tensor is 
invariant under swaps of  the first and  second pairs of arguments.
\end{NB}
Extra care has to be taken when asking which
properties of the Riemannian holonomy group are preserved:
\begin{NB}[The holonomy representation may not be irreducible]
%------------------------------------------------------------------
%
In fact, there are many instances of irreducible manifolds with metric 
connections  whose holonomy representation is \emph{not} irreducible---for
example, the $7$-dimensional Aloff-Wallach space $N(1,1)=\SU(3)/S^1$ or the 
$5$-dimensional Stiefel manifolds $V_{4,2}=\SO(4)/\SO(2)$ (see \cite{BFGK}, 
\cite{Agri} and \cite{Agricola&F03a}).  This sheds some light on
why parallel objects are---sometimes---easier to find
for such connections. This also implies that
no analogue of de Rham's splitting theorem can hold (``A complete
simply connected Riemannian manifold with reducible holonomy representation
is a Riemannian product'').
\end{NB}
\begin{NB}[The holonomy group may not be closed]
%------------------------------------------------------------------
For the Riemannian restricted holonomy group, the argument goes as
follows: By de Rham's Theorem, one can assume that the holonomy group acts
irreducibly on each tangent space; but \emph{any} connected subgroup $G$ of
$\Orth(n)$ with this property has to be closed and hence compact. 
A counter-example for a torsion-free non-metric connection due to Ozeki can be 
found in  \cite[p.~290]{Kobayashi&N1}; similar (quite pathological) examples 
can be given for metric connections with torsion, although they seem to be
not
too interesting. It suffices to say that \emph{there is no theoretical argument
ensuring the closure of the restricted holonomy group of a metric
connection with torsion}.
\end{NB}
We are particularly interested in the vector bundle of $(r,s)$-tensors 
$T^{r,s}M$  over $M^n$, that of differential $k$-forms $\Lambda^k M$ and 
its spinor bundle $\Sigma M$ (assuming that $M$ is spin, of course). 
At some point $p\in M$, the fibers are  just $(T_pM)^r\ox (T^*_pM)^s$, 
$\Lambda^k T^*_pM $ or $\Delta_n$, the $n$-dimensional spin representation
(which has dimension $2^{[n/2]}$); an element of the fibre at some point
will be called an \emph{algebraic tensor, form, spinor} or just
\emph{algebraic vector} for short.
 Then, on any of these bundles,
\begin{enumerate}
\item the holonomy representation induces a representation of $\Hol(M;\nabla)$ 
on each fibre (the ``lifted holonomy representation'');
\item the metric connection $\nabla$ induces a connection (again denoted by
$\nabla$) on these vector bundles (the ``lifted connection'')
compatible with the induced metric (for tensors) or the induced
Hermitian scalar product (for spinors), it is thus again metric;
\item   in particular, there is a notion of ``lifted parallel transport''
consisting of isometries, and its abstract holonomy representation
on the fibres coincides with the lifted holonomy representation.
\end{enumerate}
we now formulate the general principle underlying our study.
\begin{thm}[General Holonomy Principle]\label{hol-principle}
%-----------------------------------------------------------
Let $M$ be a differentiable manifold and $E$ a  (real or complex)
vector bundle over $M$ endowed with (any!)
connection $\nabla$.  The following three properties are equivalent:
\begin{enumerate}
\item $E$ has a global section $\alpha$  invariant under parallel 
transport, i.\,e.~$\alpha(q)=P_\gamma(\alpha(p))$ for any path $\gamma$
from $p$ to $q$;
\item $E$ has a parallel global section $\alpha$, i.\,e.~$\nabla\alpha=0$;
\item at some point $p\in M$, there exists an algebraic vector 
$\alpha_0\in E_p$ which is invariant under the holonomy representation on 
the fibre.
\end{enumerate}
\end{thm}
\begin{proof}
%------------
The proof is almost philosophical. To begin, the first
and last conditions are equivalent: if $\alpha$ is a section
invariant under parallel transport $P_\gamma: E_p\ra E_q$, then, for a 
closed curve $\gamma$ through $p$, $\alpha(p)=P_\gamma(\alpha(p))$ and 
hence $\alpha$ is invariant under all holonomy transformations in $p$. 

Conversely,
let $\alpha_0\in E_p$ be a holonomy invariant algebraic vector. We then define
$\alpha(q):=P_\gamma(\alpha_0)$ for $q\in M$ and any  path $\gamma$
from $p$ to $q$. This definition is in fact path independent because,
 for any other path $\gamma'$ from $q$ to $p$, their concatenation is a 
closed loop, and $\alpha_0$ is, by assumption, invariant under parallel 
transport along closed curves. 

Finally, let $X$ be a vector field, $\gamma$ one of its integral curves
going from $p$ to $q$. Then obviously $\nabla_{\dot{\gamma}}\alpha=0$
is equivalent to $\alpha(q)=P_\gamma(\alpha(p))$, showing the
equivalence of (1) and (2).
\end{proof}
The following two consequences are immediate, but of the utmost importance.
\begin{cor}\label{cor-holonomy-principle}
%----------------------------------------
\begin{enumerate}
\item[]
\item The number of parallel global sections of $E$ coincides with
the number of trivial representations occuring in the  holonomy 
representation on the fibres. 
\item The holonomy group $\Hol(\nabla)$ is a subgroup of the
isotropy group $G_\alpha :=\{g\in\Orth(n)\,:\ g^*\alpha=\alpha \}$
of any parallel  global section $\alpha$ of $E$.
\end{enumerate}
\end{cor}
This is a powerful tool for (dis-)proving existence
of parallel objects. For example, the following is a well-known result
from linear algebra:
\begin{lem}
%----------
The determinant is an $\SO(n)$-invariant element in $\Lambda^n(\R^n)$
which is not $\Orth(n)$-invariant.
\end{lem}
\begin{figure}
\bdm
\includegraphics[width=7cm]{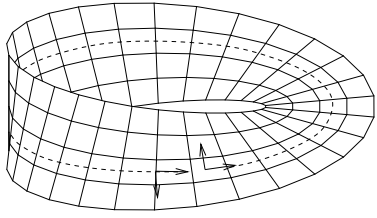}
\edm
\vspace{-4mm}
\caption{An orthonormal frame that is parallel transported along the
drawn curve reverses its orientation.}
\end{figure}
\begin{cor}
%----------
A Riemannian manifold $(M^n,g)$ is orientable if and only if the 
holonomy $\Hol(M;\nabla)$ of \emph{any metric connection} $\nabla$ is a 
subgroup of $\SO(n)$, and the volume form is then $\nabla$-parallel.
\end{cor}
\begin{proof}
%-------------
One knows that $(M^n,g)$ is orientable if and only if it admits a  nowhere
vanishing differential form $dM^n$ of degree $n$. Then pick
an orthonormal frame  
$e_1,\ldots, e_n$ in $p$ with dual $1$-forms $\sigma_1,\ldots,
\sigma_n$. Set $dM^n_p:=\sigma_1\wedge\ldots\wedge\sigma_n$ and
extend it to $M^n$ by parallel transport. Now everything follows from
the General Holonomy Principle.
\end{proof}
\begin{NB}
%---------
The property $\nabla dM^n=0$ for $ dM^n=\sigma_1\wedge\ldots\wedge\sigma_n$
can also be seen directly from the  formula
\bdm
\nabla_X (dM^n)(e_1,\ldots,e_n)\ =\ X(1) - \sum_{i=1}^n 
dM^n(e_1,\ldots,\nabla_Xe_i,\ldots,e_n),
\edm
since a metric connection satisfies $g(\nabla e_i,e_j)+ g(e_i,\nabla e_j)=0$,
so in particular $\nabla e_i$ has no $e_i$-component and all summands
on the right hand side vanish. 
\end{NB}
\begin{NB}
%---------
In fact, an arbitrary connection $\nabla$ admits a $\nabla$-parallel
$n$-form (possibly with zeroes) if and only if
\bdm
\sum_{i=1}^n g(\kr(U,V)e_i,e_i)\ =\ 0
\edm
for any orthonormal frame $e_1,\ldots, e_n$. This property is  weaker 
than the skew-adjointness of $\kr(U,V)$ that holds for all metric connections;
the holonomy is a subgroup of $\SL(n,\R)$. In 1924, J.\,A.~Schouten
called such connections ``inhaltstreue \"Ubertragungen'' 
(\emph{volume-preserving connections}), see \cite[p.\,89]{Schouten}. This terminology seems not to 
have been used  anymore afterwards\footnote{A \emph{D'Atri space}
is a Riemannian manifold whose local geodesic symmetries are
volume-preserving. Although  every naturally reductive space is a D'Atri 
space \cite{DAtri75}, the two notions are only loosely related.}.
\end{NB}
The existence of parallel objects imposes restrictions on the curvature
of the connection. For example, if a connection $\nabla$ admits a parallel
spinor $\psi$, we obtain by contracting the identity
\bdm
0 \ =\ \nabla\nabla\psi \ =\ \sum_{i,j=1}^n \kr^\nabla (e_i,e_j) 
e_i\cdot e_j \cdot \psi
\edm
the following integrability condition:
\begin{prop}\label{curv-parallel}
%---------------------------------
Let $(M^n,g)$ be a Riemannian spin manifold, $\nabla$ a metric connection
with torsion $T\in\Lambda^3(M^n)$.
A $\nabla$-parallel spinor $\psi$ satisfies
\bdm
\left[\frac{1}{2} X\haken dT+\nabla_X T \right]\cdot\psi\ =\ 
\Ric^\nabla(X)\cdot\psi
\edm
In particular, the existence of a $\nabla^g$-parallel spinor $(T=0)$
implies Ricci-flatness.
\end{prop}
\smallskip
Before considering general metric connections with torsion on manifolds,
it is worthwile to investigate the flat case $\R^n$ endowed with
its standard Euclidean metric and metric connections with constant
torsion, for it exhibits already some characteristic features of the more
general situation. Unless otherwise stated, these results can be found in
\cite{Agricola&F03a}.

The exterior algebra $\Lambda (\R^n)$ and the Clifford algebra
$\mathrm{Cl}(\R^n)$ are---as vector spaces---equivalent 
$\SO(n)$-representations, and they both act on the  complex vector space 
$\Delta_n$ of  $n$-dimensional spinors.  The Clifford 
algebra is an associative algebra with an underlying 
Lie algebra structure,
\bdm
[ \alpha, \beta ] \ = \ \alpha \cdot \beta \, - \, \beta \cdot \alpha , 
\quad \alpha, \beta \in \mathrm{Cl}(\R^n) .
\edm
We denote the corresponding Lie algebra by $\mathfrak{cl}(\R^n)$. The 
Lie algebra $\so(n)$ of the special orthogonal group is a subalgebra
of $\mathfrak{cl}(\R^n)$, 
\bdm
\so(n) \ = \ \mathrm{Lin} \big\{X \cdot Y  : X,
Y \in \R^n \ \text{and} 
\  \langle X  , \, Y \rangle \ = \ 0 \big\} \, \subset \, 
\mathfrak{cl}(\R^n) \, .
\edm
Consider an algebraic $k$-form $T \in \Lambda^k(\R^n)$ and 
denote by $G_T$ the group of all orthogonal transformations of $\R^n$ 
preserving the form $T$, by $\g_T$ its Lie algebra.
As described in Example \ref{exa-higher-forms}, we consider the spin
connection  acting on spinor 
fields $\psi : \R^n \rightarrow \Delta_n$ by the formula
\bdm
\nabla_X \psi \ := \ \nabla^{g}_X \psi \, + \, \frac{1}{2}  
(X \haken \T) \cdot \psi .
\edm 
For a $3$-form $T \in \Lambda^3(\R^n)$, the spinorial covariant 
derivative $\nabla$ is induced by the linear metric 
connection with torsion tensor $T$.
 For a general exterior form $T$, we introduce
a new Lie algebra $\g_T^*$ that is a subalgebra of
$\mathfrak{cl}(\R^n)$.
\begin{dfn}
%---------
Let $T$ be an exterior form on $\R^n$. The Lie algebra $\g_T^*$ is the 
subalgebra of $\mathfrak{cl}(\R^n)$ generated by all elements 
$X \haken T$, where $X \in \R^n$ is a vector.
\end{dfn}
\noindent
The Lie algebra $\g_T^*$ is invariant under
the action of the isotropy group $G_T$. The derived
algebra $\big[\g_{T}^*  ,  \g_{T}^* \big]$ is the Lie algebra generated 
by all curvature transformations
of the spinorial connection $\nabla$. It is the Lie algebra 
of the infinitesimal holonomy group of the spinorial covariant 
derivative $\nabla^{\T}$ (see \cite[Ch. II.10]{Kobayashi&N1}):
\begin{dfn} 
%----------
Let $T$ be an exterior form on $\R^n$. The Lie algebra
$\h_{T}^*  := \big[\g_{T}^*,\g_{T}^*\big]  \subset  \mathfrak{cl}(\R^n)$
is called the \emph{infinitesimal holonomy algebra} of the exterior form $T$.
It is invariant under the action of the isotropy group $G_{T}$.
\end{dfn} 
For a $3$-form $T$, the Lie algebras $\g_{T}^*,\h_{T}^*\subset\mathfrak{so}(n)$
are subalgebras of the orthogonal Lie algebra, reflecting the fact that the 
corresponding spinor derivative  is induced by a linear metric
connection. In fact, this result still holds for $k$-forms satisfying
$k + \binom{k-1}{2} \equiv 0 \ \mathrm{mod}\, 2$. Furthermore, the
General Holonomy Principle (Theorem \ref{hol-principle}) implies:
\begin{prop} 
%-----------
There exists a non-trivial $\nabla$-parallel spinor field $\psi : \R^n 
\rightarrow \Delta_n$ if and only if there exists a constant spinor $\psi_0 \in \Delta_n$ such that
$\mathfrak{h}_{\T}^* \cdot \psi_0 =  0$. In particular,
any $\nabla$-parallel spinor field is constant 
for a perfect Lie algebra $\g_{T}^* \,( \g_{T}^*=\h_T^*)$.
\end{prop} 
\begin{exa} 
%----------
Any $2$-form $\T \in \Lambda^2(\R^n)$ of rank $2k$ is
equivalent to $A_1 \cdot e_{12} + \cdots + A_k \cdot e_{2k-1,2k}$. The 
Lie algebra $\mathfrak{g}_{\T}^*$ is generated by
the elements $e_1, e_2, \cdots , e_{2k-1}, e_{2k}$. It is isomorphic
to the Lie algebra $\spin(2k+1)$. In particular, if $n = 8$ then 
$\Delta_8 = \R^{16}$ is a real, $16$-dimensional and the spinorial holonomy
algebra of a generic $2$-form in eight variables is the unique $16$-dimensional
irreducible representation of $\spin(9)$.
\end{exa}
\begin{exa}
%---------
Consider the $4$-forms $T_1:=e_{1234}\in \Lambda^4(\R^{n})$ for $n\geq 4$
and $T_2 = e_{1234} + e_{3456} \in \Lambda^4(\R^{m})$ for $m\geq 6$. A
straightforward 
computation yields that $\g_{T_1}^*$ and $\g_{T_2}^*$ are isomorphic  to the
pseudo-orthogonal Lie algebra $\so(4,1)$ embedded in a non-standard way
and the Euclidean Lie algebra $\mathfrak{e}(6)$, respectively. 
\end{exa} 
\begin{exa}
%---------
Consider the volume form $T = e_{123456}$ in $\R^n$ for $n\geq 6$.
The subalgebra $\g_{T}^*$ of $\mathrm{Cl}(\R^n)$
is isomorphic to the compact Lie algebra $\spin(7)$.
\end{exa} 
If $T$ is a $3$-form, more can be said. For example,
$\g_{T}^*$ is always semisimple and the following proposition 
shows that  it cannot be contained in the unitary Lie algebra 
$\mathfrak{u}(k) \subset \mathfrak{so}(2k)$.
This latter result is in sharp contrast to the situation on arbitrary 
manifolds, where such $3$-forms occur for almost Hermitian structures.
\begin{prop}\label{invform}  
%--------------------------
Let $\T$ be a $3$-form in $\R^{2k}$ and suppose that there
exists a $2$-form $\Omega$ such that $\Omega^{k} \neq 0$ and
$[\g_{T}^*,  \Omega  ] = 0$. Then $T$ is zero, $T = 0$.
\end{prop} 
Moreover, only constant spinors are parallel:
\begin{thm}\label{pareucl} 
%-------------------------
Let $T \in \Lambda^3(\R^n)$ be a $3$-form. If there exists a non trivial 
spinor $\psi \in \Delta_n$ such that $\g_{T}^* \cdot \psi = 0$, 
then $T = 0$. In particular, $\nabla$-parallel spinor fields are 
$\nabla^g$-parallel and thus constant. 
\end{thm} 
\begin{proof}
%------------
The proof  is remarkable in as much as it is of purely
algebraic nature. Indeed, it is a consequence of the 
following  formulas concerning the action of exterior forms of different 
degrees on spinors (see Appendix \ref{formulas} for the first, the other
two are simple calculations in the Clifford algebra):
\bdm
2\sigma_T \,=\, - T^2 +||T||^2 ,\quad
2\sigma_T \cdot\psi \,=\, \sum_{k=1}^n (e_k\haken T)\cdot
(e_k\haken T)  +  3  ||T||^2 ,  \quad 3  T \, = \, 
\sum_{k=1}^n e_i \cdot (e_i \haken T ) . 
\edm
For they imply that a  spinor $\psi$ with $(e_i\haken T)\cdot\psi=0$
for all $e_i$ has to satisfy $||T||^2 \psi=0$, so $T$ must vanish.
\end{proof}
This result also applies to flat tori $\R^n/\Z^n$, as the torsion form $T$ 
is assumed to be constant. Later, we shall prove a suitable generalization 
on compact spin manifolds with $\Scal^g\leq 0$, see 
Theorem~\ref{parallelflach2}. 
%
%
%-----------------------------------------------------------------------------
\section{Geometric stabilizers}\label{geom-stab}
%-----------------------------------------------------------------------------
%
By the General Holonomy Principle,  geometric representations with invariant 
objects are a natural source for parallel objects. This leads to the 
systematic investigation of \emph{geometric stabilizers}, which we shall 
now discuss.
\subsection{$\U(n)$ and $\SU(n)$ in dimension $2n$}\label{stab-Un}
%-----------------------------------------------------------------
A Hermitian metric $h(V,W)= g(V,W)-ig(JV,W)$ is invariant under
$A\in \End(\R^{2n})$ if and only if $A$ preserves the Riemannian
metric $g$ and 
the K\"ahler form $\Omega(V,W):= g(JV,W)$. Thus $\U(n)$ is  embedded
in $\SO(2n)$ as
\bdm
\U(n)\ =\ \{A\in \SO(2n)\ |\ A^*\Omega=\Omega\} .
\edm
To fix ideas, choose a skew-symmetric endomorphism $J$ of $\R^{2n}$ with
square $-1$ in the normal form
\bdm
J\ =\ \diag \big( j,j,j,\ldots\big) \text{ with } j\ =\ 
\left[\begin{array}{cc} 0 & -1 \\ 1 & 0 \end{array}\right].
\edm
Then a complex $(n\x n)$-matrix $A=(a_{ij})\in \U(n)$ is realized
as a real  $(2n\x 2n)$-matrix with $(2\x 2)$-blocks  
$\left[\begin{array}{rr} \Re a_{ij} & -\Im a_{ij} \\ \Im a_{ij}& \Re a_{ij} 
\end{array}\right]$.
An adapted orthonormal frame is one such that $J$ has the given normal form;
$\U(n)$ consists then exactly of those endomorphisms transforming
adapted orthonormal frames  into adapted orthonormal frames.
Allowing now complex coefficients, one obtains an $(n,0)$-form $\Psi$
by declaring
\bdm
\Psi \ :=\ (e_1+ie_2)\wedge\ldots\wedge (e_{n-1}+ i e_{2n})\ =:\ 
\Psi^+ + i\Psi^-
\edm 
in the adapted frame above. An element $A\in \U(n)$ acts on   $\Psi$
by multiplication with $\det A$.
\begin{lem}
%----------
Under the restricted action of $\U(n)$,  $\Lambda^{2k}(\R^{2n})$ 
contains the trivial  representation once; it is generated by
$\Omega,\Omega^2,\ldots,\Omega^n$.
\end{lem}
The action of  $\U(n)\subset \Orth(2n)$ cannot be lifted
to an action of $\U(n)$ inside $\Spin(2n)$---reflecting the fact that
not every K\"ahler manifold is spin. For the following arguments though,
it is enough to consider $\un(n)$ inside $\spin(2n)$. 
It then appears that $\un(n)$ has no invariant spinors, basically because 
$\un(n)$ has a one-dimensional center, generated precisely by $\Omega$
after identifying  $\Lambda^2(\R^{2n})\cong \so(2n)$. Hence
one-dimensional  $\un(n)$-representations are usually not trivial.
More precisely, the complex $2n$-dimensional spin representation 
$\Delta_{2n}$ splits into
two irreducible components $\Delta^\pm_{2n}$ described in terms of 
eigenspaces of $\Omega\in \un(n)$.
%
%\bdm
%\Delta^+_{2n}\big|_{\U(n)}\ \cong \ \Lambda^n \C^n\oplus
%\Lambda^{n-2} \C^n\oplus\ldots ,\quad
%\Delta^-_{2n}\big|_{\U(n)}\ \cong \ \Lambda^{n-1} \C^n\oplus
%\Lambda^{n-3} \C^n\oplus\ldots.
%\edm
%
Set (see \cite{Dirac-Buch-00} and \cite{Kirchberg86} for details on this
decomposition of spinors)
\bdm
S_r\ =\ \{\psi \in \Delta_{2n}\ :\ \Omega\psi=i(n-2r)\psi \}, \quad
\dim S_r\ =\ \binom{n}{r}, \quad 0\leq r\leq n.
\edm
$S_r$ is isomorphic to the space of $(0,r)$-forms with values in $S_0$ 
(which explains the dimension),
\bdm
S_r\ \cong\ \overline{\Lambda^{0,r}}\ox S_0.
\edm
Since the spin representations decompose as
\bdm
\Delta^+_{2n}\big|_{\un(n)}\ \cong \ S_{n} \oplus S_{n-2} \oplus\ldots ,\quad
\Delta^-_{2n}\big|_{\un(n)}\ \cong \ S_{n-1}\oplus S_{n-3} \oplus\ldots 
\edm
we conclude immediately that they cannot contain a trivial 
$\un(n)$-representation for $n$ odd. For $n=2k$ even, $\Omega$ has
eigenvalue zero on $S_k$, but this space is an irreducible representation
of dimension $\binom{2k}{k}\neq 1$, hence not trivial either. 
The representations $S_0$ and $ S_n$ are one-dimensional,
but again not trivial under $\un(n)$. If one restricts further to
$\su(n)$, they are indeed:
\begin{lem}
%----------
The spin representations $\Delta^\pm_{2n}$ contain no $\un(n)$-invariant
spinor. If one restricts further to $\su(n)$, there are exactly two
invariant spinors (both in $\Delta^+_{2n}$ for $n$ even, one in each 
$\Delta^\pm_{2n}$ for $n$ odd). 
\end{lem}
All other spinors in $\Delta^\pm_{2n}$ have geometric stabilizer groups
that do not act irreducibly on the tangent representation $\R^{2n}$.
They can be described explicitly in a similar way; By  de Rham's splitting 
theorem, they do not appear in the Riemannian setting.

To finish, we observe that the almost complex structure $J$ (and hence 
$\Omega$) can be
recovered from the invariant spinor $\psi^+\in \Delta^+_{2n}$ by
$J(X)\psi^+:=iX\cdot\psi^+$ ($X\in TM$), a formula well known from the 
investigation of Killing spinors on $6$-dimensional nearly K\"ahler 
manifolds (see \cite{Grunewald90} and \cite[Section 5.2]{BFGK}). 
\begin{NB}\label{embedding}
%---------------------------
The discussion of geometric stabilizers would  not be
complete without the explicit realization of these subalgebras inside
$\so(n)$ or $\spin(n)$. We illustrate this by describing $\un(n)$
inside $\so(2n)$.  Writing elements  $\omega\in \so(2n)$ as 
$2$-forms with respect to some orthonormal and $J$-adapted basis, 
$\omega=\sum \omega_{ij}e_i\wedge e_j$ for $1\leq i < j\leq 2n $, the 
defining equations for $\un(n)$
inside $\so(2n)$ translate into the conditions
\bdm
\omega_{2i-1,2j-1}\ =\ \omega_{2i,2j}\ \text{ and }\
\omega_{2i-1,2j}\ =\ - \omega_{2i,2j-1}\ \text{ for }\  
1\,\leq\, i\, <\, j\,\leq n.
\edm   
The additional equation picking out $\su(n)\subset\un(n)$ is
\bdm
\omega_{12}+\omega_{34}+\ldots+\omega_{2n-1,2n}\ =\ 0.
\edm
Of course, the equations get more involved for  complicated embeddings
of higher codimension (see for example \cite{Agricola&F03a} for the 
$36$-dimensional $\spin(9)$ inside the $120$-dimensional $\so(16)$), but 
they can easily be mastered 
with the help of standard linear algebra computer packages.
\end{NB}
\begin{NB}
%---------
The group $\Sympl(n)\subset\SO(4n)$ can be deduced from the previous 
discussion:
$\Sympl(n)$ with quaternionic entries $a+bj$ is embedded into $\SU(2n)$ by 
$(2\x 2)$-blocks  
$\left[\begin{array}{rr} a & b \\ -\bar{b}& \bar{a} \end{array}\right]$,
and  $\SU(2n)$ sits in $\SO(4n)$ as before. We shall not treat
$\Sympl(n)$- and quaternionic geometries in this expository article 
(but see \cite{Swann89}, \cite{Swann91}, \cite{Alekseevsky&M&P98}, 
\cite{Grantcharov&P00}, \cite{Alexandrov03}, \cite{Poon&S03}, 
\cite{Alekseevky&C05}, \cite{Martin-Cabrera&S04} for a first acquaintance).
\end{NB}
\subsection{$\U(n)$ and $\SU(n)$ in dimension $2n+1$}
%---------------------------------------------------
These $G$-structures arise from contact structures and are 
remarkable inasmuch they manifest a genuinely non-integrable behaviour---they
do not occur in Berger's list because the action of $\U(n)$ on $\R^{2n+1}$
is not irreducible, hence any manifold with this action as Riemannian 
holonomy representation
splits by de Rham's theorem. Given an almost contact metric manifold
$(M^{2n+1}, g,\xi,\eta,\vphi)$, we may construct an adapted local
orthonormal frame by choosing any $e_1\in\xi^\perp$ and setting 
$e_2=\vphi(e_1)$ (as well as fixing $e_{2n+1}=\xi$ once and for all); 
now choose any $e_3$ perpendicular to $e_1,e_2, e_{2n+1}$ and set
again $e_4=\vphi(e_3)$ etc. With respect to such a basis, $\vphi$ is given
by 
\bdm
\vphi\ =\ \diag(j,j,\ldots,j,0)\ \text{ with } j\ =\ 
\left[\begin{array}{cc} 0 & -1 \\ 1 & 0 \end{array}\right]
\edm
and we conclude that the structural group of $M^{2n+1}$ is reducible
to $\U(n)\x \{1\}$. If we denote the fundamental form by $F$ (see Example 
\ref{exa-contact}), then
\bdm
\U(n)\x \{1\}\ =\  \{A\in \SO(2n+1)\ |\ A^* F=F\} .
\edm
The $\U(n)\x \{1\}$-action on $\R^{2n+1}$ inherits invariants
from the $\U(n)$-action on $\R^{2n}$ in a canonical way;
one then just needs to check that no new one appears. Hence, we can conclude: 
\begin{lem}
%----------
Under the  action of $\U(n)$,  $\Lambda^{2k}(\R^{2n+1})$ 
contains the trivial  representation once; it is generated by
$F,F^2,\ldots,F^n$.
\end{lem}
The action of  $\U(n)\subset \Orth(2n+1)$ cannot, in general,  be lifted
to an action of $\U(n)$ inside $\Spin(2n+1)$. As in the almost Hermitian
case, let's thus study the $\un(n)$ action on $\Delta_{2n+1}$. The
irreducible $\Spin(2n+1)$-module $\Delta_{2n+1}$ splits into
$\Delta_{2n}^+\oplus \Delta_{2n}^-$ under the restricted action of
$\Spin(2n)$, and it decomposes accordingly into
\bdm
\Delta_{2n+1}\ =\ S_0\oplus\ldots S_{n},\quad
S_r\ =\ \{\psi \in \Delta_{2n+1}\ :\ F\psi=i(n-2r)\psi \}, \quad
\dim S_r\ =\ \binom{n}{r}, \quad 0\leq r\leq n.
\edm
Hence,  $\Delta_{2n+1}$ can be identified  with 
$\Delta_{2n}^+\oplus \Delta_{2n}^-$, yielding finally the following
result:
\begin{lem}
%----------
The spin representation $\Delta_{2n+1}$ contains no $\un(n)$-invariant
spinor. If one restricts further to $\su(n)$, there are precisely two
invariant spinors.
\end{lem}
\subsection{$G_2$ in dimension $7$}
%-----------------------------------
While invariant $2$-forms exist in all even dimensions
and lead  to the rich variety of almost Hermitian structures,
the geometry of $3$-forms played a rather exotic role in classical
Riemannian  geometry until the nineties, as it occurs only in apparently
random dimensions,  most notably dimension seven. That $G_2$ is the relevant
simple Lie group is a classical, although unfortunately not so well-known
result from invariant theory. A mere dimension count shows
already this effect (see Table \ref{dim-G2}):
%
%\bigskip
\begin{table}
\begin{center}
\setlength{\extrarowheight}{4pt}
\begin{tabular}{|c|c|c|}\hline
$n$ & $\dim \GL(n,\R) - \dim \Lambda^3\R^n$ & $\dim \SO(n)$ \\ \hline\hline
3 & $9-1=8$ & 3\\ \hline
4 & $16-4=12$ & 6 \\ \hline
5 & $25-10=15$ & 10 \\ \hline
6 & $36-20=16$ & 15\\ \hline
7 & $49-35=14$ & 21\\ \hline
8 & $64-56=8$ & 28 \\ \hline
\end{tabular}
\bigskip
\end{center}
\caption{Dimension count for possible geometries defined by $3$-forms.}
\label{dim-G2}
\end{table}
the stabilizer of a generic $3$-form $\omega^3$
\bdm
G^n_{\omega^3}\ :=\ \{ A\in \GL(n,\R)\ |\ \omega^3=A^*\omega^3\}
\edm
cannot be contained in the orthogonal group  for $n\leq 6$,
it must lie in some group between $\SO(n)$ and $\SL(n,\R)$ (for $n=3$, we
even have $G^3_{\omega^3}=\SL(3,\R)$). The case
 $n=7$ is the first dimension where $G^n_{\omega^3}$
can sit in $\SO(n)$. That this is indeed the case
was shown as early as  1907 in the doctoral dissertation of Walter
Reichel in Greifswald, supervised by F.~Engel (\cite{Reichel07}).
More precisely, he computed a system of invariants for a $3$-form in
seven variables and showed that there are exactly two open 
$\GL(7,\R)$-orbits of $3$-forms. The stabilizers of 
any representatives $\omega^3$ and $\tilde{\omega}^3$ of these
orbits are $14$-dimensional simple Lie groups of rank two,
one  compact and the other non-compact: 
\bdm
G^7_{\omega^3}\ \cong\ G_2\ \subset\ \SO(7),\quad
G^7_{\tilde{\omega}^3}\ \cong\ G^*_2\ \subset\ \SO(3,4).
\edm
Reichel also showed the corresponding embeddings of Lie algebras
by explicitly writing down seven equations for
 the coefficients of $\so(7)$ resp.~$\so(3,4)$ (see Remark \ref{g2-in-so7}).
As in the case of almost Hermitian geometry, every author has his or
her favourite normal $3$-form with isotropy group $G_2$, 
for example,
\bdm
\omega^3\ :=\ e_{127}+e_{347}-e_{567}+e_{135}-e_{245}+e_{146}+e_{236}.
\edm
An element of the second orbit with stabilizer the split form $G^*_2$ of
$G_2$ may be obtained by reversing any of the signs in $\omega^3$.
\begin{lem}
%----------
Under $G_2$, one has the decomposition
$\Lambda^3(\R^7)\cong \R\oplus \R^7\oplus S_0(\R^7)$, where
$\R^7$ denotes the $7$-dimensional standard representation
given by the embedding $G_2\subset\SO(7)$ and  $S_0(\R^7)$ denotes
the traceless symmetric endomorphisms of $\R^7$ (of dimension $27$).
\end{lem}
Now let's consider the spinorial picture, as $G_2$ can indeed be lifted 
to a subgroup of $\Spin(7)$. From a purely representation theoretic 
point of view,
this case is trivial: $\dim \Delta_7=8$ and the only irreducible
representations of $G_2$ of dimension $\leq 8$ are the trivial
and the $7$-dimensional representations. Hence $8=1+7$ yields:
\begin{lem}\label{G2-on-spinors}
%--------------------------------
Under the restricted action of $G_2$, the $7$-dimensional spin representation 
$\Delta_7$ decomposes as $\Delta_7\cong \R\oplus \R^7$. 
\end{lem}
This Lemma has an important consequence: the `spinorial'
characterization of $G_2$-manifolds.
\begin{cor}\label{G2-equiv-par-spin}
%-----------------------------------
Let $(M^7,g)$ be a Riemannian manifold, $\nabla$ a metric connection on its
spin bundle. Then there exists a $\nabla$-parallel spinor if and only if
$\Hol(\nabla)\subset G_2$.  
\end{cor}
One direction follows from the fact that $G_2$ is the stabilizer of an 
algebraic spinor, the converse from Lemma \ref{G2-on-spinors}.
 
In fact, the invariant $3$-form and the invariant algebraic spinor
$\psi$ are equivalent data. They are related (modulo an irrelevant constant)
by 
\be\label{spinor-form}
\omega^3(X,Y,Z)\ =\ \langle X\cdot Y\cdot Z\cdot \psi,\psi\rangle .
\ee
We now want to ask which subgroups $G\subset G_2$ admit other invariant
algebraic spinors. Such a subgroup has to appear on Berger's list and its
induced action on $\R^7$ (viewed as a subspace of $\Delta_7$) has to contain
one or more copies of the trivial representation. Thus, the only
possibilities are $\SU(3)$ with $\R^7\cong\R\oplus\C^3$ (standard 
$\SU(3)$-action on $\C^3$) and $\SU(2)$ 
with $\R^7\cong 3\cdot\R\oplus\C^2$ (standard $\SU(2)$-action on $\C^2$). 
Both indeed occur, with a total of $2$ resp.~$4$ invariant spinors.  
\begin{NB}
%---------
A modern account of Reichel's results can be found in the article
\cite{Westwick81} by R.~Westwick; it is interesting (although it seems not to
have had any further influence) that J.\,A.~Schouten also rediscovered these
results in 1931 \cite{Schouten31}. A classification of $3$-forms
is still possible in dimensions $8$ (\cite{Gurevich35}, \cite{Gurevich35}, 
\cite{Djokovic83}) and $9$ 
(\cite{Vinberg&E88}), although the latter one is already of inexorable 
complexity. Based on these results, J.~Bure\v{s} and J.~Van\v{z}ura started
recently the investigation of so-called \emph{multisymplectic structures}
(\cite{Vanzura01}, \cite{Bures&V03},  \cite{Bures04}).
\end{NB}
\begin{NB}\label{g2-in-so7}
%--------------------------
$\g_2$ inside $\spin(7)$ is a good example for illustrating how to 
obtain the defining equations of stabilizer subalgebras with the aid
of the computer  (see Remark \ref{embedding}); 
one has just to be aware that they depend not only on the orthonormal basis 
but also on a choice of spin representation. To this aim, fix a realization 
of the spin
representation $\Delta_n$ and a representative $\psi$ of the orbit of spinors 
whose stabilizer is the group $G$ we are interested in. As usual, identify 
the Lie  algebra $\spin(n)$ with the elements of the form 
$\omega=\sum_{i<j}\omega_{ij} e_i\cdot e_j$ inside the Clifford algebra
$\mathrm{Cl}(n)$. Replacing $e_i, e_j$ by the chosen representative
matrices, $\omega\cdot\psi=0$ is equivalent to a set of
equations for the coefficients $\omega_{ij}$; see for example 
\cite[p. 261]{FKMS} for an explicit realization of $\g_2\subset\spin(7)$. 
\end{NB}
\subsection{$\Spin(7)$ in dimension $8$}\label{stab-spin7}
%---------------------------------------------------------
As we just learned from $G_2$ geometry, $\Spin(7)$ has an irreducible
$8$-dimensional
representation isomorphic to $\Delta_7$, hence it can be viewed as a 
subgroup of $\SO(8)$, and it does lift to   $\Spin(8)$.
By restricting to $\SO(7)$, $\Spin(7)$ certainly also has a 
$7$-dimensional representation. What is so special in this dimension
is that $\Spin(7)$ has \emph{two} conjugacy classes in $\SO(8)$ that
are interchanged by means of the triality automorphism; hence
the decomposition of the spin representation depends on the (arbitrary)
choice of one of these classes.   
\begin{lem}\label{spin7}
%-----------------------
Under the restricted action of $\Spin(7)$, the $8$-dimensional spin 
representations decompose as $\Delta_8^+\cong \R^8\cong \Delta_7$
and $\Delta_8^-\cong \R\oplus\R^7$ for one choice of 
$\Spin(7)\subset\SO(8)$; the other choice swaps $\Delta_8^+$ and 
$\Delta_8^-$.  
\end{lem}
In particular, there is exactly one invariant spinor in $\Delta^+_8$. 
Again,  it corresponds one-to-one to an invariant form, of degree
$4$ in this case:
\bdm
\beta^4(X,Y,Z,V)\ =\ \langle X\cdot Y\cdot Z\cdot V\cdot\psi,\psi\rangle .
\edm
Yet, $\Spin(7)$-geometry in dimension eight is not just an enhanced
version of $G_2$-geometry in dimension seven. Because
$\dim \GL(8,\R)=64< 70=\dim \Lambda^4(\R^8)$, 
there are no dense open orbits under the action of $\GL(8,\R)$. Thus,
there is no result in invariant theory similar to that of Reichel for $G_2$
in the background.

Let's fix the first choice for embedding $\Spin(7)$ in $ \SO(8)$ made in 
Lemma \ref{spin7}. A second invariant spinor can either be in $\Delta_8^+$
or in $\Delta_8^-$. If it is in $\Delta_8^+$, we are asking for a subgroup
$G\subset \Spin(7)$ whose action on $\R^7$ contains the
trivial representation once---obviously, $G_2$ is such a group.
Under $G_2$,  $\Delta_8^+$ and $\Delta_8^-$ are isomorphic, 
\bdm
\Delta_8^+\big|_{G_2}\ \cong\ \Delta_8^-\big|_{G_2}\ \cong\ \R\oplus \R^7.
\edm
Thus, $\SU(3)\subset G_2\subset \Spin(7)$ and $\SU(2)\subset G_2\subset
\Spin(7)$ are two further admissible groups with $2+2$ and $4+4$ 
invariant spinors. On the other hand, if we impose a second invariant
spinor to live in $\Delta_8^-$, we need a subgroup $G\subset\Spin(7)$ that
has partially trivial action on $\R^7$, but not on  $\Delta_8^+\cong\R^8$.
A straightforward candidate is $G=\Spin(6)$ with its standard embedding
and $\R^7=\R^6\oplus\R$; the classical isomorphism $\Spin(6)=\SU(4)$
shows that $G$ acts irreducibly on $\Delta_8^+\cong\C^4$. The group $\SU(4)$
in turn has subgroups $\Sympl(2)=\Spin(5)$ and $\SU(2)\times\SU(2)=\Spin(4)$ 
that still act irreducibly on $\Delta_8^+$, and act on $\Delta_8^-$ by
\bdm
\Delta_8^-\big|_{\Sympl(2)}\ \cong\ 3\cdot\R\oplus\R^5,\quad
\Delta_8^-\big|_{\SU(2)\times\SU(2)}\ \cong\ 4\cdot\R\oplus\R^4.
\edm
The results are summarized in Table \ref{par-spin-dim8}. The resemblance
between Tables \ref{par-spin-dim8} and \ref{par-spin-dim10} in Section
\ref{physmot} is no coincidence. A convenient way to describe $\Spin(9,1)$
is to start with $\Spin(10)$ generated by elements $e_1,\ldots, e_{10}$
and acting irreducibly on $\Delta^+_{10}$. The
vector spaces $\Delta^+_{10}$ and $\Delta^+_{9,1}$ can be identified, and
$\Spin(9,1)$ can be generated by $e_i^*:=e_i$ for $i=1,\ldots,9$ and
$e_{10}^*:=i\,e_{10}$. Elements $\omega\in \spin(9,1)$ can thus be written
\bdm
\omega\ =\ \sum_{1 \leq i<j\leq 10}\omega_{ij}\, \omega_{ij} e^*_i\wedge e^*_j
\ =\ \sum_{i<j\leq 9}\omega_{ij}\, e_i\wedge e_j +
i \sum_{k\leq 9} \omega_{k,10}\, e_k\wedge e_{10}
\edm
and we conclude that $\spin(9,1)$ can  be identified with 
$\spin(9)\ltimes\R^9$.
A spinor $\psi\in\Delta^+_{9,1}$ is stabilized by
an element $\omega\in\spin(9,1)$ if and 
only if
\bdm
0\ =\ \sum_{1 \leq i<j\leq 10}\omega_{ij}\, e^*_i\cdot e^*_j\cdot\psi\ =\ 
\sum_{i<j\leq 9}\omega_{ij}\, e_i\cdot e_j\cdot\psi + i\sum_{k\leq 9}
\omega_{k,10}\, e_k\cdot e_{10}\cdot\psi.
\edm
In this last expression, both real and imaginary part have to 
vanish simultaneously, leading to 16 equations.  A careful look
reveals that they define $\spin(7)\ltimes\R^8$, and the statements from 
Table \ref{par-spin-dim8} imply those from Table \ref{par-spin-dim10}
since $\Delta^+_{9,1}\cong \Delta^+_8\oplus\Delta^-_8$.
\begin{table}%\label{par-spin-dim8}
\begin{center}
\setlength{\extrarowheight}{4pt}
\begin{tabular}{|c|c|c|}\hline
group & inv. spinors in $\Delta_8^+$ & inv. spinors in $\Delta_8^-$ \\ \hline\hline
$\Spin(7)$ & 0 & 1\\ \hline
$\SU(4)\cong\Spin(6)$ & 0 & 2 \\ \hline
$\Sympl(2)\cong\Spin(5)$ & 0 & 3 \\ \hline
$\SU(2)\times\SU(2)\cong\Spin(4)$ & 0 & 4 \\ \hline
$G_2$ & 1 & 1\\ \hline
$\SU(3)$ & 2 & 2\\ \hline
$\SU(2)$ & 4 & 4\\ \hline
$\{e\}$ & 8 & 8\\ \hline
\end{tabular}
\bigskip
\end{center}
\caption{Possible stabilizers of invariant spinors in dimension $8$.}
\label{par-spin-dim8}
\end{table}
\begin{NB}[Weak $\mathrm{PSU}(3)$-structures]
%---------------------------------------------
Recently, Hitchin observed that $3$-forms can be of interest in
$8$-dimensional geometry as well (\cite{Hitchin01}). The canonical 
$3$-form on the Lie algebra $\su(3)$ spans an open orbit under
$\GL(8,\R)$, and the corresponding $3$-form on $\SU(3)$ is parallel
with respect to the Levi-Civita connection of the biinvariant metric. The 
Riemannian holonomy
reduces to $\SU(3)/\Z_2=:\mathrm{PSU}(3)$. More generally, manifolds
modelled on this group lead to the investigation of closed and coclosed
$3$-forms that are not parallel, see also \cite{Witt06}.
\end{NB}
Sorting out the technicalities that we purposely avoided, one obtains
Wang's classification of Riemannian parallel spinors. By de Rham's
theorem, only irreducible holonomy representations occur for the
Levi-Civita connection. From Proposition \ref{curv-parallel},
we already know that these manifolds are Ricci-flat.
\begin{thm}[Wang's Theorem, {\cite{Wang89}}]\label{thm-wang}
%-----------------------------------------------------------
Let $(M^n,g)$ be a complete, simply connected, irreducible Riemannian
manifold of dimension $n$. Let $N$ denote the dimension of the space of
parallel spinors with respect to the Levi-Civita connection.
If $(M^n,g)$ is non-flat and $N>0$, then one of the following holds:
\begin{enumerate}
\item $n=2m$ ($m\geq 2$), the holonomy representation is
the vector representation of $\SU(m)$ on $\C^m$, and $N=2$ 
(``Calabi-Yau case'').
\item $n=4m$ ($m\geq 2$), the holonomy representation is
the vector representation of $\Sympl(m)$ on $\C^{2m}$, and $N=m+1$
(``hyper-K\"ahler case'').
\item $n=7$,  the holonomy representation is the unique $7$-dimensional
representation of $G_2$, and $N=1$ (``parallel $G_2$- or Joyce case'').
\item $n=8$,  the holonomy representation is the  spin
representation of $\Spin(7)$, and $N=1$ (``parallel $\Spin(7)$- or Joyce 
case'').
\end{enumerate}
\end{thm}
%
%-----------------------------------------------------------------------------
\section{A unified approach to non-integrable geometries}
%-----------------------------------------------------------------------------
%
\subsection{Motivation}
%-----------------------------------------------------------------------------
%
For $G$-structures defined by some tensor $\mathcal{T}$, it has been
for a long time customary to classify the possible types of structures
by the isotypic decomposition under $G$ of the covariant derivative 
$\nabla^g\mathcal{T}$. The integrable case is described by
$\nabla^g\mathcal{T}=0$, all other classes of non-integrable $G$-structures
correspond to  combinations of non-vanishing contributions
in the  isotypic decomposition and 
are described by some differential equation in $\mathcal{T}$. This was
carried out in detail for example for almost Hermitian manifolds 
(Gray/Hervella \cite{Gray&H80}), for $G_2$-structures in dimension $7$ 
(Fern\'andez/Gray \cite{FG}), for $\Spin(7)$-structures in dimension 
$8$ (Fern\'andez \cite{Fernandez86}) and for almost contact metric structures 
(Chinea/Gonzales \cite{Chinea&G90}). 

In this section, we shall present a simpler and unified approach
to non-integrable geometries. 
The theory of principal fibre bundles suggests  that the difference
$\Gamma$ between the Levi-Civita connection and the canonical $G$-connection
induced on the $G$-structure is a good measure for how much the given
$G$-structure fails to be integrable. By now, $\Gamma$ is widely known
as the \emph{intrinsic torsion} of the $G$-structure (see 
Section \ref{since-1980} for references). Although this is a 
``folklore'' approach, it is still not as popular as it could be.
Our presentation will follow the main lines of \cite{Fri2}. We will
see that it easily reproduces the classical results cited above with much less
computational work whilst having the advantage of being applicable 
to geometries not defined by tensors. Furthermore, it allows a uniform
and clean description of those classes of geometries admitting
$G$-connections with totally skew-symmetric torsion, and led to the discovery
of new interesting geometries.
%
%------------------------------------------------------------------------
\subsection{$G$-structures on Riemannian manifolds}\label{G-structures}
%------------------------------------------------------------------------
%
Let $G \subset \SO(n)$ be a closed subgroup of the orthogonal group
and decompose the Lie algebra $\so(n)$ into the Lie algebra $\g$ of $G$
and its orthogonal complement $\m$, i.\,e. $\so(n)=\g\oplus\m$.
Denote by $\mbox{pr}_{\g}$ and $\mbox{pr}_{\m}$ the 
projections onto $\g$ and $\m$, respectively. 
Consider an oriented Riemannian manifold $(M^n, g)$ and denote 
its frame bundle by $\mathcal{F}(M^n)$; it is a principal 
$\SO(n)$-bundle over $M^n$. By definition, a $G$-structure on $M^n$ 
is a reduction 
$\mathcal{R} \subset \mathcal{F}(M^{n})$ of the frame bundle to the subgroup 
$G$. The Levi-Civita connection is a $1$-form $Z$ on 
$\mathcal{F}(M^{n})$ with values in the Lie algebra $\so(n)$.
We restrict the Levi-Civita connection to $\mathcal{R}$ and decompose 
it with respect to the splitting $\g\oplus\m$:
\bdm
Z\big|_{T(\mathcal{R})} \ := \ Z^* \, \oplus \ \Gamma  .
\edm
Then, $Z^*$ is a connection in the principal $G$-bundle $\mathcal{R}$
and $\Gamma$ is a tensorial $1$-form of type Ad, i.\,e. a $1$-form on $M^{n}$ 
with values in the associated bundle $\mathcal{R} \times_G \m$.
By now, it has become standard to call $\Gamma$ the 
\emph{intrinsic torsion} of the $G$-structure (see Section \ref{since-1980}
for references). The 
$G$-structure $\mathcal{R}$ on $(M^n,g)$ is called
\emph{integrable} if $\Gamma$ vanishes, for this means that it is preserved by
the Levi-Civita connection and that $\Hol(\nabla^g)$ is a subgroup of $G$.
All $G$-structures with $\Gamma\neq 0$
are called \emph{non-integrable}; the basic classes of 
non-integrable $G$-structures are defined---via the decomposition of 
$\Gamma$---as the irreducible 
$G$-components of the representation $\R^{n} \otimes \m$. 
For an orthonormal frame $e_1, \ldots , e_n$ adapted to the reduction 
$\mathcal{R}$, the connection
forms $\omega_{ij} := g(\nabla^g e_i , e_j)$ of the Levi-Civita connection
define a $1$-form $\Omega : = (\omega_{ij})$ with values in the Lie algebra
$\so(n)$ of all skew-symmetric matrices. The form $\Gamma$ can then be
computed as the $\m$-projection of $\Omega$,
\bdm
\Gamma \ = \ \mbox{pr}_{\m}(\Omega) \ = \ 
\mbox{pr}_{\m}(\omega_{ij}).
\edm
Interesting is the case in which $G$ happens to be the
isotropy group of some tensor $\mathcal{T}$. Suppose that there is a 
faithful representation $\vrho : \SO(n) \rightarrow 
\SO(V)$ and a tensor $\mathcal{T} \in V$ such that
\bdm
G \ = \ \big\{ g \in \SO(n) : \vrho(g)\mathcal{T} = \mathcal{T}\big\}  .
\edm
The Riemannian covariant derivative of $\mathcal{T}$ is then given by the 
formula
\bdm
\nabla^g \mathcal{T} \ = \ \vrho_*(\Gamma)(\mathcal{T}) \, ,
\edm
where $\vrho_* : \so(n) \rightarrow \so(V)$ is the differential
of the representation. As a tensor, 
$\nabla^g\mathcal{T}$ is an element
of $\R^n \otimes V$. The algebraic $G$-types of $\nabla^g\mathcal{T}$
define the algebraic $G$-types of $\Gamma$ and vice versa. Indeed, 
we have
\begin{prop}[{\cite[Prop. 2.1.]{Fri2}}]
%--------------------------------------
The $G$-map
\bdm
\R^n \otimes \m \longrightarrow \R^n \otimes \mathrm{End}(V) \longrightarrow 
\R^n \otimes V
\edm given by $\Gamma \rightarrow \rho_*(\Gamma)(\mathcal{T})$ is injective.
\end{prop}
An easy argument in representation theory shows that for $G\neq \SO(n)$,
the $G$-representation $\R^n$ does always appear as  summand 
in the $G$-decomposition of $\R^n \otimes \m$.  
Geometrically, this module accounts precisely  for conformal transformations
of $G$-structures.  Let $(M^n,g,\mathcal{R})$ be a 
Riemannian manifold with a fixed geometric
structure and denote by $\hat{g} := e^{2f} \cdot g$ a conformal transformation
of the metric. There is a natural identification of the frame bundles
\bdm
\mathcal{F}(M^n, g) \ \cong \mathcal{F}(\hat{M}^n, \hat{g})
\edm
and a corresponding $G$-structure $\mathcal{\hat{R}}$. At the 
infinitesimal level, the conformal change is defined by the $1$-form 
$df$, corresponding to an $\R^n$-part in $\Gamma$.

We shall now answer the question under which conditions a given $G$-structure
admits a metric connection $\nabla$  with skew-symmetric torsion preserving 
the structure. For this, consider for any orthonormal basis $e_i$ of
$\m$ the $G$-equivariant map
\bdm
\Theta: \Lambda^3(\R^n) \longrightarrow \R^{n} \otimes \m, \quad
\Theta(T) \ := \ \sum_i (e_i \haken T) \otimes e_i.
\edm 
\begin{thm}[{\cite[Prop. 4.1]{Friedrich&I1}}]\label{thm2}
%-------------------------------------------------------
A $G$-structure $\mathcal{R} \subset \mathcal{F}(M^{n})$ of a 
Riemannian manifold admits a connection $\nabla$ with skew-symmetric 
torsion if and only if the $1$-form $\Gamma$ belongs to the image 
of $\Theta$, 
\bdm
2 \, \Gamma \ = \ -  \Theta(T)\ \text{ for some }T\in \Lambda^3(\R^n) . 
\edm
In this case the $3$-form $T$ is the torsion form of the connection.
\end{thm}
%
%almost contact metric structures, $\mbox{G}_2$-structures and $\Spin(7)$-structures (see \cite{Friedrich&I1}, \cite{Friedrich&I2}, \cite{Friedrich&I3}, 
%\cite{Iv} and \cite{Agricola1}, \cite{Agricola2}). 
%
\begin{dfn}
%----------
A metric $G$-connection $\nabla$ with torsion $T$ as in Theorem
\ref{thm2} will be called a \emph{characteristic connection}
and denoted by $\nabla^c$, $T=:T^c$ is called the 
\emph{characteristic torsion}. 
By construction, the holonomy $\Hol(\nabla^c)$ is a subgroup of $G$.
\end{dfn}
Thus, not every $G$-structure admits a characteristic connection.
If that is the case, $T^c$ is unique for all geometries we have investigated 
so far, and  it can easily be expressed in terms of the geometric data
(almost complex structure etc.). Henceforth, we shall just speak of 
\emph{the} characteristic connection. Due to its properties, it is an
excellent substitute for the Levi-Civita connection, which in
these situations is not adapted to the underlying geometric structure.
\begin{NB}
%---------
The canonical connection $\nabla^c$ of a naturally reductive homogeneous 
space is an example of a characteristic connection that satisfies 
in addition $\nabla^c T^c=\nabla^c\kr^c=0$; in this sense, 
geometric structures admitting a characteristic connection such that
$\nabla^c T^c=0$ constitute a natural generalization of naturally reductive
homogeneous spaces. As a consequence of the General Holonomy Principle 
(Corollary~\ref{cor-holonomy-principle}), $\nabla^c T^c=0$  implies that
the holonomy group $\Hol(\nabla^c)$ lies in the stabilizer $G_{T^c}$
of $T^c$.
\end{NB}
With this technique, we shall now describe special classes of non-integrable 
geometries, some new and others previously encountered. We order them by 
increasing dimension. 
%
%------------------------------------------------------
\subsection{Almost contact metric structures}
%------------------------------------------------------
%
At this stage, almost contact metric structures challenge
any expository paper because of the large number of classes.
Qualitatively, the situation is as follows.  The first classifications of
these structures proceeded in analogy to the Gray-Hervella set-up
for almost Hermitian manifolds (see Section~\ref{a-hermitian}) by examining 
the space of tensors with the same
symmetry properties as the covariant derivative of the fundamental
form $F$ (see Section \ref{exa-contact}) and decomposing it under the action 
of the structure group $G=\U(n)\x \{1\}$ using invariant theory.
Because the
$G$-action on $\R^{2n+1}$ is already not irreducible, this space decomposes 
into four $G$-irreducibles for $n=1$, into $10$ summands for $n=2$ and into 
$12$ for $n\geq 3$, leading eventually to $2^4$, $2^{10}$ and $2^{12}$
possible classes of almost contact metric structures (\cite{Aleksiev&G86},
\cite{Chinea&G90}, \cite{Chinea&M92}). Obviously, most of these
classes do not carry names and are not studied, and the result being what 
it is, the investigation of such structures is burdened by technical details
and assumptions. From the inner logic of non-integrable
geometries, it makes not so much sense to base their investigation
on the covariant derivative $\nabla^g F$ of $F$ or some other fundamental
tensor, as the Levi-Civita connection does not preserve the geometric
structure. This accounts for the technical complications that
one faces when following this approach.

For this section, we decided to restrict our attention to dimension five, 
this being the most relevant for the investigation of non-integrable geometries
(in dimension seven, it is reasonable to study contact structures 
simultaneously with $G_2$-structures). Besides, this case will illustrate the 
power of the intrinsic torsion concept outlined above. We shall
use throughout that $\R^5=\R^4\oplus\R$ with standard $\U(2)$-action on
the first term and trivial action on the second term. Let us look
at the decompositions of the orthogonal Lie algebras in dimension $4$ and
$5$. First, we have
\bdm
\so(4)\ =\ \Lambda^2(\R^4)\ =\ \un(2)\oplus \n^2.
\edm
Here, $\n^2$ is $\U(2)$-irreducible, while $\un(2)$ splits
further into $\su(2)$ and the span of $\Omega$. Combining this remark
with the characterization of these subspaces via the complex structure $J$ 
defining $\un(2)$, we obtain
\bdm
\un(2)\ =\ \{\omega\in \Lambda^2(\R^4)\,:\, J\omega=\omega\}
 \ = \ \su(2)\oplus \R\cdot\Omega ,\quad
\n^2\ =\ \{\omega\in \Lambda^2(\R^4)\,:\, J\omega=-\omega\}.
\edm
In particular, $\Lambda^2(\R^4)$ is the sum of three $\U(2)$-representations
of dimensions $1$, $2$ and $3$. For $\so(5)$, we deduce
immediatly
\bdm
\so(5)\ =\ \Lambda^2(\R^4\oplus\R)\ =\ \Lambda^2(\R^4)\oplus \R^4
\ =\ \un(2)\oplus (\R^4\oplus \n^2)\ =:\ \un(2)\oplus \m^6. 
\edm
 Thus,
the intrinsic torsion $\Gamma$ of a $5$-dimensional almost metric contact 
structure is an element of the representation space
\bdm
\R^5\ox\m^6\ =\  (\R^4\oplus\R)\ox (\R^4\oplus\n^2)\ =\
\n^2\oplus\R^4\oplus (\R^4\ox\n^2) \oplus(\R^4\ox\R^4). 
\edm
The last term splits further into trace-free symmetric, trace and
antisymmetric part, written for short as
\bdm
\R^4\ox\R^4\ =\ S^2_0(\R^4)\oplus \R\oplus \Lambda^2(\R^4).
\edm
The $9$-dimensional representation $S^2_0(\R^4)$ is again a sum
of two irreducible ones of dimensions $3$ and $6$, but we do not need 
this here. To decompose the representation $\R^4\ox\n^2$, we observe that the 
$\U(2)$-equivariant map $\Theta:\Lambda^3(\R^4)\ra\R^4\ox\n^2$ 
(see Section \ref{G-structures}) has $4$-dimensional irreducible image 
isomorphic to $\Lambda^3\R^4$ (which is again of dimension $4$); its 
complement is an inequivalent irreducible $\U(2)$-representation
of dimension $4$ which we call $V_4$. Consequently,
\bdm
\R^5\ox\m^6\ =\ \R\oplus \n^2\oplus\R^4\oplus S^2_0(\R^4)\oplus
\Lambda^2(\R^4)\oplus \Lambda^3(\R^4)\oplus V_4.
\edm
Taking into account the further splitting of 
$S^2_0(\R^4)\oplus\Lambda^2(\R^4)$, this is the sum of $10$ 
irreducible $\U(2)$-representations as claimed. On the other side,
\bdm
\Lambda^3(\R^5)\ =\ \Lambda^3(\R^4\oplus\R)\ =\
\Lambda^2(\R^4)\oplus \Lambda^3(\R^4).  
\edm
We found a unique copy of this $10$-dimensional space in the
$30$-dimensional space  $\R^5\ox\m^6$. Thus, we conclude from 
Theorem \ref{thm2}:
\begin{prop}
%-----------
A $5$-dimensional almost metric contact structure $(M^5,g\xi,\eta,\vphi)$
admits a unique characteristic connection if and only if its intrinsic torsion
is of class $\Lambda^2(\R^4)\oplus \Lambda^3(\R^4)$.
\end{prop}
In dimension $5$, skew-symmetry of the Nijenhuis tensor $N$ implies 
that it has to be zero, hence in the light of the more general Theorem
\ref{exa-contact}, the almost metric contact manifolds of class
$\Lambda^2(\R^4)\oplus \Lambda^3(\R^4)$ should coincide with the almost 
metric contact  structures with $N=0$ and $\xi$ a Killing vector field. That 
this is indeed
the case follows from the  classifications cited above. This class includes
for example all \emph{quasi-Sasakian manifolds} ($N=0$ and  $dF=0$),
see \cite{Kirichenko&R02}.
\begin{exa}\label{Heisenberg-Sasaki}
%-----------------------------------
Consider $\R^5$ with $1$-forms
\bdm
2 e_1\ =\ dx_1,\quad 2 e_2\ =\ dy_1,\quad 2 e_3\ =\ dx_2,\quad
2 e_4\ =\ dy_2,\quad 4 e_5\ =\ 4\eta \ =\ dz -y_1dx_1-y_2dx_2,
\edm
metric $g=\sum_i e_i\ox e_i$, and almost complex structure $\vphi$ defined 
in $\langle \xi\rangle^\perp$ by 
\bdm
\vphi(e_1)=e_2,\quad \vphi(e_2)=-e_1,\quad \vphi(e_3)=e_4,\quad
\vphi(e_4)=-e_3,\quad \vphi(e_5)=0.
\edm
Then $(\R^5,g,\eta,\vphi)$ is a Sasakian manifold, and the
torsion of its characteristic connection is of type $\Lambda^2(\R^4)$
(\cite[Example 3.D]{Fino94}) and explicitly given by
\bdm
T^c\ =\ \eta\wedge d\eta\ =\ 2(e_1\wedge e_2+ e_3\wedge e_4)\wedge e_5.
\edm
This example is in fact a left-invariant metric on a $5$-dimensional
Heisenberg group with $\Scal^g=-4$ and 
$\Scal^{\nabla^c}=\Scal^g-3 ||T||^2 /2=-16$ (see 
Theorem~\ref{curvature-identities}).
In a left-invariant frame, spinors are simply functions 
$\psi : \R^5 \rightarrow \Delta_5$
with values in the $5$-dimensional spin representation. 
In \cite{Friedrich&I2}, it is shown that there exist two
$\nabla^c$-parallel spinors $\psi_i$  with the additional property 
$F\cdot\psi_i=0\ (i=1,2)$. This implies $T\cdot \psi_i=0$, an equation of 
interest in superstring theory (see Section~\ref{super}). It turns out that
these spinors are constant, hence the same result holds for all
compact quotients $\R^5/\Gamma$ ($\Gamma$ a discrete subgroup).
\end{exa}
We recommend the articles \cite{Fino94} and \cite{Fino95} for a detailed 
investigation of the representation theory of almost metric contact 
structures (very much in the style of the book \cite{Salamon01})---in 
particular, the decomposition of the space of possible torsion tensors 
$\mathcal{T}$ of metric connections (see Proposition \ref{classessum}) under 
$\U(n)$ is being related to the possible classes for the intrinsic torsion.
%
%----------------------------------------------------------------
\subsection{$\SO(3)$-structures in dimension $5$}
%----------------------------------------------------------------
%
These structures were discovered by Th.~Friedrich in a 
systematic investigation of  possible $G$-structures for interesting 
non-integrable geometries (see \cite{Fri2}); until that moment,
it was generally believed that contact structures were the only remarkable
$G$-structures in dimension $5$.

The group $\SO(3)$ has a unique, real, irreducible
representation in dimension $5$. We consider the corresponding non-standard 
embedding
$\SO(3) \subset \SO(5)$ as well as the decomposition
\bdm
\so(5) \ = \ \so(3) \oplus \m^7 .
\edm
It is well known that the $\SO(3)$-representation $\m^7$ is the
unique, real, irreducible representation of $\SO(3)$ in dimension $7$. 
We decompose the tensor product into irreducible components
\bdm
\R^5 \otimes \m^7 \ = \ \R^3 \oplus \R^5 \oplus \m^7 \oplus
\mbox{E}^{9} \oplus \mbox{E}^{11} .
\edm
There are five basic types of $\SO(3)$-structures on $5$-dimensional 
Riemannian
manifolds. The symmetric spaces $\SU(3)/\SO(3)$ and $\SL(3,\R)/\SO(3)$
are  examples of $5$-dimensional Riemannian manifolds with an integrable 
$\SO(3)$-structure ($\Gamma = 0$). On the other hand, $3$-forms on $\R^5$
decompose into
\bdm
\Lambda^3(\R^5) \ = \ \R^3 \oplus \m^7  .
\edm
In particular, a conformal change of an $\SO(3)$-structure does \emph{not} 
preserve the property that the structure admits a connection with totally 
skew-symmetric torsion.

M.~Bobienski and P.~Nurowski investigated  $\SO(3)$-structures in
their articles \cite{Bobienski&N05} and \cite{Bobienski06}. In particular, 
they found a
ternary symmetric form describing the reduction to $\SO(3)$ and 
constructed many examples of non-integrable $\SO(3)$-structures
with non-vanishing intrinsic torsion. Recently, P.~Nurowski suggested a 
link to Cartan's work on isoparametric surfaces in spheres,
and predicted the existence of similar geometries in dimensions $8$, $14$ and 
$26$; we refer the reader to \cite{Nurowski06} for details. The case of
$\SO(3)$-structures  illustrates that new classes of non-integrable
geometries are still to be discovered beyond the well-established ones,
and that their study reveals deeper connections between areas which used to be
far from each other.
%
%--------------------------------------------------------------------------
\subsection{Almost Hermitian manifolds in dimension $6$}
\label{a-hermitian}
%--------------------------------------------------------------------------
%
We begin with the Gray-Hervella classification of almost
Hermitian manifolds and the consequences for the characteristic connection.
Although most of these results hold in all even dimensions, we shall henceforth
restrict our attention to the most interesting case, namely dimension $6$.

Let us consider a $6$-dimensional almost Hermitian manifold 
$(M^6, g,J)$, corresponding to a $\U(3)$-structure inside $\SO(6)$.
We decompose the Lie algebra into $\so(6) = \un(3) \oplus \m$
and remark that the $\U(3)$-representation in $\R^6$ is the real 
representation underlying $\Lambda^{1,0}$. Similarly, $\m$ is the real 
representation underlying $\Lambda^{2,0}$. We decompose the complexification
under the action of $\U(3)$:
\bdm
\Big(\R^6 \otimes \m\Big)^{\C} \ = \ \Big(\Lambda^{1,0} \otimes \Lambda^{2,0}
\oplus \Lambda^{1,0} \otimes \Lambda^{0,2}\Big)_{\R}^{\C}  .
\edm
The symbol $( \ldots )_{\R}^{\C}$ means that we understand the complex
representation as a real representation and complexify it. Next we split
the complex $\U(3)$-representations
\bdm
\Lambda^{1,0} \otimes \Lambda^{2,0} \ = \ \C^3 \otimes \Lambda^2(\C^3) \ 
= \ \Lambda^{3,0} \oplus \mbox{E}^8  ,
\edm
\bdm
\Lambda^{1,0} \otimes \Lambda^{0,2} \ = \ \C^3 \otimes 
\Lambda^2(\overline{\C}^3) \ = \ \C^3 \otimes \Lambda^2(\C^3)^* \ = \
(\C^3)^* \oplus \mbox{E}^6  .
\edm
$\mbox{E}^6$ and $\mbox{E}^8$ are irreducible $\U(3)$-representations
of complex dimensions $6$ and $8$, respectively. Finally we obtain
\bdm
\R^6 \otimes \m \ = \ \Lambda^{3,0} \oplus \mbox{E}^8 \oplus
\mbox{E}^6 \oplus (\C^3)^* \ =:\ \W^{(2)}_1 \oplus  \W^{(16)}_2 \oplus  
\W^{(12)}_3 \oplus \W^{(6)}_4.
\edm
Consequently, $\R^6 \otimes \m$ splits into four irreducible 
representations of real dimensions $2,  16, 12$ and $6$, that is, there
are four basic classes and a total of $16$ classes of $\U(3)$-structures on
$6$-dimensional Riemannian manifolds, a result known as 
Gray/Hervella-classification (\cite{Gray&H80}). Recently, F.~Mart\'{\i}n
Cabrera established the defining differential equations for these classes
solely in terms of the intrinsic torsion (see \cite{Martin-Cabrera05}), as we
shall state them for $G_2$-manifolds in the next section. 
In case we restrict the structure group to $\SU(3)$, 
the orthogonal complement $\su(3)^\perp$ is now $7$- instead of 
$6$-dimensional, and we obtain
\bdm
\R^6 \otimes \su(3)^\perp \ = \ \W_1 \oplus  \W_2 \oplus  \W_3 \oplus \W_4
\oplus \W_5,
\edm
where $\W_5$ is isomorphic to $\W_4\cong (\C^3)^*$. Furthermore, 
$\W_1$ and $\W_2$ are not irreducible anymore, but they split into
$\W_1=\W_1^+\oplus \W_1^-=\R\oplus\R$ and 
$\W_2=\W_2^+\oplus \W_2^-=\su(3)\oplus\su(3)$ (see \cite{Chiossi&S02},
\cite{Martin-Cabrera05}). Table \ref{U3-str} summarizes some remarkable 
classes of $\U(3)$-structures in dimension $6$. Most of these have by
now well-established names, while there is still some confusion for others;
these can be recognized by the parentheses indicating the different
names to be found in the literature. In the last column,
we collected characterizations of these classes (where several are listed,
these are to be understood as equivalent characterizations, not as 
simultaneous requirements).
Observe that we included in the last line a remarkable class of 
$\SU(3)$-structures, the
so called \emph{half-flat}  $\SU(3)$-structures
($\Psi^+$ is the real part of the $(3,0)$-form defined by $J$, see 
Section~\ref{stab-Un}). The name is chosen in order to suggest that half of
all $\W$-components vanish for these structures. Relying on results of
\cite{Hitchin01}, S.~Chiossi and S.~Salamon
described in \cite{Chiossi&S02} explicit metrics with Riemannian
holonomy $G_2$ on the product of any half-flat $\SU(3)$-manifold
with a suitable interval. A construction of half-flat 
$\SU(3)$-manifolds as $T^2$-principal fibre bundles over K\"ahlerian 
$4$-manifolds goes back to Goldstein 
and Prokushkin \cite{Goldstein&P02} and was generalized by Li, Fu and Yau
\cite{Li&Y05}, \cite{Fu&Y05}.
We refer to Section \ref{square-dirac} for 
examples of half-flat $\SU(3)$-structures on nilmanifolds. 
%
%--------------------------------
\begin{table}
\begin{center}
\setlength{\extrarowheight}{4pt}
\begin{tabular}{|c|c|c|}\hline 
name & class & characterization\\ \hline \hline
nearly K\"ahler manifold & $\W_1$ & \begin{tabular}{l} 
a) $(\nabla^g_X J)(X)=0$ \\ b) 
$N$ skew-sym.~and $\tau^2(d\Omega)=-9 d\Omega$ \\
c) $\exists$ real Killing spinor \end{tabular}\\[1mm] \hline 
almost K\"ahler manifold & $\W_2$ & $d\Omega = 0$ \\[1mm]  \hline
\begin{tabular}{c} balanced (almost Hermitian) or\\
(Hermitian) semi-K\"ahler m. \end{tabular}
& $\W_3$ & $N=0$ and $\delta\Omega=0$\\[1mm]\hline
locally conformally K\"ahler m. & $\W_4$ & 
$N=0$ and $d\Omega=\Omega\wedge\theta$ ($\theta$: Lee form)\\[1mm]  \hline
quasi-K\"ahler manifold & $\W_1 \oplus \W_2$ & $\nabla^g_X\Omega(Y,Z)+
\nabla^g_{JX}\Omega(JY,Z)=0$ \\[1mm]  \hline
Hermitian manifold & $\W_3 \oplus \W_4$ & a) $N=0$,\ \ b) 
$\tau^2(d\Omega)=- d\Omega$\\[1mm]  \hline
\begin{tabular}{c} (almost-)semi-K\"ahler or \\ (almost) cosymplectic m. 
\end{tabular} & $\W_1\oplus \W_2\oplus \W_3$ & a) $\delta\Omega=0$,\ \  b) 
$\Omega\wedge d\Omega=0$ \\ \hline
KT- or $\G_1$-manifold & $\W_1\oplus \W_3\oplus \W_4$ & \begin{tabular}{l} a)
$N$ is skew-symmetric\\ b) $\exists$ char. connection $\nabla^c$ 
 \end{tabular} \\[1mm]  \hline \hline
half-flat  $\SU(3)$-manifold & $\W_1^-\oplus \W_2^-\oplus \W_3$ &
$\Omega\wedge d\Omega=0$ and $d\Psi^+=0$
\\[1mm]  \hline
\end{tabular}
\bigskip
\end{center}
\caption{Some types of $\U(3)$- and $\SU(3)$-structures in dimension six.}
\label{U3-str}
\end{table}
%----------------------
%
\begin{thm}
%----------
An almost Hermitian $6$-manifold $(M^6,g,J)$ admits a characteristic
connection $\nabla^c$ if and only if it is of class 
$\W_1\oplus \W_3\oplus \W_4$, i.\,e.~if its Nijenhuis tensor $N$ is 
skew-symmetric. Furthermore, $\nabla^c$ is unique and given by the
expression
\bdm
g(\nabla^c_X Y, Z) \ :=\  
g(\nabla^g_X Y,Z) + \frac{1}{2}\left[N(X,Y,Z) +d\Omega(JX,JY,JZ)\right]
\edm
\end{thm}
\begin{proof}
%------------
In order to apply Theorem \ref{thm2}, we need the decomposition of $3$-forms
into isotypic $\U(3)$-representations,
\bdm
\Lambda^3(\R^6) \ = \ \Lambda^{3,0} \oplus \mbox{E}^6  \oplus
(\C^3)^* \ =\ \W_1\oplus \W_3\oplus \W_4.
\edm
Therefore, the image of  $\Theta$ consists of the sum 
$\W_1\oplus \W_3\oplus \W_4$ and $\Theta$ is injective, i.\,e.,
there exists at most one characteristic torsion form. 
From the Gray-Hervella classification,
we know that almost Hermitian manifolds without  
$\W_2$-part (so-called \emph{$\G_1$-manifolds} or \emph{KT-manifolds},
standing for `K\"ahler with torsion', although not K\"ahler) are 
characterized by the property that their Nijenhuis tensor $N$ is 
skew-symmetric. In Lemma \ref{lem-almost-Herm-conn}, it was shown that
the stated connection fulfills all requirements, hence  by 
uniqueness it coincides with the characteristic connection.
\end{proof}
$\Lambda^3(\R^6)$ admits another decomposition. The map
\bdm
\tau:\ \Lambda^3(\R^6)\lra\Lambda^3(\R^6),\quad
\tau(T)\ :=\ \sum_{i=1}^6(e_i\haken\Omega)\wedge (e_i\haken T)
\edm
is $\U(3)$-equivariant. Its square $\tau^2$ is diagonalizable with
eigenspaces $\W_1$ (eigenvalue $-9$) and $\W_3\oplus\W_4$ (eigenvalue $-1$).
This explains the second characterization of these two classes in Table \ref{U3-str}.
\begin{NB}[Parallel torsion]
%---------------------------
In Example~\ref{exa-almost-hermitian}, it had been observed that
the characteristic torsion of nearly K\"ahler manifolds is always
parallel (Kirichenko's Theorem). Another interesting class of
almost Hermitian $\G_1$-manifolds with this property are the
so-called \emph{generalized Hopf structures}, that is, 
locally conformally K\"ahler manifolds (class $\W_4$, sometimes abbreviated
lcK-manifolds) with parallel
\emph{Lee form} $\theta:=\delta\Omega\circ J\neq 0$ (in fact, 
$\nabla^g\theta=0$ and $\nabla^c T^c=0$ are equivalent conditions 
for $\W_4$-manifolds). Besides the classical Hopf manifolds, they include
for example total spaces of flat principal $S^1$-bundles over compact 
$5$-dimensional Sasaki manifolds (see \cite{Vaisman76}, \cite{Vaisman79}
for details); generalized Hopf structures are never Einstein.
We recommend the book by S.~Dragomir and L.~Ornea \cite{Dragomir&O98}
as a general reference and the articles \cite{Belgun00}, \cite{Fujiki&P05}
for  complex lcK-surfaces.

In his thesis, N. Schoemann investigates almost hermitian structures with
parallel  skew-symmetric torsion in dimension $6$. A full classification 
of the possible 
algebraic types of the torsion form is worked out, and based on this a 
systematic description
of the possible geometries is given. In addition numerous new examples
are constructed (and, partially, classified) on naturally
reductive spaces (including compact spaces
with closed torsion form) and on nilmanifolds (see \cite{Alexandrov&F&S04} 
and \cite{Schoemann06}).
\end{NB}
\begin{NB}[Almost K\"ahler manifolds]
%------------------------------------
The geometry of almost K\"ahler manifolds is strongly related to famous
problems in differential geometry.  W.~Thurston was the first  to
construct an explicit compact symplectic manifold with $b_1=3$, hence
that does not admit a K\"ahler structure \cite{Thurston76}. E.~Abbena
generalized this example and gave a natural associated metric which
makes it into an almost  K\"ahler non-K\"ahler manifold \cite{Abbena84};
many more examples of this type have been constructed since then.

In 1969, S.\,I.~Goldberg conjectured that a compact almost K\"ahler-Einstein
manifold is K\"ahler \cite{Goldberg69}. In this generality, the 
conjecture is still open. 
K.~Sekigawa proved it under the assumption of
non-negative scalar curvature  \cite{Sekigawa87}, and it is known that
the conjecture is false for non-compact manifolds: P.~Nurowski and  
M.~Przanowski gave the first example of a $4$-dimensional  Ricci-flat 
almost-K\"ahler non-K\"ahler manifold \cite{Nurowski&P97}, J.~Armstrong
showed some non-existence results \cite{Armstrong98}, while V.~Apostolov, 
T.~Dr\u{a}ghici 
and A.~Moroianu constructed non-compact counterexamples to the conjecture in 
dimensions $\geq 6$ \cite{Apostolov&D&M01}. Different partial results with
various additional curvature assumptions are now available. The 
integrability conditions for almost K\"ahler manifolds were studied in full 
generality in \cite{Apostolov&A&D02} and \cite{Kirchberg05}.
\end{NB}
\begin{NB}[Nearly K\"ahler manifolds]
%------------------------------------
We close this section with some additional remarks on nearly
K\"ahler manifolds. Kirichenko's Theorem ($\nabla^c T^c=0$) implies that 
$\Hol(\nabla^c)\subset \SU(3)$, that the first Chern class of the tangent 
bundle $c_1(TM^6,J)$ vanishes,  $M^6$ is spin and that the metric 
is Einstein.
The only known examples are homogeneous metrics on $S^6, \C\P^3, S^3\x S^3$
and on the flag manifold $F(1,2)=\U(3)/\U(1)\x \U(1)\x \U(1)$, 
although (many?) more are expected to exist. It was shown 
that these exhaust all nearly K\"ahler manifolds
that are locally homogeneous  (see \cite{Butruille05}) or
satisfying $\Hol(\nabla^c)\neq \SU(3)$ (see \cite{BelgunMoroianu}).  
A by now classical result asserts that a $6$-dimensional spin manifold
admits a real Killing spinor if and only if it is nearly
K\"ahler (see \cite{Friedrich&G85}, \cite{Grunewald90} and \cite{BFGK}). 
Finally, more recent structure theorems justify why nearly K\"ahler
manifolds are only interesting in dimension $6$: any complete
simply connected nearly  K\"ahler manifold is locally a Riemannian product
of K\"ahler manifolds, twistor spaces over  K\"ahler manifolds and 
$6$-dimensional nearly  K\"ahler manifolds 
(see \cite{Nagy02a}, \cite{Nagy02b}).
\end{NB}
%
%-------------------------------------------------------
\subsection{$G_2$-structures in dimension $7$}
%------------------------------------------------------
%
We consider $7$-dimensional Riemannian manifolds equip\-ped with a 
$G_2$-structure. Since  $G_2$ is the isotropy
group of a $3$-form $\omega$ of general type, a $G_2$-structure
is a triple $(M^7, g, \omega)$ consisting of a $7$-dimensional Riemannian
manifold and a $3$-form $\omega$ of general type at any point.  
We decompose the $G_2$-representation (see \cite{Friedrich&I1})
\bdm
\R^7 \otimes \m \ = \ \R \oplus \Lambda^2_{14} \oplus
\Lambda^3_{27}\oplus  \R^7 \ =:\ \X^{(1)}_1\oplus\X^{(14)}_2 \oplus  
\X^{(27)}_3 \oplus \W^{(7)}_4,
\edm
and, consequently, there are again four basic classes and a total of $16$
classes $G_2$-structures (namely, parallel $G_2$-manifolds and $15$
non-integrable $G_2$-structures).  This result is known as
the Fern\'andez/Gray-classification of $G_2$-structures (see \cite{FG}); some 
important classes are again summarized in tabular form, see 
Table \ref{G2-str}. The different
classes of $G_2$-structures can be characterized by differential
equations. They can be written in a unified way as
\bdm
d\omega\ =\ \lambda\cdot *\omega+\frac{3}{4}\theta\wedge\omega+*\tau_3,\quad
\delta\omega\ =\ - * d * \omega \ =\ -*(\theta\wedge *\omega)
+*(\tau_2\wedge\omega),
\edm
where $\lambda$ is a scalar function corresponding to the $\X_1$-part
of the intrinsic torsion $\Gamma$, $\tau_2$, $\tau_3$ are $2$- resp.~$3$-forms
corresponding to its $\X_2$ resp.~$\X_3$-part and $\theta$ is a $1$-form
describing its ~$\X_4$-part, which one sometimes calls the
\emph{Lee form} of the $G_2$-structure. This accounts for some of
the characterizations listed in Table \ref{G2-str}.
For example, a $G_2$-structure is of type $\X_1$ 
(\emph{nearly parallel} $G_2$-structure) if and only if there exists a number 
$\lambda$ (it has to be constant in this case) such that
$ d \omega =  \lambda * \omega$
holds. Again, this condition is equivalent to the existence of a real
Killing spinor  and the metric has to be  Einstein \cite{FKMS};
more recently, the Riemannian curvature properties of arbitrary 
$G_2$-manifolds have been discussed in detail by R.~Cleyton and S.~Ivanov
\cite{Cleyton&I06a}. $G_2$-structures of type $\X_1\oplus\X_3$ 
(\emph{cocalibrated} $G_2$-structures) are characterized by the
condition that the $3$-form is coclosed, $\delta \omega^3 = 0$. 
%
%------------------------------
\begin{table}
\begin{center}
\setlength{\extrarowheight}{4pt}
\begin{tabular}{|c|c|c|}\hline 
name & class & characterization\\ \hline \hline
nearly parallel $G_2$-manifold & $\X_1$ & \begin{tabular}{l} a) 
$d\omega=\lambda\, *\omega$ for some $\lambda\in\R$ \\
b) $\exists$ real Killing spinor \end{tabular}\\[1mm] \hline 
\begin{tabular}{c} almost parallel or closed (or\\
calibrated symplectic) $G_2$-m.\end{tabular} 
& $\X_2$ & $d\omega = 0$ \\[1mm]  \hline
 balanced  $G_2$-manifold & $\X_3$ & $\delta\omega=0$ and $d\omega\wedge\omega=0$ 
 \\[1mm]  \hline
locally conformally parallel $G_2$-m. & $\X_4$ & 
\begin{tabular}{c} $d\omega=\frac{3}{4}\theta\wedge\omega $ and\\
 $d*\omega = \theta\wedge *\omega$  for some $1$-form $\theta$\end{tabular}
\\[1mm]  \hline
\begin{tabular}{c}
cocalibrated (or semi-parallel\\ or cosymplectic ) $G_2$-manifold\end{tabular}
& $\X_1 \oplus \X_3$ & $\delta\omega=0$ \\[1mm]  \hline
\begin{tabular}{c} locally conformally (almost)\\
  parallel $G_2$-manifold \end{tabular} & $\X_2 \oplus \X_4$ &
$d\omega=\frac{3}{4}\theta\wedge\omega $ \\[1mm]  \hline
$G_2T$-manifold & $\X_1\oplus \X_3\oplus \X_4$ & \begin{tabular}{l} a)
$d*\omega = \theta\wedge *\omega$  for some $1$-form $\theta$ \\
b) $\exists$ char. connection $\nabla^c$  \end{tabular} \\[1mm]  \hline
\end{tabular}
\bigskip
\end{center}
\caption{Some types of $G_2$-structures in dimension seven.}
\label{G2-str}
\end{table}
%-------------------------------- 
%
Under the restricted action of $G_2$, one obtains the following
isotypic decomposition of $3$-forms on $\R^7$:
\bdm
\Lambda^3(\R^7) \ = \ \R  \oplus \Lambda^3_{27}\oplus \R^7\ =\
\X_1\oplus\X_3\oplus\X_4.
\edm
This explains the first part of the following theorem
and the acronym `$G_2T$-manifolds' for this class: it stands for `$G_2$ with 
(skew) torsion'. The explicit
formula for the characteristic torsion may be derived directly from
the properties of $\nabla^c$.
\begin{thm}[{\cite[Thm. 4.8]{Friedrich&I1}}]\label{thm-G2-char-torsion}
%----------------------------------------------------------------------
A $7$-dimensional manifold $(M^7,g,\omega)$ with a fixed
$G_2$-structure $\omega\in\Lambda^3(M^7)$  admits a characteristic
connection $\nabla^c$ if and only if it is of class 
$\X_1\oplus \X_3\oplus \X_4$, i.\,e.~if there exists a $1$-form
 $\theta$ such that $d*\omega = \theta\wedge *\omega$.
 Furthermore, $\nabla^c$ is unique and given by the
expression
\bdm
\nabla^c_X Y \ :=\ \nabla^g_X Y + \frac{1}{2}\left[-*d\omega - \frac{1}{6}
g(d\omega,*\omega)\omega+*(\theta\wedge\omega)
\right].
\edm
$\nabla^c$ admits (at least) one parallel spinor.
\end{thm}
This last remarkable property is a direct consequence of our investigation of
geometric stabilizers, as explained in Corollary \ref{G2-equiv-par-spin}.
For a nearly parallel $G_2$-manifold, the $\nabla^c$-parallel spinor
coincides with the Riemannian Killing spinor and the manifold turns out to be
$\nabla^c$-Einstein \cite{Friedrich&I1}.
Some subtle effects occur when more spinors enter the play, as we shall now
explain. First, we recall the fundamental theorem on Killing spinors in
dimension $7$:
\begin{thm}[{\cite{FKMS}}]
%-------------------------
A $7$-dimensional simply connected compact Riemannian spin manifold $(M^7,g)$ 
admits
\begin{enumerate}
\item one real Killing spinor if and only if it is a nearly parallel $G_2$-manifold;
\item two real Killing spinors if and only if it is a Sasaki-Einstein manifold;
\item three real Killing spinors if and only if it is a $3$-Sasaki manifold.
\end{enumerate}
\end{thm}
Furthermore, $3$ is also the maximal possible number of Killing spinors
for $M^7\neq S^7$. On the other side, the characteristic connection
$\nabla^c$ of a $G_2T$-manifold has $2$ resp.~$4$ parallel spinors if its 
holonomy reduces further to $\SU(3)$ resp.~$\SU(2)$. But there is no general
argument identifying Killing spinors with parallel spinors:
the characteristic connection of a Sasaki-Einstein manifold does not
necessarily admit parallel spinors (see \cite{Friedrich&I1}, 
\cite{Friedrich&I2}), a  $3$-Sasaki manifold does not even admit a
characteristic connection in any reasonable sense (each Sasaki structure
has a characteristic connection, but it does not preserve the other two
Sasaki structures), see Section \ref{exa-3-Sasaki} and \cite{Agricola&F03a}. 
This reflects the fact that  Sasaki-Einstein and $3$-Sasaki manifolds
do not fit too well into the general framework of $G$-structures.
\begin{NB}[Parallel torsion]
%---------------------------
For a nearly parallel $G_2$-manifold, the explicit formula from
Theorem \ref{thm-G2-char-torsion} implies 
that the characteristic torsion $T^c$ is proportional to $\omega$, hence it
is trivially $\nabla^c$-parallel. For the larger class of
cocalibrated $G_2$-manifolds (class $\X_1\oplus\X_3$), the case of parallel 
characteristic torsion has  been investigated systematically  by Th.~Friedrich
(see \cite{Friedrich06}). Again, many formerly unknown examples have been
constructed, for example, from deformations of $\eta$-Einstein Sasaki 
manifolds, from $S^1$-principal fibre bundles over $6$-dimensional K\"ahler 
manifolds or from naturally reductive spaces.
\end{NB}
%
%
%----------------------------------------------------
\subsection{$\Spin(7)$-structures in dimension $8$}
%----------------------------------------------------
%
Let us consider $\Spin(7)$-structures on $8$-dimensional Riemannian
manifolds. The subgroup $\Spin(7) \subset \SO(8)$ is the real 
$\Spin(7)$-representation $\Delta_7 = \R^8$, the complement 
$\m = \R^7$ is the standard $7$-dimensional representation and the 
$\Spin(7)$-structures on an $8$-dimensional Riemannian manifold $M^8$ 
correspond to the irreducible components of the tensor product
\bdm
\R^8 \otimes \m \ = \ \R^8 \otimes \R^7 \ = \ \Delta_7 \otimes \R^7 \ = \ 
\Delta_7 \oplus \mbox{K}  \ = \ \R^8 \oplus \mbox{K} ,
\edm
where $\mbox{K}$ denotes the kernel of the Clifford multiplication 
$\Delta_7 \otimes \R^7 \to \Delta_7$. It is well known that 
$\mbox{K}$ is an irreducible $\Spin(7)$-representation, i.e.~there are two 
basic classes of $\Spin(7)$-structures (the Fern\'andez classification of
$\Spin(7)$-structures, see \cite{Fernandez86}).
For $3$-forms, we find the isotypic decomposition
\bdm
\Lambda^3(\R^8) \ = \ \Delta_7 \oplus \mbox{K},
\edm
showing that $\Lambda^3(\R^8)$ and  $\R^8 \otimes \m$ are isomorphic.
Theorem \ref{thm2} yields immediately that \emph{any $\Spin(7)$-structure on an
$8$-dimensional Riemannian manifold admits a unique connection with 
totally skew-symmetric torsion}. The explicit formula for its characteristic
torsion may be found in \cite{Ivanov04}.
%
%-----------------------------------------------------------------------------
\section{Weitzenb\"ock formulas for Dirac operators with torsion}
%-----------------------------------------------------------------------------
%
\subsection{Motivation}
%-----------------------------------------------------------------------------
 The question whether or
not the characteristic connection of a $G$-structure admits parallel  
tensor fields  differs radically from the corresponding problem  for 
the Levi-Civita connection. In particular, one is interested in  the 
existence of parallel spinor fields,
interpreted in superstring theory as  supersymmetries of the
model. The main analytical tool for the investigation of parallel spinors
is the Dirac operator and several remarkable identities for it.
We discuss two identities for the square of the Dirac operator. While
the first one is straightforward and merely of  
computational difficulty, the second relies on 
comparing the Dirac operator corresponding to the connection
with torsion  $T$ with the spinorial Laplace operator corresponding
to the connection with torsion  $3\,T$. Such an argument has been used in
the literature at several places. The first was probably
S.~Slebarski (\cite{SlebarskiI87}, \cite{SlebarskiII87}) who noticed
that on a naturally reductive space, the connection with torsion one-third
that of the canonical connection behaves well under fibrations; 
S.~Goette applied this property to the computation of the $\eta$-invariant
on homogeneous spaces \cite{Goette99}.
J.-P.~Bismut used such a rescaling  for proving an index theorem for Hermitian 
manifolds \cite{Bis}. It is implicit in Kostant's work on a `cubic Dirac 
operator', which can be understood as an identity in the Clifford algebra for 
the symbol of the Dirac operator of the rescaled canonical
connection on a naturally reductive space (\cite{Kostant99}, \cite{Agri}). 
\subsection{The square of the Dirac operator and parallel spinors}
\label{square-dirac}
%---------------------------------------------------------------------------
%
Consider a Riemannian spin manifold $(M^n,g, T)$ with a $3$-form 
$T\in\Lambda^3(M^n)$ as well as the one-parameter family of linear metric 
connections with skew-symmetric torsion ($s\in\R$),
\bdm
\nabla^s_X Y\ :=\ \nabla^{g}_X Y + 2  s \, T(X,Y,-)\,.
\edm
In particular, the superscript $s=0$ corresponds to the Levi-Civita 
connection and $s=1/4$ to the connection with torsion $T$ considered before.
As before, we shall also sometimes use the superscript $g$ to denote
the Riemannian quantities corresponding to $s=0$.
These connections can all be lifted to connections on the spinor
bundle $\Sigma M^n$, where they take the expression
\bdm
\nabla^s_X \psi\ :=\ \nabla^{g}_X \psi + s ( X\haken T) \cdot 
\psi\,.
\edm
Two important elliptic operators may be defined on  $\Sigma M^n$,
namely, the Dirac operator and the spinor Laplacian associated with
the connection $\nabla^s$:
\bdm
D^s\ :=\ \sum_{k=1}^n e_k\cdot\nabla^s_{e_k}\ =\ D^0+3s T,\quad
\Delta^s(\psi)\ =\ (\nabla^s)^*\nabla^s\psi \ =\ -\sum_{k=1}^n \nabla^s_{e_k}
\nabla^s_{e_k} \psi +\nabla^s_{\nabla^{g}_{e_i} e_i}\psi\,.
\edm
By a result of Th.~Friedrich and S.~Sulanke \cite{Friedrich&S79}, the Dirac
operator $D^\nabla$ associated with any metric connection 
$\nabla$ is formally self-adjoint if and only if the
$\nabla$-divergence 
$\mathrm{div}^\nabla(X) := \sum_i g(\nabla_{e_i}X,e_i)$ of any 
vector field $X$ coincides with its Riemannian $\nabla^g$-divergence.
Writing $\nabla=\nabla^g+ A$, this is manifestly equivalent to
$\sum_i g(A(e_i,X),e_i)=0$ and trivially satisfied for metric connections
with totally skew-symmetric torsion\footnote{One checks that it also holds for
metric connections with vectorial torsion, but not for connections of
Cartan type $\mathcal{A}'$.}.

Shortly after P.~Dirac introduced the Dirac operator, E.~Schr\"odinger noticed
the existence of a remarkable formula for its square
\cite{Schroedinger32}. Of course, since the concept of spin manifold had not 
yet been established, all arguments of that time were of local nature, but
contained already all important ingredients that would be established in a more
mathematical way later. By the sixties and the seminal work of
Atiyah and Singer on index theory for elliptic differential operators,
Schr\"odinger's article was almost forgotten and the formula rediscovered
by A.~Lichnerowicz \cite{Lichnerowicz63}. In our notation, the 
\emph{Schr\"odinger-Lichnerowicz  formula} states that  
\bdm
(D^0)^2 \ =\ \Delta^0 +\frac{1}{4}\Scal^0.
\edm
Our goal is to derive useful relations for the square of $D^s$.
In order to state the first formula, let us introduce the
first order differential operator
\be\label{operator-D}
\D^s\psi \ :=\ \sum_{k=1}^n (e_k\haken\mathrm{T})\cdot \nabla^s_{e_k}\psi \ =\
\D^0\psi + s  \sum_{k=1}^n (e_k\haken \T)\cdot (e_k\haken \T)
\cdot\psi.
\ee
\begin{thm}[{\cite[Thm 3.1, 3.3]{Friedrich&I1}}]\label{FI-SL-AC}
%---------------------------------------------------------------
Let $(M^n,g,\nabla^s)$ be an $n$-dimensional Riemannian spin
manifold with a metric connection $\nabla^s$ of skew-symmetric
torsion $4 s \cdot T$. Then, the square of the Dirac operator
$D^s$ associated with  $\nabla^s$ acts on an arbitrary spinor field 
$\psi$ as
\be\label{FI-SL}
(D^s)^2\psi \ =\ \Delta^s(\psi) + 3  s \, d T \cdot\psi
- 8  s^2\, \sigma_{T}\cdot\psi+
2 s \, \delta T \cdot\psi - 4  s \, \D^s\psi + \frac{1}{4}\,\Scal^s 
\cdot \psi.
\ee
Furthermore, the anticommutator of $D^s$ and $T$ is
\be\label{D-omega-anticomm}
D^s \circ T + T \circ D^s\ =\ dT+\delta T- 8 s \, \sigma_{T}
-2 \, \D^s.
\ee
\end{thm}
$\Scal^s$ denotes the scalar curvature of the connection $\nabla^s$. Remark
that $\Scal^0 = \Scal^g$ is the usual scalar curvature of the underlying
Riemannian manifold $(M^n,g)$ and that the relation
$\Scal^s=\Scal^0-24 s^2 ||T||^2$ holds. Moreover, the divergence
$\delta T$ can be taken with respect to any connection $\nabla^s$
from the family, hence we do not make a notational difference between them
(see Proposition~\ref{delta-form}).

This formula for $(D^s)^2$ has the disadvantage of still
containing a first order differential operator with uncontrollable spectrum
as well as several $4$-forms that are difficult to treat algebraically,
hence it is not suitable for deriving vanishing theorems.
It has however a nice application in the study of $\nabla^s$-parallel
spinors for \emph{different} values of $s$. As motivation, let's 
consider the following example:
\begin{exa}\label{Lie-grp} 
%-------------------------
Let $G$ be a simply connected Lie group, $g$ a biinvariant metric and consider
the torsion form $T(X,Y,Z) := g([X,Y],Z)$. The connections $\nabla^{\pm 1/4}$
are flat \cite{Kobayashi&N2}, hence they both admit  non-trivial
parallel spinor fields.
\end{exa}
Such a property for the connections with torsions $\pm T$ is required
in some superstring models. 
Theorem~\ref{FI-SL-AC} now implies that there cannot be many values $s$
admitting $\nabla^s$-parallel spinors.
\begin{thm}[{\cite[Thm. 7.1.]{Agricola&F03a}}]\label{parameter}  
%-------------------------------------------------------------
Let $(M^n,g,T)$ be a compact spin manifold, $\nabla^s$ the family of
metric connections defined by $T$ as above. For any $\nabla^s$-parallel spinor 
$\psi$, the following formula holds: 
\bdm
64 \, s^2 \int_{M^n} \langle \sigma_{T} \cdot \psi \, , \, \psi 
\rangle \, + \, \int_{M^n} \Scal^s \cdot ||\psi||^2 \ = \ 0  .
\edm
If the mean value of $\langle \sigma_{T}  \cdot \psi \, , \, \psi \rangle$ 
does not vanish,  the parameter $s$ is given by 
\bdm
s \ = \ \frac{1}{8} \int_{M^n} \langle dT \cdot \psi \, , \, \psi 
\rangle \Big/ 
\int_{M^n} \langle \sigma_{T} \cdot \psi \, , \, \psi \rangle . 
\edm
If the mean value of $\langle \sigma_{T}  \cdot \psi \, , \, \psi \rangle$ 
vanishes,  the parameter $s$ depends only on the Riemannian scalar
curvature and on the length of the torsion form, 
\bdm
0 \ = \ \int_{M^n} \mathrm{Scal}^s \ = \ \int_{M^n} \Scal^g \, - \,
24s^2 \int_{M^n} ||T||^2  .
\edm
Finally, if the $4$-forms $d T$ and $\sigma_{T}$ are proportional
(for example, if $\nabla^{1/4}T=0$),
there are at most three parameters with $\nabla^s$-parallel spinors.   
\end{thm}
\begin{NB}\label{NB-dT-prop-sigma}
%---------------------------------
The property that $d T$ and $\sigma_{T}$ are proportional is more
general than requiring parallel torsion. For example, it holds for the whole 
family of connections $\nabla^t$ on naturally reductive spaces discussed in 
Sections~\ref{exa-nat-red} and \ref{kostant}, but its torsion is 
$\nabla^t$-parallel only for $t=1$.
\end{NB}
\begin{exa} 
%----------
On the $7$-dimensional Aloff-Wallach space $N(1,1)=\SU(3)/S^1$,
one can construct a non-flat connection 
such that $\nabla^{s_0}$ \emph{and} $\nabla^{-s_0}$ admit parallel spinors
for suitable $s_0$, hence showing that both cases from 
Theorem~\ref{parameter} can actually occur in non-trivial situations.
On the other hand, it can be shown that on a $5$-dimensional
Sasaki manifold,  only the characteristic connection $\nabla^c$ can have
parallel spinors \cite{Agricola&F03a}.
\end{exa}
Inspired by the homogeneous case (see Section \ref{kostant}),
we were looking for an alternative comparison of $(D^s)^2$ with the
Laplace operator of some \emph{other} connection $\nabla^{s'}$ from the same
family. Since $(D^s)^2$ is a symmetric second order differential
operator with metric principal symbol, a very general result by
P.\,B.~Gilkey claims that there exists a connection $\nabla$
and an endomorphism $E$ such that $(D^s)^2=\nabla^*\nabla+E$ \cite{Gilkey75}.
Based on the results of Theorem~\ref{FI-SL-AC}, one shows:
\begin{thm}[Generalized Schr\"odinger-Lichnerowicz formula, 
{\cite[Thm. 6.2]{Agricola&F03a}}]\label{new-weitzenboeck}
%--------------------------------------------------------------------
The spinor Laplacian $\Delta^s$ and the square of the Dirac operator
$D^{s/3}$ are related by
\bdm
(D^{s/3})^2\ =\ \Delta^s + s \, d T 
+\frac{1}{4}\,\Scal^g - 2s^2 \, ||T||^2.
\edm
\end{thm}
We observe that $D^{s/3}$ appears basically by quadratic completion.
A first consequence is a non linear version of Theorem \ref{pareucl}. 
\begin{thm}[{\cite[Thm. 6.3]{Agricola&F03a}}]\label{parallelflach2}  
%------------------------------------------------------------------
Let $(M^n, g,  T)$ be a compact, Riemannian spin
manifold of non positive scalar curvature, $\Scal^g \leq 0$. If there
exists a solution $\psi \neq 0$  of the equations
\bdm
\nabla_X\psi \ = \ \nabla^g_X \psi +\frac{1}{2}(X\haken T)\cdot\psi\ =\ 0,
\quad \langle dT\cdot\psi,\psi\rangle\ \leq\ 0,
\edm
the $3$-form and the scalar curvature vanish, $T = 0 = \Scal^g$, 
and $\psi$ is parallel with respect to the Levi-Civita connection.
\end{thm}
Theorem \ref{parallelflach2} applies, in particular, to Calabi-Yau or 
Joyce manifolds, where we know that $\nabla^g$-parallel spinors exist 
by Wang's Theorem (Theorem \ref{thm-wang}). Let us perturb the connection 
$\nabla^g$ by a suitable $3$-form (for example, a closed one). Then the
new connection $\nabla$ does not admit $\nabla$-parallel spinor 
fields: the Levi-Civita connection and its parallel spinors are thus, 
in some sense, rigid. Nilmanifolds are a
second family of examples where the theorem applies. A further family of
examples arises from certain naturally reductive spaces with torsion form
$T$ proportional to the torsion form of the canonical connection, see
\cite{Agri}.
From the high energy physics point of view, a parallel spinor is interpreted 
as a supersymmetry transformation. Hence the physical problem behind the
above question (which in fact motivated our investigations) is really 
whether a free ``vacuum solution'' can also carry a non-vacuum 
supersymmetry, and how the two are related. 

Naturally, Theorem \ref{parallelflach2} raises the question to which extent
compactness is really necessary. We shall now show that it
is by using the equivalence between the inclusion $\Hol(\nabla)\subset G_2$ and
the existence of a $\nabla$-parallel spinor for a metric connection
with skew-symmetric torsion known from Corollary \ref{G2-equiv-par-spin}.
For this, it is sufficient to find a $7$-dimensional Riemannian 
manifold $(M^7,g)$ whose Levi-Civita  connection has a parallel spinor 
(hence is  Ricci-flat, in particular), but also admits a $\nabla$-parallel 
spinor for some other metric connection with skew-symmetric torsion. 

Gibbons  {\sl et al.}~produced non-complete metrics with Riemannian 
holonomy $G_2$ in \cite{GLPS}.
Those metrics have among others the interesting feature of admitting 
a $2$-step nilpotent isometry group $N$ acting on orbits of codimension one.
By \cite{Ch-F} such metrics are locally conformal to homogeneous metrics on
rank-one solvable extension of $N$, and the induced $\SU(3)$-structure
on $N$ is half-flat. In the same paper all half-flat $\SU(3)$
structures on $6$-dimensional nilpotent Lie groups whose
rank-one solvable extension is endowed with a conformally parallel 
$G_2$ structure were classified. Besides the torus,
there are exactly six instances, which we considered in relation to
the problem posed. It turns out that four metrics of the six only carry 
integrable $G_2$ structures, thus reproducing the
pattern of the compact situation, whilst one admits complex 
solutions, a physical interpretation for which is still lacking.
The remaining solvmanifold $(\mathrm{Sol},g)$---which has \emph{exact} 
Riemannian holonomy
$G_2$---provides a positive answer to both questions posed above,
hence becoming the most interesting.  The Lie algebra associated to this
solvmanifold  has structure equations
\bdm
\ba{l}
 [e_i, e_7] = \tfrac 35 m e_i,\ i = 1, 2, 5,\quad
 [e_j, e_7] = \tfrac 65 m e_j,\ j = 3, 4, 6,\\[2pt]
 [e_1, e_5] = - \tfrac 25 m e_3,\ [e_2, e_5] = - \tfrac 25 m e_4,\ 
 [e_1, e_2] = - \tfrac 25 m e_6.
\ea
\edm
The homogeneous metric it bears can be also seen as a $G_2$ metric on
the product $\R\times \mathbb{T}$, where $ \mathbb{T}$ is the total
space of a $T^3$-bundle over another $3$-torus. For the sake of an easier
formulation of the result, we denote by $\nabla^T$ the metric connection with
torsion $T$.

\begin{thm}[{\cite[Thm. 4.1.]{Agri&C&F05}}]\label{main-result}
%-------------------------------------------------------------
The  equation $\nabla^T\Psi=0$
admits $7$ solutions for some $3$-form, namely:
\begin{itemize}
\item[a)] A two-parameter family of pairs 
$(T_{r,s},\Psi_{r,s})\in \Lambda^3(\mathrm{Sol})\x \Sigma (\mathrm{Sol})$ 
such that $\nabla^{T_{r,s}}\Psi_{r,s}=0$;\\
  for $r=s$ the torsion $T_{r,r}=0$ and $\Psi_{r,r}$ is a multiple of 
the $\nabla^g$-parallel spinor.\smallbreak
\item[b)] Six `isolated' solutions occuring in pairs, 
$(T^\eps_i,\Psi^\eps_i)\in \Lambda^3(\mathrm{Sol})\x \Sigma(\mathrm{Sol})$ for 
$i=1,2,3$ and $\eps=\pm$.
\end{itemize}
All these $G_2$ structures admit exactly one parallel spinor, and for 
\begin{itemize}
\item[] $ |r|\neq |s|$:\ $\omega_{r,s}$ is of general type 
  $\R\oplus S^2_0\R^7\oplus \R^7$,\smallbreak
\item[] $r=s$:\ $\omega_{r,r}$ is $\nabla^g$-parallel,\smallbreak
\item[] $r=-s$:\  the $G_2$ class has no $\R$-part.
\end{itemize}
Here, $\omega_{r,s}$ denotes the defining $3$-form of the $G_2$-structure,
see eq.~$(\ref{spinor-form})$.
\end{thm}
\begin{NB}
%---------
A routine computation establishes that 
$\langle dT\cdot\Psi,\Psi\rangle < 0$ for all  solutions found in
Theorem \ref{main-result}, except for the integrable case $r=s$ of
solution a) where it vanishes trivially since $T=0$. 
\end{NB}
\begin{NB}
%---------
The interaction between explicit Riemannian 
metrics with holonomy $G_2$ on non-compact manifolds and the non-integrable
$G_2$-geometries as investigated with the help of connections with torsion
was up-to-now limited to ``cone-type arguments'', i.\,e.~a
non-integrable structure on some manifold was used to
construct an integrable structure on a higher dimensional manifold
(like its cone, an so on). It is thus a  natural question  whether the same
Riemannian manifold $(M,g)$ can carry structures of both type 
\emph{simultaneously}. This appears to be a remarkable property, of which the
above example is the only known instance. 
To emphasize this, consider that the projective space $\C
P^3$ with the well-known K\"ahler-Einstein structure and the nearly
K\"ahler one inherited from
triality does not fit the picture, as they refer to different metrics.
\end{NB}
\subsection{Naturally reductive spaces and Kostant's cubic Dirac operator}\label{kostant}
%------------------------------------------------------------------------- 
%
On arbitrary manifolds, only Weitzenb\"ock formulas that express $D^2$
through the Laplacian are available. On homogeneous spaces, it makes
sense to look for expressions for $D^2$ of \emph{Parthasarathy type},
that is, in terms of Casimir operators. 
Naturally reductive spaces $M^n=G/H$ with their family of metric 
connections ($X,Y\in\m$)
\bdm
\nabla^t_X Y\ :=\ \nabla^g_X Y -\frac{t}{2}[X,Y]_\m
\edm
were in fact investigated prior to the more general case described in
the previous section. As symmetric spaces are good toy models for
integrable geometries, homogeneous non-symmetric spaces are a very
useful field for `experiments' in non-integrable geometry.
Furthermore, many examples of such geometries are in fact homogeneous. 
We will show that
the main achievement in \cite{Kostant99} was to realize that, for
the parameter value $t=1/3$, the square of $D^t$ may be expressed in
a very simple way in terms of Casimir operators and scalars only
(\cite[Thm 2.13]{Kostant99}, \cite[10.18]{Sternberg99}). It is a remarkable
generalization of the classical Parthasarathy formula for $D^2$ on
symmetric spaces (formula~(\ref{Kos-Parth-D2-symm}) in this article, see
\cite{Parthasarathy72}). We shall speak of  the \emph{generalized 
Kostant-Parthasarathy formula} in the sequel.
S.\ Slebarski  used the connection $\nabla^{1/3}$  to prove a "vanishing 
theorem" for the kernel of the twisted Dirac operator, which can be easily
recovered from Kostant's formula (see \cite[Thm 4]{Landweber00}).
His articles \cite{SlebarskiI87} and \cite{SlebarskiII87} contain
several formulas of Weitzenb\"ock type for $D^2$, but none of them is of 
Parthasarathy type.

In order to exploit the full power of harmonic analysis, it is necessary
to extend the naturally reductive metric $\lan\ , \ \ran$ on $\m$ to the whole
Lie algebra $\g$ of $G$. By a classical theorem of B.~Kostant, there exists
a unique $\Ad(G)$ invariant, symmetric,
non degenerate, bilinear form $Q$ on $\g$  such that 
 \bdm
 Q(\h\cap\g,\m)\ =\ 0\ \text{ and }\ Q|_{\m}\ =\ \lan\ , \ \ran\, 
 \edm
if $G$ acts effectively on $M^n$ and $\g=\m+[\m,\m]$, which we will
tacitly assume from now on \cite{Kostant56}. In general, $Q$ does not have to
be positive definite; if it is, the metric is called 
\emph{normal homogeneous}. Assume furthermore that 
there exists a homogeneous spin structure on $M$, i.\ e.,
a lift $\Adtilde:\ H\ra\Spin(\m)$ of the isotropy representation such 
that the diagram
%
%-------------------------------------------
\begin{diagram}
  &                       &  \Spin(\m)\\
  &\ruTo^{\Adtilde} & \uTo_{\lambda}\\
H & \rTo^{\Ad} & \SO(\m)\\
\end{diagram}
%-------------------------------------------
%
commutes, where $\lambda$ denotes the spin covering. Moreover,   
denote by $\adtilde$ the corresponding lift into $\spin(\m)$ of
the differential $\ad:\h\ra\so(\m)$ of $\Ad$. Let
$\kappa:\, \Spin(\m)\ra \GL(\Delta_{\m})$ be the spin representation,
and identify sections of the spinor bundle 
$\Sigma M^n=G\x_{\Adtilde}\Delta_{\m}$
with functions $\psi:\ G\ra\Delta_{\m}$ satisfying
 \bdm 
 \psi(gh)\ =\ \kappa(\Adtilde(h^{-1}))\psi(g)\,.
 \edm
The Dirac operator takes for $\psi\in \Sigma M^n$ the form
\bdm
D^t\psi\ =\ \sum_{i=1}^n e_i(\psi)+ \frac{1-t}{2}\, H\cdot\psi,
\edm
where $H$ is the third degree element in the Clifford algebra  
$\mathrm{Cl}(\m)$ of $\m$ induced from the torsion,
\bdm
H\ :=\ \frac{3}{2}\sum_{i<j<k} \lan [e_i,e_j]_\m,e_k\ran\,
e_i\cdot e_j\cdot e_k.
\edm 
This fact 
suggested  the name "cubic Dirac operator" to B.\ Kostant.
Two expressions appear over and over again for naturally reductive spaces:
these are the $\m$- and $\h$-parts of the Jacobi identity,
 \bea[*]
 \Jacm(X,Y,Z)& :=& [X,[Y,Z]_{\m}]_{\m}+[Y,[Z,X]_{\m}]_{\m}+
 [Z,[X,Y]_{\m}]_{\m}\, ,\\
 \Jach(X,Y,Z)& :=& [X,[Y,Z]_{\h}]\,+\,[Y,[Z,X]_{\h}]\,+\, [Z,[X,Y]_{\h}]\,.
 \eea[*]
Notice that the summands of $\Jach(X,Y,Z)$ automatically lie in $\m$
by the assumption that $M$ is reductive. The Jacobi identity for $\g$
implies $\Jacm(X,Y,Z)+\Jach(X,Y,Z) =0$. In fact, since
the torsion is given by $T^t(X,Y)=-t[X,Y]_\m$, one immediately sees that
$\langle\Jacm(X,Y,Z),V\rangle$ is just $-\sigma_{T^t}(X,Y,Z,V)$ as defined before.
From the explicit formula for $T^t$ and the property $\nabla^1 T^1=0$, 
it is a routine computation to show that (see \cite[Lemma 2.3, 2.5]{Agri})
\bdm
\nabla^t_Z T^t(X,Y)\ =\ \frac{1}{2}t(t-1)\Jacm(X,Y,Z),\quad
dT^t(X,Y,Z,V)\ =\ -2 t \lan\Jacm(X,Y,Z),V\ran.  
\edm
In particular, $dT^t$ and $\sigma_{T^t}$ are always proportional (see 
Remark~\ref{NB-dT-prop-sigma}). The first formula implies 
$X\haken \nabla^t_X T^t=0$, hence $\delta^t T^t=0$ and it equals 
the Riemannian divergence $\delta^g T^t$ by Proposition \ref{delta-form}
of the Appendix. Since the $\Ad(G)$-invariant extension $Q$ of
$\lan\ , \ \ran$ is not necessarily positive definite when restricted to $\h$, 
it is more appropriate to work with dual rather than with orthonormal bases.
So  pick bases $x_i,y_i$ of $\h$ wich are dual with respect to $Q_{\h}$, 
i.\ e., $Q_{\h}(x_i,y_j)=\delta_{ij}$. 
 The (lift into the spin bundle of the) Casimir operator of the full Lie 
algebra $\g$ is now the sum of 
a second order differential operator (its $\m$-part) and a constant
element of the Clifford algebra  (its $\h$-part) 
 \bdm
 \Omega_{\g}(\psi) \ =\ - \sum_{i=1}^n e^2_i(\psi)  - 
\sum_{j=1}^{\dim\h}\adtilde(x_j)\circ \adtilde(y_j)\cdot\psi \ \ \text{ for }
\psi\in\Sigma M^n\,.
 \edm
In order to prove the generalized Kostant-Parthasarathy formula for the 
square of $D^t$,
similar technical prerequisites as in Section \ref{square-dirac} are needed,
but now expressed with respect to representation theoretical quantities
instead of analytical ones. We refer to \cite{Agri} for details and
will rather formulate the final result without detours. Observe that
the dimension restriction below ($n\geq 5$) is not essential, for
small dimensions a similar formula holds, but it looks slightly different.
\begin{thm}[Generalized Kostant-Parthasarathy formula, {\cite[Thm 3.2]{Agri}}]
\label{K-P-1}
%------------------------------------------------------------------------------
For $n\geq 5$, the square of $D^t$ is given by
 \bea[*]
 (D^t)^2\psi & =& \Omega_{\g}(\psi) + \frac{1}{2}(3t-1) \sum_{i,j,k}
 \lan [e_i,e_j]_{\m},e_k \ran e_i\cdot e_j\cdot e_k(\psi)\\
 &-& \frac{1}{2}\sum_{i<j<k<l}\langle e_i, \Jach(e_j,e_k,e_l)+ 
{\scriptsize \frac{9(1-t)^2}{4}}
 \Jacm(e_j,e_k,e_l)\rangle\cdot e_i\cdot e_j\cdot e_k\cdot e_l\cdot\psi \\
 &+&  \frac{1}{8}\sum_{i,j} Q_{\h}([e_i,e_j],[e_i,e_j])\psi
 +\frac{3(1-t)^2}{24} \sum_{i,j} Q_{\m}([e_i,e_j],[e_i,e_j])\psi\,.
 \eea[*]
\end{thm}
Qualitatively, this result is similar to equation (\ref{FI-SL}) of
Theorem \ref{FI-SL-AC}, although one cannot be deduced directly from the other.
Again, the square of the Dirac operator is written as the sum of a second order
differential operator (the Casimir operator), a first order differential 
operator, a four-fold product in the Clifford algebra and a scalar part
(recall that $\delta^t T^t=0$, hence this term has no counterpart here).
An immediate consequence is the special case $t=1/3$:
\begin{cor}[The Kostant-Parthasarathy formula for $t=1/3$]\label{K-P-2}
%----------------------------------------------------------------------
For  $n\geq 5$ and $t=1/3$, the general formula for $(D^t)^2$ 
reduces to
 \bea[*]
 (D^{1/3})^2\psi & = & \Omega_{\g}(\psi)+ \frac{1}{8} \bigg[ 
 \sum_{i,j} Q_{\h}([e_i,e_j],[e_i,e_j]) + 
 \frac{1}{3}\sum_{i,j} Q_{\m}([e_i,e_j],[e_i,e_j]) \bigg]\psi \\
 & =& \Omega_{\g}(\psi) + \frac{1}{8}
 \bigg[\Scal^{1/3}+ \frac{1}{9}\sum_{i,j}Q_{\m}([e_i,e_j],[e_i,e_j])\bigg]
 \psi\,.
 \eea[*]
\end{cor}
\begin{NB}
%---------
In particular, one  immediately recovers the classical Parthasarathy formula
for a symmetric space, since then 
all scalar curvatures coincide and $[e_i,e_j]\in \h$. In fact, 
compared with Theorem \ref{new-weitzenboeck}, Corollary \ref{K-P-2} 
has the advantage of containing no $4$-form action on the spinor
and the draw-back that the Casimir operator of a naturally reductive
space is not necessarily a non-negative operator (see Section \ref{Casimir}
for a detailed investigation of this point).
\end{NB}
As in the classical Parthasarathy formula, the scalar term as well as the
eigenvalues of $\Omega_{\g}(\psi)$ may be expressed
in representation theoretical terms if $G$  (and hence $M$) is compact. 
\begin{lem}[{\cite[Lemma 3.6]{Agri}}]\label{comp-scalar-in-K-P}
%-------------------------------------------------------------
Let $G$ be compact, $n\geq 5$, and denote by $\vrho_{\g}$
and $\vrho_{\h}$ the half sum of the positive roots of $\g$ and $\h$, 
respectively. Then the Kostant-Parthasarathy
formula for $(D^{1/3})^2$ may be restated as
 \bdm
 (D^{1/3})^2\psi \ = \ \Omega_{\g}(\psi)+ \left[Q( \vrho_{\g},\vrho_{\g})-
 Q(\vrho_{\h},\vrho_{\h})\right]\psi \ =\ \Omega_{\g}(\psi)+
 \lan \vrho_{\g}- \vrho_{\h},\vrho_{\g}- \vrho_{\h}  \ran
 \psi\,.
 \edm
In particular, the scalar term is positive independently
of the properties of $Q$.
\end{lem}
We can formulate our first conclusion from Corollary \ref{K-P-2}:
\begin{cor}[{\cite[Cor.\,3.1]{Agri}}]\label{eigenvalue-estimate}
%----------------------------------------------------------------
Let $G$ be compact. If the operator $\Omega_{\g}$ is non-negative,  the 
first eigenvalue
$\lambda_1^{1/3}$ of the Dirac operator $D^{1/3}$ satisfies the inequality
 \bdm
 \big(\lambda_1^{1/3}\big)^2\ \geq \ Q( \vrho_{\g},\vrho_{\g})-
 Q(\vrho_{\h},\vrho_{\h}) \,.
 \edm
Equality occurs if and only if there exists an algebraic spinor in
$\Delta_{\m}$  which is fixed under the lift $\kappa(\Adtilde H)$ of the isotropy
representation.
\end{cor}
\begin{NB}
%---------
This eigenvalue estimate is remarkable for several reasons. Firstly, for
homogeneous non symmetric spaces, it is sharper than the classical
Parthasarathy formula.  For a symmetric space, one classically obtains 
$\lambda_1^2\geq \Scal/8$. But since the Schr\"odinger-Lichnerowicz formula
yields immediately $\lambda_1^2\geq \Scal/4$, the lower bound in the
classical Parthasarathy formula is never attained. In contrast, there exist 
many examples of homogeneous non symmetric 
spaces with constant spinors. Secondly, it uses a lower bound wich
is \emph{always} strictly positive; for many naturally reductive metrics
with negative scalar curvature a pure curvature bound would  be of 
small interest. Our previously discussed generalizations of
the Schr\"odinger-Lichnerowicz yield no immediate eigenvalue estimate.
S.\,Goette derived in \cite[Lemma 1.17]{Goette99} an
eigenvalue estimate for normal homogeneous naturally reductive metrics,
but it is also not sharp.
\end{NB}
\begin{NB}\label{D-new-op}
%-------------------------
Since $D^t$ is a $G$-invariant differential operator on $M$ by construction,
Theorem \ref{K-P-1} implies that the linear combination
of the first order differential operator and the multiplication
by the element of degree four in the Clifford algebra
appearing in the formula for $(D^t)^2$ is again $G$ invariant for all
$t$. Hence, the first order differential operator
\bdm
\tilde{\mathcal{D}}\psi \ :=\ \sum_{i,j,k} \lan [Z_i,Z_j]_{\m},Z_k \ran 
Z_i\cdot Z_j\cdot Z_k(\psi)
\edm
has to be a  $G$-invariant differential operator, a fact that cannot be seen
directly by any simple arguments. It has no analogue on symmetric spaces and
certainly deserves  further separate investigations.
Of course, it should be understood as  a `homogeneous cousin'
of the more general operator $\mathcal{D}$ defined in equation 
(\ref{operator-D}).
\end{NB}
\subsection{A Casimir operator for characteristic connections}
\label{Casimir}
%-------------------------------------------------------------------------
%
Typically, the canonical connection of a naturally reductive homogeneous 
space $M$ can be given an alternative geometric 
cha\-rac\-teri\-zation---for example,
as the unique metric connection with skew-symmetric torsion preserving a given 
$G$-structure. Once this is done,  $D^{1/3}$, $\Scal^g$ and $||\T||^2$ are 
geometrically invariant objects, whereas $\Omega_{\g}$ still heavily relies on 
the concrete realization of the homogeneous space $M$ as a quotient. 
At the same time, the same interesting $G$-structures exist on many
non-homogeneous manifolds. Hence it was our goal to find a tool similar to
$\Omega_{\g}$ which has more intrinsic geometric meaning and which can
be used in both situations just described \cite{AgriFri2}.

We consider a Riemannian spin manifold $(M^n, g, \nabla)$ 
with a metric connection
$\nabla$ and skew-symmetric torsion $T$. Denote by 
$\Delta_{T}$ the spinor Laplacian of the connection $\nabla$.
\begin{dfn} 
%----------
The \emph{Casimir operator} of  $(M^n,g,\nabla)$ is
the differential operator acting on spinor fields by
\begin{eqnarray*}
\Omega & := & (D^{1/3})^2 \, + \, \frac{1}{8} \, (d T \, - \,
2 \, \sigma_{T}) \, + \, \frac{1}{4} \, \delta(T) 
\, - \,  \frac{1}{8} \, \Scal^{g} \, - \, \frac{1}{16} \, ||\T||^2 \\
& = & \Delta_{T} \, + \, \frac{1}{8} \, ( 3 \, d T \, 
- \, 2 \, \sigma_{T} \, + \, 2 \, \delta(T) \, + \,
\Scal ) \, .
\end{eqnarray*} 
\end{dfn}
\begin{NB} 
%---------
A naturally reductive space $M^n = G/H$ endowed with its canonical  
connection satisfies $dT=2\sigma_{T}$ and $\delta T=0$, hence
$\Omega=\Omega_{\g}$ by  Theorem~\ref{K-P-2}. For connections with
$d T\neq 2\sigma_{T}$ and $\delta T\neq 0$, the numerical factors are
chosen in such a way to yield an overall expression proportional
to the scalar part of the right hand side of equation~(\ref{FI-SL}).
\end{NB} 
\begin{exa} 
%----------
For the Levi-Civita connection $(T=0)$ of an arbitrary
Riemannian manifold, we obtain
\bdm
\Omega \ = \ (D^{g})^2 \, - \, \frac{1}{8} \, \Scal^{g} \ = \  
\Delta^{g} \, + \, \frac{1}{8} \, \Scal^{g} \, .
\edm
The second equality is just the classical Schr\"odinger-Lichnerowicz
formula for the Riemannian Dirac operator, whereas the first one is---in 
case of a symmetric space---the classical Parthasarathy formula.
\end{exa} 
\begin{exa} 
%-----------
Consider a $3$-dimensional manifold of constant scalar
curvature, a constant $a \in \R$ 
and the $3$-form $T = 2 \, a \, dM^3$. Then
\bdm
\Omega \ = \ (D^{g})^2 \, - \, a \, D^{g} \, - \, 
\frac{1}{8} \, \Scal^{g}  .
\edm 
The kernel of the Casimir
operator corresponds to eigenvalues $\lambda \in \mathrm{Spec}(D^g)$
of the Riemannian Dirac operator such that
\bdm
8 \, (\lambda^2 \, - \, a \, \lambda) \, - \, \Scal^g 
\ = \ 0 \, . 
\edm 
In particular, the kernel of $\Omega$ is in general larger then the
space of $\nabla$-parallel spinors. Indeed, such spinors exist
only on space forms. More generally, fix a real-valued smooth
function $f$ and consider the $3$-form $T := f \cdot dM^3$. If
there exists a $\nabla$-parallel spinor
\bdm
\nabla^g_X \psi \, + \, (X \haken T) \cdot \psi \ = \ 
\nabla^g_X \psi \, + \, f \cdot X \cdot \psi \ = \ 0 \, ,
\edm
then, by a theorem of A.~Lichnerowicz (see \cite{Li}), $f$ is constant and
$(M^3,g)$ is a space form.
\end{exa} 
Let us collect some elementary properties of the Casimir operator.
\begin{prop}[{\cite[Prop. 3.1]{AgriFri2}}]
%-----------------------------------------
The kernel of the Casimir operator contains all $\nabla$-parallel spinors. 
\end{prop}
\begin{proof}
By Theorem~\ref{FI-SL-AC},
one of the integrability conditions for a $\nabla$-parallel spinor field 
$\psi$ is
\bdm
\big(3 \, d T \, 
- \, 2 \, \sigma_{T} \, + \, 2 \, \delta(T) \, + \,
\Scal \big) \cdot \psi \ = \ 0 \, .   \qedhere
\edm
\end{proof} 
\noindent
If the torsion form $T$ is $\nabla$-parallel,  the formulas for the
Casimir operator simplify. Indeed, in this case we have 
(see the Appendix) 
\bdm
d \T \ = \ 2 \, \sigma_{\T} , \quad \delta(\T) \ = \ 0  \, , 
\edm
and the Ricci tensor $\mathrm{Ric}$ of $\nabla$ is symmetric.
Using the formulas of Section~\ref{square-dirac} (in particular, Theorems
\ref{FI-SL-AC} and \ref{new-weitzenboeck}), we obtain a simpler expressions 
for the Casimir operator. 
\begin{prop}[{\cite[Prop. 3.2]{AgriFri2}}]\label{casimir} 
%---------------------------------------------------------
For a metric connection with parallel torsion ($\nabla T=0$),
the Casimir operator  can equivalently be written as: 
\begin{eqnarray*}
\Omega & = & (D^{1/3})^2 \, - \, \frac{1}{16} \, \big(2 \, \Scal^{g} \, + \, 
||T||^2 \big) \
= \ \Delta_{T} \, + \, \frac{1}{16} \, \big(2 \, \Scal^{g} \, + \, ||T||^2 
\big) \, - \, \frac{1}{4} \, T^2\\
& = & \Delta_{T} \, + \, \frac{1}{8} \, \big( 2 \, dT \, + \, 
\Scal \big)  .
\end{eqnarray*}
\end{prop}
\noindent
Integrating these formulas, we obtain a vanishing theorem for
the kernel of the Casimir operator.
\begin{prop}[{\cite[Prop. 3.3]{AgriFri2}}]\label{prop1}
%------------------------------------------------------
Assume that $M$ is compact and that $\nabla$ has parallel torsion $T$.
If one of the conditions 
\bdm
2 \, \Scal^g \ \leq \ - \, ||\T||^2 \quad \mathrm{or} \quad 
2 \, \Scal^g \ \geq \ 4 \, \T^2 \, - \, ||\T||^2,
\edm
holds,  the Casimir operator is non-negative in $\mathrm{L}^2(S)$.
\end{prop}
\begin{exa}
For a naturally reductive space $M=G/H$, the first condition can never hold,
since by Lemma \ref{comp-scalar-in-K-P}, 
 $2\, \Scal^g+||\T||^2$ is strictly positive. In concrete examples,
the second condition typically singles out the normal homogeneous 
metrics among the naturally reductive ones. 
\end{exa}
\begin{prop}[{\cite[Prop. 3.4]{AgriFri2}}]\label{prop2}
%------------------------  ----------------------------
If the torsion form is $\nabla$-parallel,  the Casimir operator 
$\Omega$ and the square of the 
Dirac operator $(D^{1/3})^2$ commute with the endomorphism $\T$, 
\bdm
\Omega \circ \T \ = \ \T \circ \Omega \, , \quad (D^{1/3})^2 \circ \T \ 
= \ \T \circ (D^{1/3})^2 .
\edm
\end{prop}
\noindent
The endomorphism $\T$ acts on the spinor bundle as a symmetric endomorphism
with {\it constant} eigenvalues.
\begin{thm}\label{kernel}  
%------------------------
Let $(M^n,g,\nabla)$ be a compact Riemannian spin manifold
equipped with a metric connection $\nabla$ with parallel, skew-symmetric
torsion, $\nabla T = 0$. The endomorphism $T$ and the Riemannian
Dirac operator $D^g$ act
in the kernel of the Dirac operator $D^{1/3}$. In particular, 
if, for all $\mu \in 
\mathrm{Spec}(T)$, the number $- \, \mu/4$ is not an eigenvalue of the
Riemannian Dirac operator, then the kernel of $D^{1/3}$ is trivial.
\end{thm}
\noindent
If $\psi$ belongs to the 
kernel of $D^{1/3}$ and is an eigenspinor of the endomorphism $T$, 
we have $4 \cdot D^g \psi = - \, \mu \cdot \psi$ , $\mu \in 
\mathrm{Spec}(T)$.  
Using the estimate of the eigenvalues of the Riemannian Dirac operator
(see \cite{Friedrich80}),
we obtain an upper bound for the minimum $\Scal_{\mathrm{min}}^g$ Riemannian 
scalar curvature in case that the kernel of the operator 
$D^{1/3}$ is non trivial.
\begin{prop}\label{prop3}    
Let $(M^n,g,\nabla)$ be a compact Riemannian spin manifold
equipped with a metric connection $\nabla$ with parallel, skew-symmetric
torsion, $\nabla T = 0$. If the kernel of the Dirac operator $D^{1/3}$ is
non trivial, then the minimum of the Riemannian scalar curvature is bounded
by
\bdm
\max\big\{ \mu^2 \, : \, \mu \in \mathrm{Spec}(T) \big\} \ \geq \ 
\frac{4n}{n-1} \, \Scal_{\mathrm{min}}^g \, .
\edm 
\end{prop}
\begin{NB}
%--------
If $ (n-1) \, \mu^2 = 4n \, \Scal^g$ is in the spectrum
of $T$ and there exists a spinor field $\psi$ in the kernel of $D^{1/3}$
such that $T \cdot \psi = \mu \cdot \psi$, then we are in the limiting case 
of the inequality
in \cite{Friedrich80}. Consequently, $M^n$ is an Einstein manifold of 
non-negative
Riemannian scalar curvature and $\psi$ is a Riemannian Killing spinor.
Examples of this type are $7$-dimensional $3$-Sasakian manifolds.
The possible torsion forms have been discussed in \cite{Agricola&F03a}, 
Section $9$.
\end{NB}
We discuss in detail what happens for $5$-dimensional Sasakian manifold.
Let $(M^5, g, \xi, \eta, \varphi)$ be a (compact) $5$-dimensional
Sasakian spin manifold with a fixed spin structure,  $\nabla^c$ 
its characteristic connection.
We orient $M^5$ by the condition that the differential of the contact
form is given by $d\eta=2(e_1\wedge e_2+e_3\wedge e_4)$, and write
henceforth $e_{ij\ldots}$ for $e_i\wedge e_j\wedge \ldots$. 
Then we know that
\bdm
\nabla T^c \ = \ 0  , \quad T^c \ = \ \eta \wedge d \eta \ 
= \ 2 \, (e_{12} \, + \, e_{34}) \wedge
e_5  , \quad (T^c)^2 \ = \ 8 \, - \, 8 \, e_{1234} 
\edm
and 
\bdm
\Omega \ = \ (D^{1/3})^2 \, - \, \frac{1}{8} \, \Scal^g \, 
- \, \frac{1}{2} \ = \ \Delta_{T^c} \, + \, \frac{1}{8} \,
\Scal^g \, - \, \frac{3}{2} \, + \, 2 \, e_{1234} .
\edm
We study the kernel of the Dirac operator $D^{1/3}$.
The endomorphism $ T^c$ acts in the $5$-dimensional spin representation
with eigenvalues $(-4,0,0,4)$ and, according to Theorem
\ref{kernel}, we have to distinguish two cases. If $D^{1/3} \psi = 0$ and 
$T^c \cdot \psi = 0$, the
spinor field is harmonic and the formulas of
Proposition \ref{casimir} yield in the compact case the condition
\bdm
\int_{M^5} \big( 2 \, \Scal^g \, + \, 8 \big) \, ||\psi||^2 \ 
\leq \ 0  .
\edm 
Examples of that type are
the $5$-dimensional Heisenberg group with its left invariant Sasakian
structure and its compact quotients (Example \ref{Heisenberg-Sasaki}) or 
certain $S^1$-bundles over a flat torus. 
The space of all $\nabla$-parallel spinors satisfying
$T^c \cdot \psi_0 = 0$ is a $2$-dimensional subspace of the kernel of 
the operator $D^{1/3}$ (see \cite{Friedrich&I1}, \cite{Friedrich&I2}). 
The second case for spinors in the kernel is given by $D^{1/3} \psi = 0$ 
and $T^c \cdot \psi =  \pm 4 \, \psi$. The spinor field is an eigenspinor for 
the Riemannian Dirac operator, $D^g \psi = \mp \, \psi$.  
The formulas of
Proposition \ref{casimir} and Proposition \ref{prop3} yield in the compact 
case two conditions:
\bdm
\int_{M^5} \big(\Scal^g \, - \, 12 \big) ||\psi||^2 \ 
\leq \ 0 \quad \quad \mathrm{and} \quad \quad 5 \, 
\Scal_{\mathrm{min}}^g \ \leq \ 16 .
\edm 
The paper \cite{FK} contains a construction of Sasakian manifolds admitting 
a spinor field of that algebraic type in the kernel of $D^{1/3}$. 
We describe the construction explicitly. Suppose that the 
Riemannian Ricci tensor $\Ric^g$ of a simply-connected, $5$-dimensional 
Sasakian manifold is given by the formula
\bdm
\Ric^g \ = \ - \, 2 \cdot g \, + \, 6 \cdot \eta \otimes \eta \, .
\edm
Its scalar curvature 
equals $\Scal^g = -\, 4$. In the simply-connected and compact case, 
they are total spaces of $S^1$ principal bundles over $4$-dimensional
Calabi-Yau orbifolds (see \cite{BG}).
There exist (see \cite{FK}, Theorem 6.3) two spinor fields $\psi_1$, 
$\psi_2$ such that
\bdm
\nabla^g_X \psi_1 \ = \ - \, \frac{1}{2} \, X \cdot \psi_1 \, + \, 
\frac{3}{2} \, \eta(X) \cdot \xi \cdot \psi_1 \, , \quad 
\T \cdot \psi_1 \ = \ - \, 4 \, \psi_1 ,
\edm
\bdm
\nabla^g_X \psi_2 \ = \  \frac{1}{2} \, X \cdot \psi_2 \, - \, 
\frac{3}{2} \, \eta(X) \cdot \xi \cdot \psi_2 \, , \quad 
T^c \cdot \psi_2 \ = \ 4 \, \psi_2 .
\edm
In particular, we obtain 
\bdm
D^g \psi_1 \ = \ \psi_1 \, , \ T^c \cdot \psi_1 \ = \ 
- \, 4 \, \psi_1 \, , \quad \mathrm{and} \quad
D^g \psi_2 \ = \ - \, \psi_2 , \  T^c \cdot \psi_2 \ = \ 4 \, \psi_2 ,
\edm
and therefore the spinor fields $\psi_1$ and $\psi_2$ belong to the kernel of the operator $D^{1/3}$. 
\vspace{2mm} 

\noindent
Next, we investigate the kernel of the Casimir operator. Under the action 
of the torsion form, the spinor bundle $\Sigma M^5$ splits into three
subbundles $\Sigma M^5 = \Sigma_0 \oplus \Sigma_{4} \oplus 
\Sigma_{-4}$ corresponding to the eigenvalues of $T^c$. Since
$\nabla T^c = 0$, the connection $\nabla$ preserves the splitting.
The endomorphism $e_{1234}$ acts by the formulas
\bdm
e_{1234} \ = \ 1 \quad \text{on }\  \Sigma_0, \quad e_{1234} \ = \ - \, 1 
\quad \text{on } \ \Sigma_4 \oplus \Sigma_{-4}.
\edm
Consequently, the formula 
\bdm
\Omega \ = \ \Delta_{T^c} \, + \, \frac{1}{8} \,
\Scal^g \, - \, \frac{3}{2} \, + \, 2 \, e_{1234}
\edm
shows that the Casimir operator splits into the sum $\Omega = \Omega_0 
\oplus \Omega_4 \oplus \Omega_{-4}$ of three operators 
acting on sections in $\Sigma_0$, $\Sigma_4$ and $\Sigma_{-4}$. 
On $\Sigma_0$, we have
\bdm
\Omega_0 \ = \ \Delta_{T^c} \, + \, \frac{1}{8} \, \Scal^g 
\, + \, \frac{1}{2} \ = \ (D^{1/3})^2 \, - \, 
\frac{1}{8} \, \Scal^g \, - \, \frac{1}{2} \, .
\edm
In particular, the kernel of $\Omega_0$ is trivial if 
$\Scal^g \neq -4$. The 
Casimir operator on $\Sigma_4 \oplus \Sigma_{-4}$ is given by
\bdm
\Omega_{\pm 4} \ = \ \Delta_{T^c} \, + \, 
\frac{1}{8} \, \Scal^g 
\, - \, \frac{7}{2} \ = \ (D^{1/3})^2 \, - \, 
\frac{1}{8} \, \Scal^g \, - \, \frac{1}{2} 
\edm
and a non trivial kernel can only occur if $-4 \leq \Scal^g \leq 28$.
A spinor field $\psi$ in the kernel of
the Casimir operator $\Omega$ satisfies the equations
\bdm
(D^{1/3})^2 \cdot \psi \ = \ \frac{1}{8} \, (4 \, + \, \Scal^g ) \, \psi
\, , 
\quad  T^c\cdot \psi \ = \ \pm \, 4 \, \psi \, .
\edm
In particular, we obtain
\bdm
\int_{M^5} \langle (D^g \, \pm \, 1)^2 \, \psi \, , \, \psi \rangle \ = \ 
\frac{1}{8} \int_{M^5} \big(4 \, + \, \Scal^g \big) \, ||\psi||^2 ,
\edm
and the first eigenvalue of the operator $(D^g \pm 1)^2$ is bounded by the
scalar curvature,
\bdm
\lambda_1(D^g \, \pm \, 1)^2 \ \leq \ \frac{1}{8} \, \big( 4 \, + \, 
\Scal^g_{\mathrm{max}} \big)  . 
\edm
Let us consider special classes of Sasakian manifolds. A first 
case is $\Scal^g = - \, 4$. Then the formula for the Casimir
operator simplifies, 
\bdm
\Omega_0 \ = \ \Delta_{T^c} \ = \ (D^{1/3})^2 \, , \quad 
\Omega_{\pm 4} \ = \ \Delta_{\mathrm{T}} \, - \, 4 \ = \ (D^{1/3})^2 .
\edm
If $M^5$ is compact, the kernel of the operator $\Omega_0$
coincides with the space of $\nabla$-parallel spinors in the bundle 
$\mathrm{S}_0$.  A spinor field $\psi$ in the kernel 
the operator
$\Omega_{\pm4}$ is an eigenspinor of the Riemannian Dirac
operator,
\bdm
D^g(\psi) \ = \ \mp \, \psi , \quad \mathrm{T} \cdot \psi \ 
= \ \pm \, 4 \, \psi\, .
\edm
Compact Sasakian manifolds admitting spinor fields in the kernel
of $\Omega_0$ 
are quotients of the $5$-dimensional Heisenberg group (see 
\cite{Friedrich&I2}, Theorem 4.1). Moreover, the $5$-dimensional Heisenberg 
group and its compact quotients 
admit spinor fields in the kernel of $\Omega_{\pm 4}$, too. 

\vspace{2mm}
\noindent
A second case is 
$\Scal^g = 28$. Then
\bdm
\Omega_0 \ = \ \Delta_{T^c} \, + \, 4 \ = \ (D^{1/3})^2 \, - \, 4 
\, , \quad \Omega_{\pm 4} \ = \ \Delta_{T^c} \ = \ (D^{1/3})^2 \, 
- \, 4  .
\edm
The kernel of $\Omega_0$ is trivial and the kernel of $\Omega_{\pm4}$ coincides
with the space of $\nabla$-parallel spinors in the bundle $\mathrm{S}_{\pm4}$.
Sasakian manifolds admitting spinor fields of that type have been described
in \cite{Friedrich&I1}, Theorem 7.3 and Example 7.4.
\vspace{2mm} 

\noindent
If $ -\, 4 < \mathrm{Scal}^g < 28$, the kernel of the operator $\Omega_0$ is
trivial and the kernel of $\Omega_{\pm4}$ depends on the geometry
of the Sasakian structure. Let us discuss Einstein-Sasakian 
manifolds. Their scalar curvature equals
$ \Scal^g = 20$ and the Casimir operators are
\bdm
\Omega_0 \ = \ \Delta_{T^c} \, + \, 3 \, , \quad 
\Omega_{\pm 4} \ = \ \Delta_{T^c} \, - \, 1 \ 
= \ (D^{1/3})^2 \, - \, 3  . 
\edm
If $M^5$ is simply-connected, there exist two Riemannian 
Killing spinors (see \cite{FK})
\bdm
\nabla^g_X \psi_1 \ = \ \frac{1}{2} \, X \cdot \psi_1 , 
\quad D^g(\psi_1) \ = \ - \, \frac{5}{2} \, \psi_1 , \quad 
T^c \cdot \psi_1 \ = \ 4 \, \psi_1 , 
\edm
\bdm
\nabla^g_X \psi_2 \ = \ - \, \frac{1}{2} \, X \cdot \psi_2  , 
\quad D^g(\psi_2) \ = \ \frac{5}{2} \, \psi_2 , \quad 
T^c \cdot \psi_2 \ = \  - \, 4 \, \psi_2  .
\edm
We compute the Casimir operator
\bdm
\Omega(\psi_1) \ = \ - \, \frac{3}{4} \, \psi_1 , \quad
\Omega(\psi_2) \ = \ - \, \frac{3}{4} \, \psi_2 .
\edm
In particular, the Casimir operator of a Einstein-Sasakian manifold
has {\it negative} eigenvalues. The Riemannian Killing spinors are parallel 
sections in
the bundles $\Sigma_{\pm 4}$ with respect to the
flat connections $\nabla^{\pm}$ 
\bdm
\nabla^+_X\psi \ := \ \nabla^g_X\psi \, - \, \frac{1}{2} \, X \cdot \psi 
\quad \mathrm{in} \quad \Sigma_4 , \quad 
\nabla^-_X\psi \ := \ \nabla^g_X\psi \, + \, \frac{1}{2} \, X \cdot \psi 
\quad \mathrm{in} \quad \Sigma_{-4} .
\edm
We compare these connections with our canonical connection $\nabla$: 
\bdm
\big(\nabla^{\pm}_X \, - \, \nabla_X \big) \cdot \psi^{\pm} \ 
= \ \pm \, \frac{i}{2} \, g(X, \xi) \cdot \psi^{\pm} \, , \quad \psi^{\pm}
\in \Sigma_{\pm 4}.
\edm
The latter equation means that the bundle $\Sigma_{4} \oplus 
\Sigma_{-4}$ equipped
with the connection $\nabla$ is equivalent to the $2$-dimensional 
trivial bundle with the connection form 
\bdm
\mathcal{A} \ = \  \frac{i}{2} \, \eta \cdot 
\left[\ba{cc} -1 & 0 \\ \ \ 0 & 1 \ea \right] .
\edm
The curvature
of $\nabla$ on these bundles is given by the formula
\bdm
\mathcal{R}^{\nabla} \ = \  \frac{i}{2} \, d \eta \cdot 
\left[\ba{cc} -1 & 0 \\ \ \ 0 & 1 \ea \right]  = \ 
i \, (e_{1} \wedge e_{2} \, + \, e_{3} \wedge e_{4}) \cdot
\left[\ba{cc} 1 & \ 0 \\ 0 & -1 \ea \right] .
\edm
Since the divergence $\mathrm{div}(\xi) = 0$ of the Killing vector field
vanishes, the Casimir operator on  $\Sigma_{4} \oplus 
\Sigma_{-4}$ is the following operator acting on pairs of functions:
\bdm
\Omega_{4} \oplus \Omega_{-4} \ = \ \Delta_{\mathrm{T}} \, - \, 1 \ = \ 
\Delta \, - \, \frac{3}{4} \, + \,  
\left[\ba{cc} -\, i & 0 \\ \ \ 0 & i \ea \right] \, \xi\, .
\edm
Here $\Delta$ means the usual Laplacian of $M^5$ acting on functions and $\xi$ is the differentiation in direction of the vector field $\xi$. In particular,
the kernel of $\Omega$ coincides with solutions $f : M^5 \rightarrow 
\C$ of the equation
\bdm
\Delta(f) \, - \, \frac{3}{4} \, f \, \pm \, i \, \xi (f) \ = \ 0 \ .
\edm
The $\mathrm{L}^2$-symmetric differential operators $\Delta$ and $i \, \xi$
commute. Therefore, we can diagonalize them simultaneously. The latter 
equation is solvable if and only if there exists a common eigenfunction
\bdm
\Delta(f) \ = \ \mu \, f , \quad i \, \xi (f) \ = \ \lambda 
\, f \, , \quad 4 \, (\mu \, + \, \lambda ) \, - \, 3 \ = \ 0 \, . 
\edm
The Laplacian $\Delta$ is the sum of the \emph{non-negative} horizontal
Laplacian and the operator $(i \, \xi)^2$. Now, the conditions
\bdm
\lambda^2 \ \leq \ \mu \, , \quad 4 \, (\mu \, + \, \lambda ) \, - \, 3 \ = \ 0 
\edm
restrict the eigenvalue of the Laplacian, $ 0 \leq \mu \leq 3$. On the other
side, by the Lichnerowicz-Obata Theorem, we have $ 5 \leq \mu$,
a contradiction. In particular, we proved
\begin{thm} 
%----------
The Casimir operator of a compact $5$-dimensional 
Sasaki-Einstein manifold has trivial kernel; in particular, there
are no $\nabla^c$-parallel spinors.
\end{thm}
\noindent
The same argument estimates the eigenvalues of the Casimir operator. It turns
out that the smallest eigenvalues of $\Omega$ is negative and equals $- 3/4$.
The eigenspinors are the Riemannian Killing spinors. The next eigenvalue
of the Casimir operator is at least 
\bdm
\lambda_2(\Omega) \ \geq \ \frac{17}{4} \, - \, \sqrt{5} \approx \ 2.014 .
\edm 
In the literature, similar results for almost Hermitian $6$-manifolds
and $G_2$-manifolds admitting a characteristic connection can be found.
\subsection{Some remarks on the common sector of type II superstring theory}
\label{super}
%---------------------------------------------------------------------------
%
The mathematical model discussed in the common sector of type II superstring 
theory (also sometimes referred to as type I supergravity) consists
of a Riemannian manifold $(M^n ,  g)$, a metric connection
$\nabla$ with totally skew-symmetric torsion $T$ and a 
non-trivial spinor field $\Psi$. Putting the full Ricci tensor aside for
starters and assuming that the dilaton is constant, there are three equations 
relating these objects:
\bdm\tag{$*$}
\nabla \Psi \ = \ 0 \, , \quad \delta(T) \ = \ 0 \, , \quad
T \cdot \Psi \ = \ \mu \cdot \Psi \, .
\edm
The spinor field describes the supersymmetry of the model. It has been
our conviction throughout this article that this is the most important
of the equations, as non-existence of $\nabla$-parallel spinors
implies the breakdown of supersymmetry. Yet, interesting things can be said
if looking at all equations simultaneously.
Since $\nabla$ is a metric
connection with totally skew-symmetric torsion, the divergences
$\delta^{\nabla}(T) = \delta^g(T)$ of the torsion form coincide
(see Proposition \ref{delta-form}). We denote this unique
$2$-form simply by $\delta(T)$. The third equation is an
algebraic link between the torsion form $T$ and the spinor field
$\Psi$. Indeed, the $3$-form $T$ acts as an endomorphism in the spinor
bundle and the last equation requires that $\Psi$ is an eigenspinor for this
endomorphism. Generically, $\mu=0$ in the physical model; but the mathematical
analysis becomes more transparent if we first include this parameter.
A priori,  $\mu$ may be an arbitrary function. Since $T$ acts 
on spinors as a symmetric endomorphism, $\mu$ has to be real. Moreover,
we will see that only real, constant parameters $\mu$ are possible.
Recall that the conservation
law $\delta(T) = 0$ implies that the Ricci tensor $\Ric^{\nabla}$ of the
connection $\nabla$ is symmetric, see the Appendix. Denote by 
$\Scal^{\nabla}$ the 
$\nabla$-scalar curvature and by $\Scal^g$ the scalar curvature of the
Riemannian metric. Based on the results of Section \ref{square-dirac},
the existence of the $\nabla$-parallel spinor
field yields the so called integrability conditions, 
i.\,e.~relations
between $\mu$, $T$ and the curvature tensor of the connection $\nabla$. 
\begin{thm}[{\cite[Thm 1.1.]{Agri&F&N&P05}}]\label{AlgebraIdent}
%----------------------------------------------------------------
%
Let $(M^n,g,\nabla,T,\Psi, \mu)$ be a solution of $(*)$ and assume that
the spinor field $\Psi$ is non-trivial. Then the function $\mu$ is constant
and we have 
\bdm
||T||^2 \ = \ \mu^2 \, - \, \frac{\Scal^{\nabla}}{2} \ \geq \ 0 \, , \quad
\Scal^g \ = \ \frac{3}{2} \, \mu^2 \, + \, \frac{\Scal^{\nabla}}{4} \, . 
\edm
Moreover, the spinor field $\Psi$ is an eigenspinor of the endomorphism
defined by the $4$-form $d T$,
\bdm
d T \cdot \Psi \ = \ - \,\frac{\Scal^{\nabla}}{2} \cdot \Psi \, . 
\edm
\end{thm}
Since $\mu$ has to be constant, equation 
$T \cdot \Psi = \mu \cdot \Psi$ yields: 
\begin{cor}
For all vectors $X$,  one has
\bdm
(\nabla_X T) \cdot \Psi \ = \ 0 \, .
\edm
\end{cor}
The set of equations $(*)$ is completed in the common sector of type II 
superstring theory by the condition $\Ric^\nabla =0$ and the requirement
$\mu=0$. 
In \cite{Agri}, it had been shown that  the existence of a non-trivial 
solution of this system  implies $T = 0$ on \emph{compact} 
manifolds. 
Theorem~\ref{AlgebraIdent} enables us to prove the same result
without compactness assumption and under the much weaker curvature
assumption $\Scal^{\nabla}=0$:
\begin{cor}
%-----------
Assume that there exists a spinor field $\Psi\neq 0$ satisfying the equations
$(*)$. If $\mu=0$ and $\Scal^{\nabla}=0$, the torsion form $T$ 
has to vanish.
\end{cor}
This result  underlines the strength of the algebraic identities in 
Theorem~\ref{AlgebraIdent}. Physically, this result may either show that
the dilaton is a necessary ingredient for $T\neq 0$ (while it is not for
$T=0$) or that the set of equations is too restrictive (it is
derived from a variational principle).
\begin{NB}
%---------
In the common sector of type II string theories, the "Bianchi identity"
$dT=0$ is often required in addition. It does not affect the mathematical
structure of the equations $(*)$, hence we do not include it into our 
discussion. 
\end{NB} 
On a naturally reductive space, even more is true. The generalized
Kostant-Parthasarathy formula implies for the family of connections
$\nabla^t$:

\begin{thm}[{\cite[Thm.~4.3]{Agri}}]\label{eq-3-4-no-sol}
%--------------------------------------------------------
If the operator $\Omega_{\g}$ is non-negative and if $\nabla^t$ is not the
Levi-Civita connection, there do not exist any
non trivial solutions to the equations
 \bdm
 \nabla^t\psi \ = \ 0, \quad T^t\cdot \psi=0\,.
 \edm
\end{thm}
The last equation  in type II  string theory deals with the Ricci
tensor $\Ric^{\nabla}$ of the connection. Usually one
requires for constant dilaton  that the
Ricci tensor has to vanish (see \cite{Gauntlett}). 
The result above, however, indicates that this condition may be too
strong. Understanding
this tensor as an energy-momentum tensor, it seems to be more
convenient to impose a weaker condition, namely
\bdm
\mathrm{div}(\Ric^{\nabla}) \ = \ 0 \, .
\edm
A subtle point is however the fact that there are a priori two different 
divergence operators. The first operator  $\mathrm{div}^g$
is defined by the Levi-Civita connection of the Riemannian metric, while
the second operator  $\mathrm{div}^{\nabla}$  is defined by the connection
$\nabla$. By Lemma \ref{same-div}, they coincide if $\Ric$ is 
symmetric, that is, if $\delta T=0$. This is for example satisfied if 
$\nabla T=0$. We can then prove:
\begin{cor} \label{paralleleTorsion}
%------------------------------------
Let $(M^n,g,\nabla,T,\Psi, \mu)$ be a a manifold with metric
connection defined by $T$ and assume that there exists a
spinor $0\neq \psi\in\Sigma M^n$ such that
\bdm
\nabla \psi \ = \ 0  , \quad \nabla T \ = \ 0  , \quad
T \cdot \psi \ = \ \mu \cdot \psi. 
\edm
Then all scalar curvatures are constant and the divergence of the
Ricci tensor vanishes, $\mathrm{div}(\Ric^{\nabla}) = 0$.
\end{cor} 
This is one possible way to weaken the original set of equations 
in such a way that the curvature condition follows from the other ones,
as it is the case for $T=0$---there, the existence of a $\nabla^g$-parallel
spinor implies $\Ric^g=0$. Of course, only physics can provide a definite
answer whether these or other possible replacements are `the right ones'.

Incorporating a non-constant dilaton $\Phi\in C^\infty(M^n)$ is more subtle. 
The full set of equations reads in this case
\bdm
\Ric^\nabla + \frac{1}{2} \delta T + 2\nabla^gd\Phi\ =\ 0,\quad
\delta T\ =\ 2\, \grad(\Phi)\haken T,\quad \nabla\psi\ =\ 0, \quad
(2\, d\Phi - T)\cdot \psi\ =\ 0.
\edm
In some geometries,
it is possible to interpret it as a partial conformal change of the metric.
In dimension $5$,  this allows the proof that $\Phi$ basically has
to be constant:
\begin{thm}[{\cite{Friedrich&I2}}]
%----------------------------------
Let $(M^5,g,\xi,\eta,\vphi)$ be a normal almost contact metric structure with
Killing vector field $\xi$, $\nabla^c$ its characteristic connection and
$\Phi$ a smooth function on $M^5$. If there exists a spinor field
$\psi\in\Sigma M^5$ such that
\bdm
\nabla^c \psi\ =\ 0,\quad (2\, d\Phi - T)\cdot \psi\ =\ 0,
\edm
then the function $\Phi$ is constant.
\end{thm}
In higher dimension, the picture is less clear, basically because
a clean geometric interpretation of $\Phi$ is missing.
%
%-----------------------------------------------------------------------------
\appendix\section{Compilation of remarkable identities for connections with
skew-symmetric torsion}\label{formulas}
%-----------------------------------------------------------------------------
We collect in this appendix some more or less technical formulas
that one needs in the investigation of metric connections with
skew-symmetric torsion. In order to keep this exposition readable, we
decided to gather them in a separate section.

We tried to provide at least one reference with full proofs
for every stated result; however, no claim is made whether these
are the articles where these identities appeared for the first time.
In fact, many of them have been derived and rederived by authors
when needed, some had been published earlier but the authors had not
considered it worth to publish a proof etc.

In this section, the connection $\nabla$ is normalized as
\bdm
\nabla_X Y\ =\ \nabla_X^g Y+\frac{1}{2}T(X,Y,*),\quad
\nabla_X \psi\ =\ \nabla_X^g\psi + \frac{1}{4} (X\haken T)\cdot \psi.
\edm
It then easily follows that the Dirac operators are related by
$D^\nabla = D^g + (3/4)T$.
\begin{dfn} Recall that for any $3$-form $T$, an algebraic $4$-form
$\sigma_T$ quadratic in $T$ may be defined by
$2\,\sigma_T=\sum\limits_{i=1}^n (e_i\haken T)\wedge (e_i\haken T)$,
where $e_1,\ldots, e_n$ denotes an orthonormal frame. Alternatively,
$\sigma_T$ may be written without reference to an orthonormal frame as
\bdm
\sigma_T(X,Y,Z,V)\ =\ g(T(X,Y), T(Z,V))+g(T(Y,Z),T(X,V))+g(T(Z,X),T(Y,V)).
\edm
We first encountered $\sigma_T$ in the first Bianchi identity
for metric connections with torsion $T$ (Theorem \ref{Bianchi-I}).
\end{dfn}
\begin{prop}[{\cite[Prop. 3.1.]{Agri}}]\label{3-form-square}
%-----------------------------------------------------------
Let $T$ be a $3$-form, and denote by the same symbol its
associated $(2,1)$-tensor. Then its square inside
the Clifford algebra has no contribution of degree $6$ and $2$, and its scalar
and fourth degree part are given by
\bdm
T^2_0\ =\ \frac{1}{6}\,\sum_{i,j=1}^n ||T(e_i,e_j)||^2
\ =:\ ||T||^2,\quad T^2_4\ =\ - \, 2 \cdot \sigma_{T}. 
\edm
\end{prop}

\begin{lem}[{\cite[Lemma 2.4.]{Agri}}]\label{ext-diff}
%-----------------------------------------------------
If $\omega$ is an $r$-form and $\nabla$ any connection with torsion, then
 \bea[*]
(d\omega)(X_0,\,\ldots,X_r)& =& \sum_{i=0}^r (-1)^i (\nabla_{X_i}\omega)
(X_0,\ldots,\hat{X}_i,\ldots,X_r) \\
&- & \sum_{0\leq i<j\leq r}(-1)^{i+j} 
\omega(T(X_i,X_j),X_0,\ldots,\hat{X}_i,\ldots,\hat{X}_j,\ldots,X_r)\,.
 \eea[*]
\end{lem}
\begin{cor}[{\cite{Ivanov&P01}}]
%--------------------------------
For a metric connection $\nabla$ with torsion $T$,
the exterior derivative of $T$ is given by
\bdm
dT(X,Y,Z,V)\ =\ \stackrel{X,Y,Z}{\sigma}\left[(\nabla_X T)(Y,Z,V) \right]
-(\nabla_V T)(X,Y,Z)+2\,\sigma_T(X,Y,Z,V).
\edm
In particular, $dT=2\sigma_T$ if $\nabla T=0$.
\end{cor}
\begin{prop}[{\cite[Prop. 5.1.]{Agricola&F03a}}]\label{delta-form}
%-----------------------------------------------------------------
Let $\nabla$ be a connection with skew-symmetric torsion
and define the $\nabla$-divergence of a differential form $\omega$ as
\bdm
\delta^{\nabla}\omega\ :=\ - \,\sum_{i=1}^n e_i \haken \nabla_{e_i}\omega\, .
\edm 
Then, for any exterior form $\omega$, the following formula holds:
\bdm
\delta^{\nabla} \omega \ = \ \delta^g \omega \, - \, \frac{1}{2} \cdot
\sum_{i,j=1}^n (e_i \haken e_j \haken T) \wedge (e_i \haken e_j \haken \omega)
\ .
\edm
In particular, for the torsion form itself, we obtain $\delta^{\nabla} T = 
\delta^g T=:\delta T$.
\end{prop}
\begin{cor} 
%----------
If the torsion form $T$ is $\nabla$-parallel, then its divergence vanishes, 
\bdm
\delta^g T \ = \ \delta^{\nabla} T \ = \ 0  .
\edm
\end{cor} 
We define the  divergence for a $(0,2)$-tensor $S$ as 
$\mathrm{div}^{\nabla}(S)(X)\ :=\  \sum_{i} (\nabla_{e_i}S)(X,e_i)$
and denote by $\mathrm{div}^g$ the divergence with respect to the
Levi-Civita connection $\nabla^g$. Then 
\bdm
\mathrm{div}^g(S)(X) \, -  \, \mathrm{div}^{\nabla}(S)(X) \ = \ 
- \, \frac{1}{2} \sum_{i,j=1}^n S(e_i \, , \, e_j) \, 
\T(e_i \, , \, X \, , \, e_j) \ = \ 0
\edm
because $S$ is symmetric while $T$ is antisymmetric, and we conclude 
immediately:
\begin{lem}[{\cite[Lemma 1.1]{Agri&F&N&P05}}]\label{same-div}
%------------------------------------------------------------
If $\nabla$ is a metric connection with totally skew-symmetric
torsion and $S$  a symmetric $2$-tensor, then 
\bdm
\mathrm{div}^g(S) \ = \ \mathrm{div}^{\nabla}(S)  .
\edm
\end{lem}
\begin{thm}[{\cite{Ivanov&P01}}]\label{curvature-identities}
%------------------------------------------------------------
The Riemannian curvature quantities and the $\nabla$-curvature
quantities are related by
\bea[*]
\kr^g(X,Y,Z,V) & =& \kr^\nabla (X,Y,Z,V) - \frac{1}{2}(\nabla_X T)(Y,Z,V)
+\frac{1}{2}(\nabla_Y T)(X,Z,V)\\
&& -\frac{1}{4}g(T(X,Y),T(Z,V)) -\frac{1}{4}\sigma_T(X,Y,Z,V)\\
\Ric^g(X,Y) &=& \Ric^\nabla(X,Y)  +\frac{1}{2}\delta T(X,Y)
-\frac{1}{4}\sum_{i=1}^{\dim M}g(T(e_i,X),T(e_i,Y))\\  
\Scal^\nabla &=& \Scal^g - \frac{3}{2} \, ||T||^2
\eea[*]
In particular, $\Ric^\nabla$ is symmetric if and only if $\delta T=0$,
\bdm
\Ric^\nabla(X,Y)-\Ric^\nabla(Y,X)\ =\ -\delta T (X,Y).
\edm
\end{thm}

%---------------------------------------------------------------------------
%
% Gesamtliteraturverzeichnis Uebersichtsartikel
%
% I. Agricola, WS 2005/06
%  
%---------------------------------------------------------------------------
%\chapter*{Bibliography}
%-----------------------------------------------------------------------
%\renewcommand{\bibname}{ }

%

%
\end{document}